\documentclass[12pt]{article}

\usepackage{amsfonts,amssymb,amsmath,latexsym,hyperref,fancyhdr}
\usepackage{graphicx}
\usepackage{epstopdf}
\usepackage{algorithmic}
\usepackage[british]{babel}
\usepackage{enumerate,bbm,xcolor,theorem,enumitem,verbatim,stmaryrd}
\usepackage[left=2cm,top=2cm,right=2cm,bottom=2cm]{geometry}

\numberwithin{equation}{section}
\pdfoutput=1

\ifpdf
  \DeclareGraphicsExtensions{.eps,.pdf,.png,.jpg}
\else
  \DeclareGraphicsExtensions{.eps}
\fi

\catcode`\<=\active \def<{
\fontencoding{T1}\selectfont\symbol{60}\fontencoding{\encodingdefault}}
\catcode`\>=\active \def>{
\fontencoding{T1}\selectfont\symbol{62}\fontencoding{\encodingdefault}}
\newcommand{\mathd}{\mathrm{d}}

\newcommand{\nobracket}{}
\newcommand{\noplus}{}

\newcommand{\tmem}[1]{{\em #1\/}}
\newcommand{\tmmathbf}[1]{\ensuremath{\boldsymbol{#1}}}
\newcommand{\tmop}[1]{\ensuremath{\operatorname{#1}}}
\newcommand{\tmstrong}[1]{\textbf{#1}}

\newenvironment{proof}{\noindent\textbf{Proof\ }}{\hspace*{\fill}$\Box$\medskip}

\newtheorem{lemma}{Lemma}[section]
\newtheorem{corollary}{Corollary}[section]

{\theorembodyfont{\rmfamily}\newtheorem{remark}{Remark}[section]}

\definecolor{myred}{RGB}{160,0,0}
\definecolor{mygreen}{RGB}{0,160,0}
\definecolor{myblue}{RGB}{0,0,160}
\hypersetup{
    colorlinks = true,
    linkcolor = {myred},
    citecolor = {mygreen},
    urlcolor  = {myblue},
}

\begingroup 
\makeatletter
\gdef\eqna@origamp{&} 
\catcode`\&\active 
\gdef\eqna@newamp{%
  \ifx\@currenvir\eqna@currenvir 
    \eqna@onlyfirstamp\let\eqna@onlyfirstamp\@empty 
  \else 
    \eqna@origamp 
  \fi
}
\gdef\eqna@hook{
  \let\eqna@currenvir\@currenvir 
  \catcode`\&\active 
  \let&\eqna@newamp 
  \let\eqna@onlyfirstamp\eqna@origamp 
  }
\catcode`\*11 
\gdef\eqnarray{\eqna@hook\align} 
\gdef\eqnarray*{\eqna@hook\align*} 
\global\let\endeqnarray\endalign
\global\let\endeqnarray*\endalign*
\endgroup

\title{A surprising observation on the quarter-plane diffraction problem}

\author{Raphael C. Assier$^{*}$ and I. David Abrahams$^{\dagger}$\\
\footnotesize{$^{*}$ Department of Mathematics, University of Manchester, Oxford Road, Manchester, {\rm M13 9PL}, UK}\\
\footnotesize{$^{\dagger}$ Isaac Newton Institute, University of Cambridge, 20 Clarkson Road, Cambridge CB3 0EH, UK}
}

\usepackage{amsopn}

\ifpdf
\hypersetup{
  pdftitle={A surprising observation},
  pdfauthor={R. C. Assier and I. D. Abrahams}
}
\fi

\newcommand{\RED}{} 

\newcommand{\COM}{} 

\definecolor{darkgreen}{rgb}{0,0.5,0}
\newcommand{\GRE}[1]{{\color{darkgreen}{#1}}} 

\begin{document}

\maketitle

\begin{abstract}
  {In this paper, we revisit Radlow's highly original attempt at a double Wiener--Hopf solution to the \COM{canonical} problem of wave diffraction by a quarter-plane. Using a constructive approach, we reduce the problem to two equations, one containing his somewhat controversial ansatz, and an additional compatibility equation. We then show that despite Radlow's ansatz being erroneous, it gives surprisingly accurate results in the far-field, in particular for the spherical diffraction coefficient. This unexpectedly good \COM{result} is established by comparing it to results obtained by the recently established modified Smyshlyaev formulae.}
\end{abstract}



\section{Introduction}\label{sec:intro}

Since the middle of the twentieth century, the intrinsically three-dimensional
canonical problem of wave diffraction by a quarter-plane has attracted a great deal of
attention, with many different mathematical techniques employed in seeking useful solutions.

This diffraction problem, a natural extension to Sommerfeld's famous
half-plane problem {\cite{Sommerfeldgerman,Sommerfeldtranslation}, represents one of the building blocks of the
	geometrical theory of diffraction (GTD, {\cite{GTD}}). Its far-field
	behaviour is very rich, including a set of primary and secondary edge-diffracted waves as well as a spherical wave emanating from the corner of the
	quarter-plane. The primary and secondary edge waves can be described
	analytically using the GTD, see for example {\cite{Assier2012b}}. Other
	techniques such as ray asymptotic theory on a surface of a sphere
	{\cite{Shanin2012}} or a Sommerfeld-Malyuzhinets integral approach {\cite{Lyalinov2013,Lyalinov2015}} also lead to the same results.
	However, the spherical wave is more problematic. In particular, one of the
	remaining challenges is to obtain a simple (easy to evaluate) closed-form
	expression for its diffraction coefficient.
	
	By considering the quarter-plane as a degenerated elliptic cone, the field can
	be expressed as a spherical wave multipole series involving Lam{\'e} functions {\cite{kraus,satterwhite,hansen}}. However these series
	are poorly convergent in the far-field and as such cannot lead to the
	sought-after diffraction coefficient. A review of this approach and attempts to
	accelerate the series convergence are described in {\cite{blume}}.
	
	A different and more recent way of considering this problem, based on the use
	of spherical Green's functions, has been introduced in 
	{\cite{smy1,Smyshlyaev1993,babich95}} and led to an integral formula for
	the spherical diffraction coefficient. However, this solution is not valid for
	all incidence/observation directions and requires a numerical treatment and
	some regularisation of Abel-Poisson type in order \COM{for it} to be evaluated {\cite{nonsingular}}.
	
	Building on this type of approach, a hybrid numerical-analytical method,
	which partially solves the acoustic quarter-plane problem in the Dirichlet
	case has been introduced in {\cite{shanin1,shanin2}}. The
	main advantage of this method compared to the one mentioned previously is that
	in this case the formulae giving the diffraction coefficient, known as the
	{\tmem{Modified Smyshlyaev Formulae}} (MSF), are `naturally convergent' in the
	sense that they do not require any special treatment to regularise or
	accelerate convergence. The method is based on planar and spherical edge
	Green's functions and on the theory of {\tmem{embedding formulae}}, introduced
	in {\cite{williams}} and further developed in {\cite{embedding}} for
	example. This method has been extensively described, adapted to the Neumann
	case and implemented in {\cite{Assier2012}}. We will use this method as a
	benchmark in the present paper; its implementation relies on an
	\textit{a priori} knowledge of the eigenvalues of the Laplace-Beltrami operator on a
	sphere with a slit. A detailed spectral analysis of this operator is given in
	{\cite{Assier2016}}. In particular, it gives a rapid way of evaluating the
	diffraction coefficient for a wide range of \COM{incident wave and} observer directions, but
	is not valid for all such directions. As discussed in {\cite{Assier2012b}}, a
	reason behind the limits of the MSF validity is the existence of secondary
	edge-diffracted waves.
	
	Another attempt, crucial to the present work, was published by Radlow in two
	successive papers {\cite{radlow,Radlow1965}}. The method is based
	on a Wiener--Hopf {\cite{Noble,Abrahams}} approach in two complex
	variables, and the author obtains a closed-form solution in Fourier space. 
	In the latter
	paper, an ansatz for the solution is proposed and a non-constructive intricate
	proof of its validity is given. This ansatz has long be\RED{en} known to be erroneous
	(see e.g. {\cite{meister}}), since it \COM{is shown to give} the wrong tip behaviour.
	The correct tip behaviour should include an eigenvalue of the Laplace-Beltrami
	operator (see {\cite{jansen}} for example). The technical reason as to why
	Radlow's proof is incorrect has been given fairly recently in {\cite{albani}}, in particular
	the field \COM{corresponding to} his ansatz does not satisfy the correct boundary
	condition. For a more extensive literature review on the use of functions of
	two complex variables in diffraction theory, \RED{the reader can be referred} to the introduction
	of {\cite{Assier2018a}}.
	
	In the present work, we revisit Radlow's approach and offer a formally exact solution from which we show that
	his ansatz appears constructively in a natural way. However, there is an extra term, which proves that Radlow's ansatz cannot be the
	true solution. The extra term is complicated, and contains integrals of as-yet unknown functions. The calculation/approximation of this term will be the subject of
	future work. However, while doing this work, we came across what we can refer
	to as a {\tmem{surprising observation}}. Serendipity made us compare the
	spherical diffraction coefficient calculated with Radlow's ansatz, i.e. setting th\COM{e} additional term to zero, to the one
	calculated using the MSF approach. It turns out, as we will show, that
	the two are \COM{very close} (at least in the Dirichlet case). Some hints can
	be found in the literature regarding the accuracy of Radlow's ansatz
	compared to full numerical computations {\cite{Albertsen1997,Tew1980}}, though, never before were the diffraction coefficients
	compared like for like.
	
	In Section \ref{sec:formulation}, the problem is formulated, and symmetries
	are exploited. In Section \ref{sec:Fourier}, the machinery required to work in
	Fourier space for two complex variables is introduced, the Wiener-Hopf
	functional equation is derived, and the inverse transform form of the solution
	is written down. Throughout this work, and starting from this section, we will
	use the phase portrait technique (see {\cite{Wegert2012}}) to visualise
	functions of a complex variable. This visualisation technique will play an
	important role in our reasoning. In Section \ref{sec:factK}, we present a way
	of factorising the Wiener-Hopf kernel into four factors with known analyticity
	properties. We write each factor as a modified Cauchy integral, in the form that allows easy implementation and fast evaluation. In Section
	\ref{sec:successiveWH}, two successive Wiener-Hopf procedures are performed,
	leading to the theoretical core of the present work: \COM{the two equations (\ref{eq:genericWHeq1})--(\ref{eq:genericWHeq2})} linking the
	main unknowns of the problem. The first equation involves Radlow's ansatz and
	an additional term, while the second equation, which we call the
	{\tmem{compatibility equation}}, may provide a way to find the unknown additional term. \COM{T}he diffraction coefficient
	is related to the solution of the Wiener-Hopf problem.
	Finally, in Section
	\ref{sec:results}, we compare the diffraction coefficient obtained by the MSF
	technique to that obtained assuming that Radlow's ansatz is correct. As we shall show, the two
	are surprisingly in very close agreement.
	
	\section{Formulation}\label{sec:formulation}
	
	\subsection{Geometry, governing equation and incident wave}
	
	Let us consider the three-dimensional $(x_1, x_2, x_3)$ space, and the quarter plane $\tmop{QP}$ defined by
	\begin{eqnarray}
	\tmop{QP} & = & \left\{ \tmmathbf{x}= (x_1, x_2, x_3) \in \mathbb{R}^3,
	\text{ such that } x_1 \geqslant
	0, x_2 \geqslant 0 \text{ and } x_3 = 0
	\right\}, 
	\end{eqnarray}
	and illustrated in Figure \ref{fig:qpplusspher}. We aim to solve the
	three-dimensional wave equation
	\begin{eqnarray}
	\frac{\partial^2 \mathfrak{u}_{\tmop{tot}}}{\partial t^2} = c^2 \Delta
	\mathfrak{u}_{\tmop{tot}} & \text{ and } & \frac{\partial^2
		\mathfrak{u}}{\partial t^2} = c^2 \Delta \mathfrak{u}, 
	\label{eq:timewaveeq}
	\end{eqnarray}
	in $\mathbb{R}^3 \backslash \tmop{QP}$ for the total velocity potential
	$\mathfrak{u}_{\tmop{tot}}$ and the scattered velocity potential $\mathfrak{u}$,
	when the quarter-plane is subject to an incident plane wave
	$\mathfrak{u}_{\tmop{in}} = e^{i (\tmmathbf{k} \cdot \tmmathbf{x}- \Omega t)}$,
	so that we can write $\mathfrak{u}_{\tmop{tot}} =\mathfrak{u}_{\tmop{in}}
	+\mathfrak{u}$. $\Omega$ represents the radian frequency of the incident wave, $c$ is
	the speed of sound and $\tmmathbf{k}$ is the incident wavevector, such that
	the wavenumber $k = | \tmmathbf{k} |$ is given by $k = \Omega / c$. To be consistent with Radlow, we take the total field to satisfy the Dirichlet (soft) boundary condition
	$\mathfrak{u}_{\tmop{tot}} = 0$ on $\tmop{QP}$. As is usual in scattering
	problems, let us make the hypothesis of time-harmonicity, assuming that all
	time-dependent quantities involved have a time-dependency consisting solely in
	a multiplicative factor $e^{- i \Omega t}$. We can then introduce the
	quantities $u_{\tmop{tot}} (\tmmathbf{x})$, $u_{\tmop{in}} (\tmmathbf{x})$ and
	$u (\tmmathbf{x})$, defined by $\mathfrak{u}_{\tmop{tot}} (\tmmathbf{x}, t) =
	\text{Re}(u_{\tmop{tot}} (\tmmathbf{x}) e^{- i \Omega t})$,
	$\mathfrak{u}_{\tmop{in}} (\tmmathbf{x}, t) = \text{Re}(u_{\tmop{in}} (\tmmathbf{x})
	e^{- i \Omega t})$ and $\mathfrak{u} (\tmmathbf{x}, t) = \text{Re}(u
	(\tmmathbf{x}) e^{- i \Omega t})$ respectively. As a consequence, the total
	field $u_{\tmop{tot}} (\tmmathbf{x})$ and the scattered field $u$ should
	satisfy the Helmholtz equation
	\begin{eqnarray}
	\Delta u + k^2 u = 0 & \text{ on } & \mathbb{R}^3 \backslash \tmop{QP} 
	\end{eqnarray}
	and $u_{\tmop{tot}}$ should satisfy the Dirichlet boundary condition
	\begin{eqnarray}
	u_{\tmop{tot}} = 0 & \text{ on } & \text{QP.}  \label{eq:Dirichletcondition}
	\end{eqnarray}
	The wavevector $\tmmathbf{k}$ is oriented in the incident direction towards
	the vertex of the quarter-plane (also the origin of our three-dimensional
	space) and as such, we can write $\tmmathbf{k}= - k\tmmathbf{\omega}_0$, where
	$\tmmathbf{\omega}_0$ represents the point of the unit sphere determining the
	incident direction. Using the spherical coordinates $(r, \theta, \varphi)$, as
	illustrated in Figure \ref{fig:qpplusspher}, we can introduce $\theta_0$ and
	$\varphi_0$, such that $\tmmathbf{\omega}_0$ corresponds to the point with
	spherical coordinates $(1, \theta_0, \varphi_0)$ and hence
	$\tmmathbf{\omega}_0$ can be represented in Cartesian coordinates by $(\sin
	(\theta_0) \cos (\varphi_0), \sin (\theta_0) \sin (\varphi_0), \cos
	(\theta_0))$.
	
	\begin{figure}[htbp]
		\centering
		\includegraphics[width=0.3\textwidth]{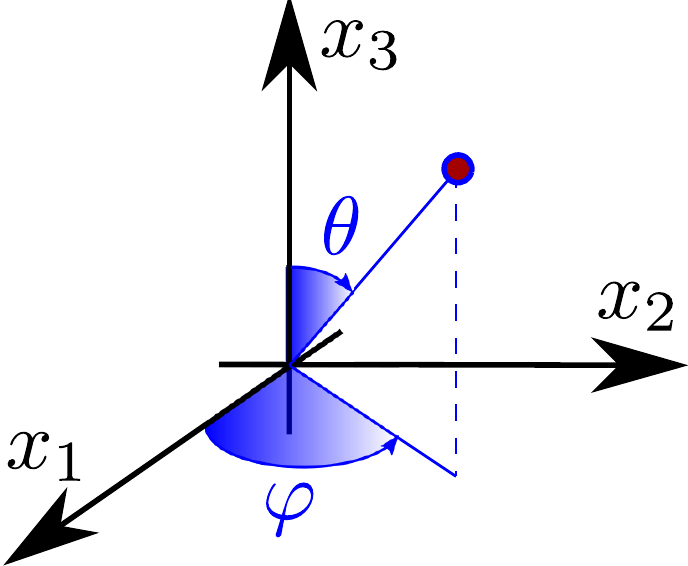}\quad\includegraphics[width=0.3\textwidth]{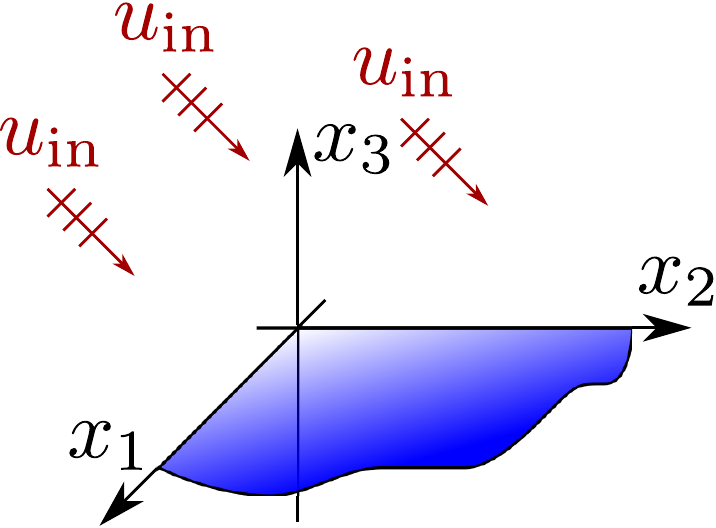}\quad  \includegraphics[width=0.3\textwidth]{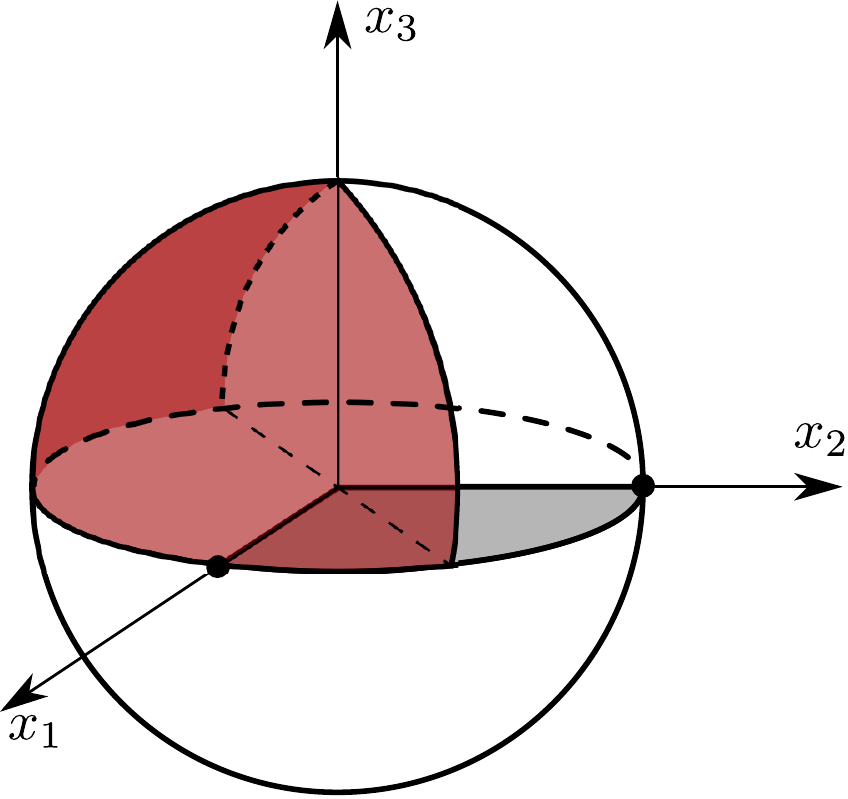}
		\caption{Spherical coordinates definition, quarter-plane illustration and geometric restriction of incidence}
		\label{fig:qpplusspher}
	\end{figure}
	
	The incident wave can hence be rewritten as
	\begin{eqnarray}
	& u_{\tmop{in}} (\tmmathbf{x}) = e^{i\tmmathbf{k} \cdot \tmmathbf{x}} =
	e^{- ik\tmmathbf{\omega}_0 \cdot \tmmathbf{x}} = e^{- i (a_1 x_1 + a_2 x_2 +
		a_3 x_3)}, &  \label{eq:incidentwave}
	\end{eqnarray}
	where $a_1 = k \sin (\theta_0) \cos (\varphi_0)$, $a_2 = k \sin (\theta_0) \sin
	(\varphi_0)$ and $a_3 = k \cos (\theta_0)$.
	
	\subsection{Edge, vertex and radiation conditions}\label{sec:edgevertexetc}
	
	In order for the problem to be well-posed, some other conditions need to be
	satisfied. These have been dealt with in \RED{detail} in \cite{Assier2016} for example, and so we will be brief. We impose the edge and vertex conditions: the energy of the field should remain bounded as we approach the edges and the vertex (i.e. no sources should be located on these), and the radiation condition: the scattered field $u$ should be outgoing in the far-field (i.e. no sources other than the incident wave at infinity).
	
	\subsection{Symmetries of the problem}\label{sub:symmetry}
	
	Let us now exploit the symmetry of the problem in order to reduce the range of
	the incident wave. First of all, due to the obvious ``vertical symmetry'' of
	the quarter-plane, it is enough to restrict the problem to incident waves
	coming from above the quarter-plane; this means that $\theta_0$ lies within $[0, \pi / 2]$. Moreover, in the $(x_1, x_2)$-plane,
	our domain is symmetric with respect to the bisector separating the
	quarter-plane into two plane sector with internal angle $\pi / 4$; i.e. it
	is possible to restrict $\varphi_0$ to belong to $[- 3 \pi / 4, \pi / 4]$,
	corresponding to the restricted zone of incidence depicted in Figure
	\ref{fig:qpplusspher}.
	

	\RED{Finally, it is well-known that the scattered field $u$ is symmetric (this can be seen by decomposing the field into its symmetric and antisymmetric parts)}, i.e., we have $u (x_1, x_2, x_3) = u (x_1, x_2, -
	x_3)$. Note that this automatically implies that $\partial u / \partial x_3$
	is an antisymmetric function. Therefore we can also restrict the observer
	region to $x_3 \geqslant 0$, i.e. $\theta \in [0, \pi / 2]$.

	
	
	
	\subsection{Jump in normal derivative across the quarter-plane}
	
	Let us consider the quantity
	\begin{eqnarray*}
		& f (x_1, x_2) = \left[ \frac{\partial u}{\partial x_3} \right]^{x_3 =
			0^+}_{x_3 = 0^-} = \frac{\partial u}{\partial x_3} (x_1, x_2, 0^+) -
		\frac{\partial u}{\partial x_3} (x_1, x_2, 0^-) . & 
	\end{eqnarray*}
	It is clear that in the part of the $x_3 = 0$ plane that does not contain
	$\tmop{QP}$, this quantity should be zero, since $u$ and its normal derivative
	are continuous. So we have that $f (x_1 < 0, x_2) = f (x_1, x_2 < 0) = 0$.  
	
	On QP, the far-field will be of the form $u = u_{\tmop{re}} + u_{\tmop{diff}}$
	on the (top) illuminated face, while it will be of the form $u = - u_{\tmop{in}} +
	u_{\tmop{diff}}$ on the bottom face. Here $u_{\tmop{re}}$ represents the
	reflected wave and is given by $u_{\tmop{re}} (x_1, x_2, x_3) = - e^{- i (a_1
		x_1 + a_2 x_2 - a_3 x_3)}$, and $u_{\tmop{diff}}$ encompasses all the
	different diffracted fields (primary and secondary edge diffraction plus corner
	diffraction), which decay at least like $1 / \sqrt{k \rho}$, where $\rho$ is
	the distance to the closest edge. Hence as both $x_1$ and $x_2$ tend to $+
	\infty$, we will have $u \sim u_{\tmop{re}}$ on the illuminated face and $u \sim -
	u_{\tmop{in}}$ on the bottom face. Hence we have
	\begin{eqnarray*}
		f (x_1, x_2) & \underset{x_1, x_2 \rightarrow + \infty}{\sim} &
		\frac{\partial u_{\tmop{re}}}{\partial x_3} (x_1, x_2, 0^+) + \frac{\partial
			u_{\tmop{in}}}{\partial x_3} (x_1, x_2, 0^-) \underset{x_1, x_2 \rightarrow + \infty}{=} \mathcal{O} (e^{- i
			(a_1 x_1 + a_2 x_2)}). 
	\end{eqnarray*}
	
	\subsection{Formulation summary}\label{sec:formulation-summary}
	
	To summarise, the scattering problem we wish to solve is the following:
	\begin{eqnarray}
	u_{\tmop{tot}} (\tmmathbf{x}) =  u_{\tmop{in}} (\tmmathbf{x}) + u
	(\tmmathbf{x}), & \text{ } & u_{\tmop{in}} (\tmmathbf{x})=e^{- i (a_1 x_1 + a_2 x_2 + a_3 x_3)}, \nonumber\\
	\Delta u + k^2 u = 0 \text{ on }\mathbb{R}^3 \backslash \text{QP}, &\text{ }& u_{\tmop{tot}} (\tmmathbf{x}) = 0 \text{ on QP}, \nonumber\\
	f (x_1, x_2) & \underset{x_{1, 2} \rightarrow \infty}{=} & \mathcal{O} (e^{-
		i (a_1 x_1 + a_2 x_2)})  \label{eq:finQ1}\\
	f (x_1, x_2) & = & 0 \begin{array}{l}
	\tmop{for} \hspace{1em} (x_1, x_2) \in Q_2 \cup Q_3 \cup Q_4
	\end{array},  \label{eq:finQ234}
	\end{eqnarray}
	subject to the vertex, edge and radiation conditions. The $Q_i$ are the
	different quadrants of the equatorial $(x_1, x_2)$-plane, illustrated in
	Figure \ref{fig:defdiaglog}, and defined by
	\begin{eqnarray*}
		Q_1 = \left\{ (x_1, x_2), x_1 \geqslant 0 \text{ and }
		x_2 \geqslant 0 \right\}, & \text{ } & Q_2 = \left\{ (x_1, x_2), x_1
		\leqslant 0 \text{ and } x_2 \geqslant 0 \right\},\\
		Q_3 = \left\{ (x_1, x_2), x_1 \leqslant 0 \text{ and }
		x_2 \leqslant 0 \right\}, &\text{ } & Q_4 = \left\{ (x_1, x_2), x_1
		\geqslant 0  \text{ and } x_2 \leqslant 0 \right\} .
	\end{eqnarray*}
	
	\COM{It is convenient at this point to defer the solution of this boundary value problem to Sections \ref{sec:successiveWH} and \ref{sec:results}. In the following Sections, \ref{sec:Fourier} and \ref{sec:factK}, it will be helpful to the general reader to first introduce the mathematical machinery to be employed later. This will, of necessity, be rather tedious; hence, a more experienced reader may choose to start with Section \ref{sec:successiveWH} and refer back to the earlier sections as required. }
	
	
	\section{Transformation in Fourier space}\label{sec:Fourier}
	
	\subsection{Some useful functions}\label{sec:usefulfunctions}
	
	In order to be able to define precisely quantities of interest in the
	following section, we need to introduce a few intermediate functions, as well
	as some useful notations. Let $\log (z)$ and $\sqrt{z}$ be the default complex
	logarithm and square root used by most mathematical software (e.g.
	Mathematica, Matlab, etc.). They correspond to the usual principal value of the logarithm and square root on the positive real axis and have a branch cut on the negative real axis.
	Let us now define a slightly different version of the logarithm: the function
	$\overset{\swarrow}{\log}$, that will be used first in \COM{S}ection
	\ref{sec:factK-oK+o}, defined by $\overset{\swarrow}{\log} (z) = \log \left( e^{- \frac{i \pi}{4}} z
	\right) + \frac{i \pi}{4}$,
	so that this is a logarithm in the sense that $\exp (
	\overset{\swarrow}{\log} (z)) = z$, it coincides with the usual real
	logarithm on the positive real axis, and has a branch cut extending diagonally down
	from the branch point $z = 0$, as illustrated in Figure \ref{fig:defdiaglog}.
	
	\begin{figure}[htbp]
		\centering
		\includegraphics[width=0.2\textwidth]{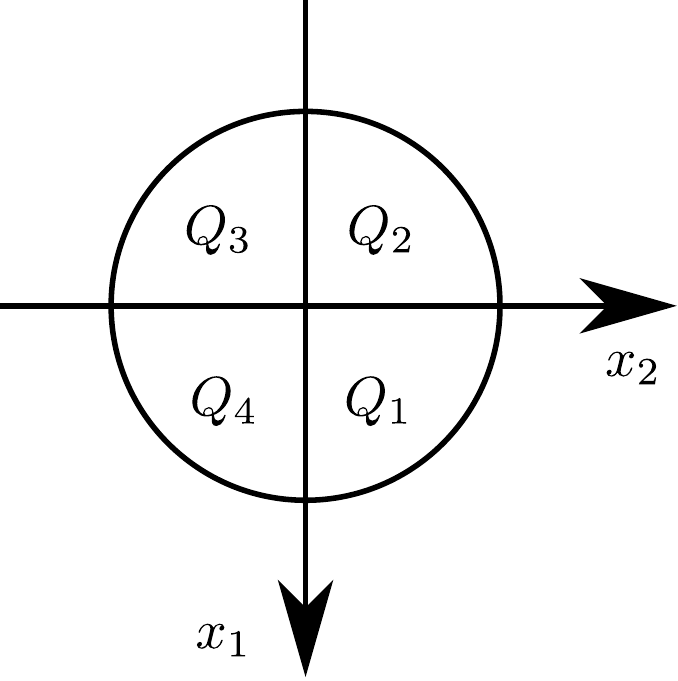} \quad \includegraphics[width=0.35\textwidth]{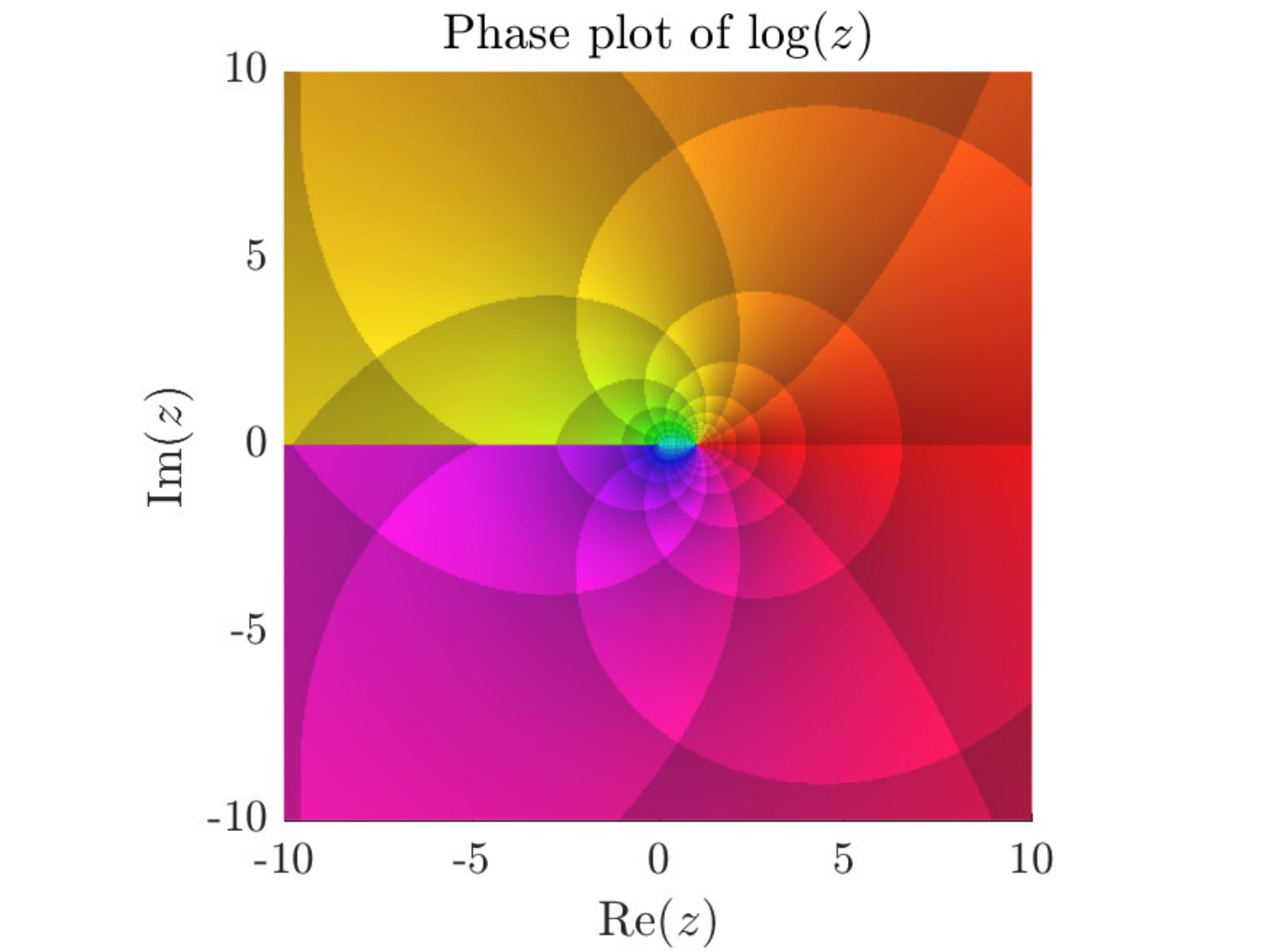} \quad \includegraphics[width=0.35\textwidth]{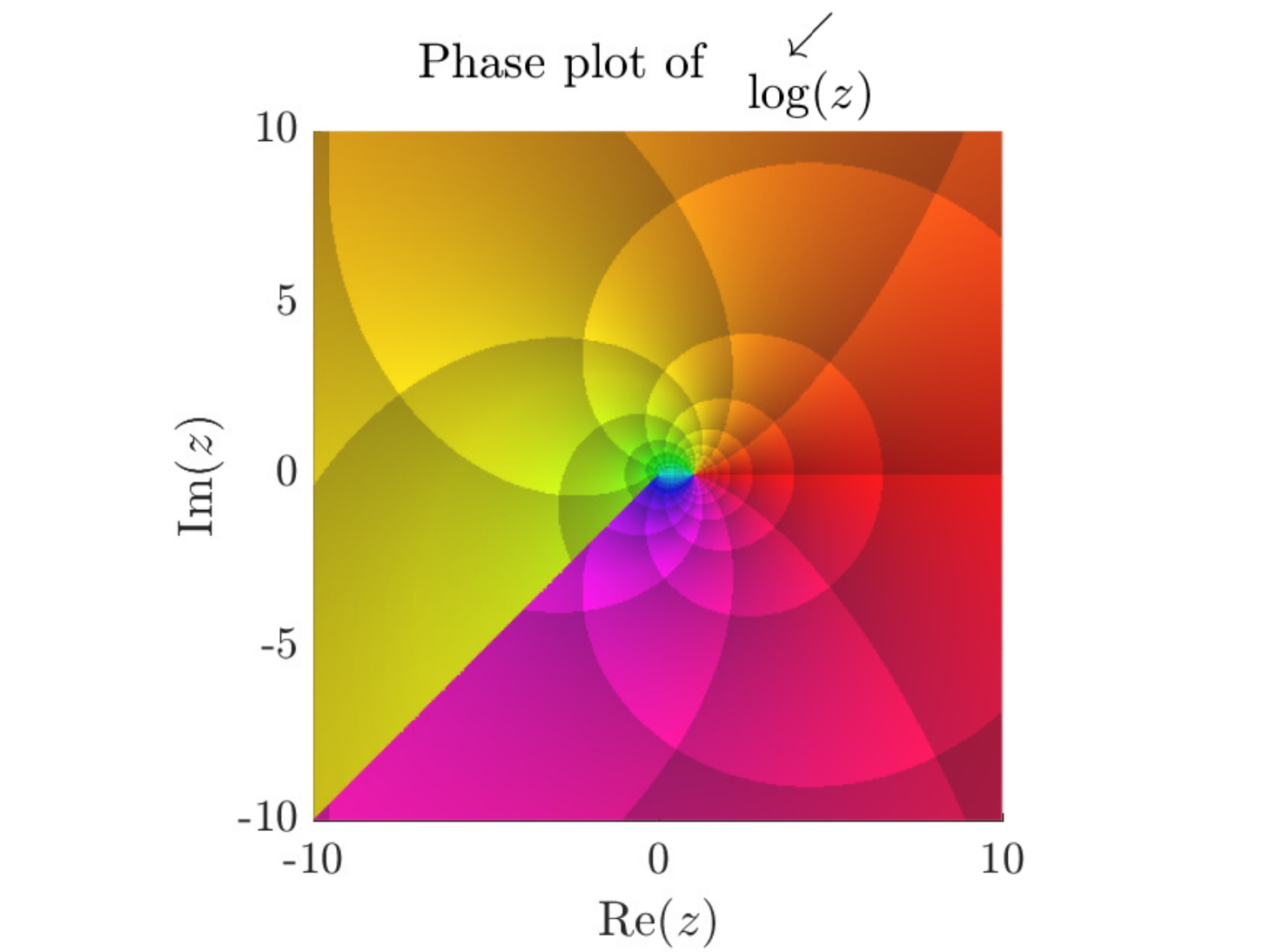}
		\caption{The quadrants $Q_i$ and phase portraits of the functions $\log (z)$ and
			$\overset{\swarrow}{\log} (z)$}
		\label{fig:defdiaglog}
	\end{figure}
	
	Let us now define the function $\sqrt[\downarrow]{z}$, that will be used
	extensively throughout this work, by $\sqrt[\downarrow]{z} = e^{i \frac{\pi}{4}} \sqrt{- iz}$ 
	so that this is a square root in the sense that $\left( \sqrt[\downarrow]{z}
	\right)^2 = z$, it coincides with the usual real square root on the positive
	real axis, and has a branch cut on the negative imaginary axis, as shown in
	Figure \ref{fig:differentrootsandkappa}. Building on this, we can define the
	function $\kappa (\mathfrak{K}, z)$ for any $\mathfrak{K}$ such that $\text{Im}
	(\mathfrak{K}) \geqslant 0$ and $\text{Re}(\mathfrak{K}) > 0$ by
	\begin{eqnarray}
	\kappa (\mathfrak{K}, z) & = & \sqrt[\downarrow]{\mathfrak{K}- z} 
	\sqrt[\downarrow]{\mathfrak{K}+ z}. \label{eq:defkapparev2}
	\end{eqnarray}
	The function $\kappa$ satisfies $(\kappa (\mathfrak{K}, z))^2 =\mathfrak{K}^2 - z^2$ with the principal Riemann sheet chosen such that $\kappa
	(\mathfrak{K}, 0) =\mathfrak{K}$. It has two branch cuts in the complex $z$ plane, one starting at the branch
	point $z =\mathfrak{K}$ and extending vertically upwards, and one starting at the
	branch point $z = -\mathfrak{K}$ and extending vertically downwards\footnote{The
		arrow notations used throughout this paper have the main objective of giving the
		reader an easy way to implement this work on a computer. One should also
		note that even if $\kappa$ is defined with two down-arrow functions, one of
		its branches is vertical upwards. This is due to the fact that the argument
		within one of the down-arrow function is $- z$.} as can be visualised in Figure
	\ref{fig:differentrootsandkappa}.
	
	\begin{figure}[htbp]
		\centering
		\includegraphics[width=0.3\textwidth]{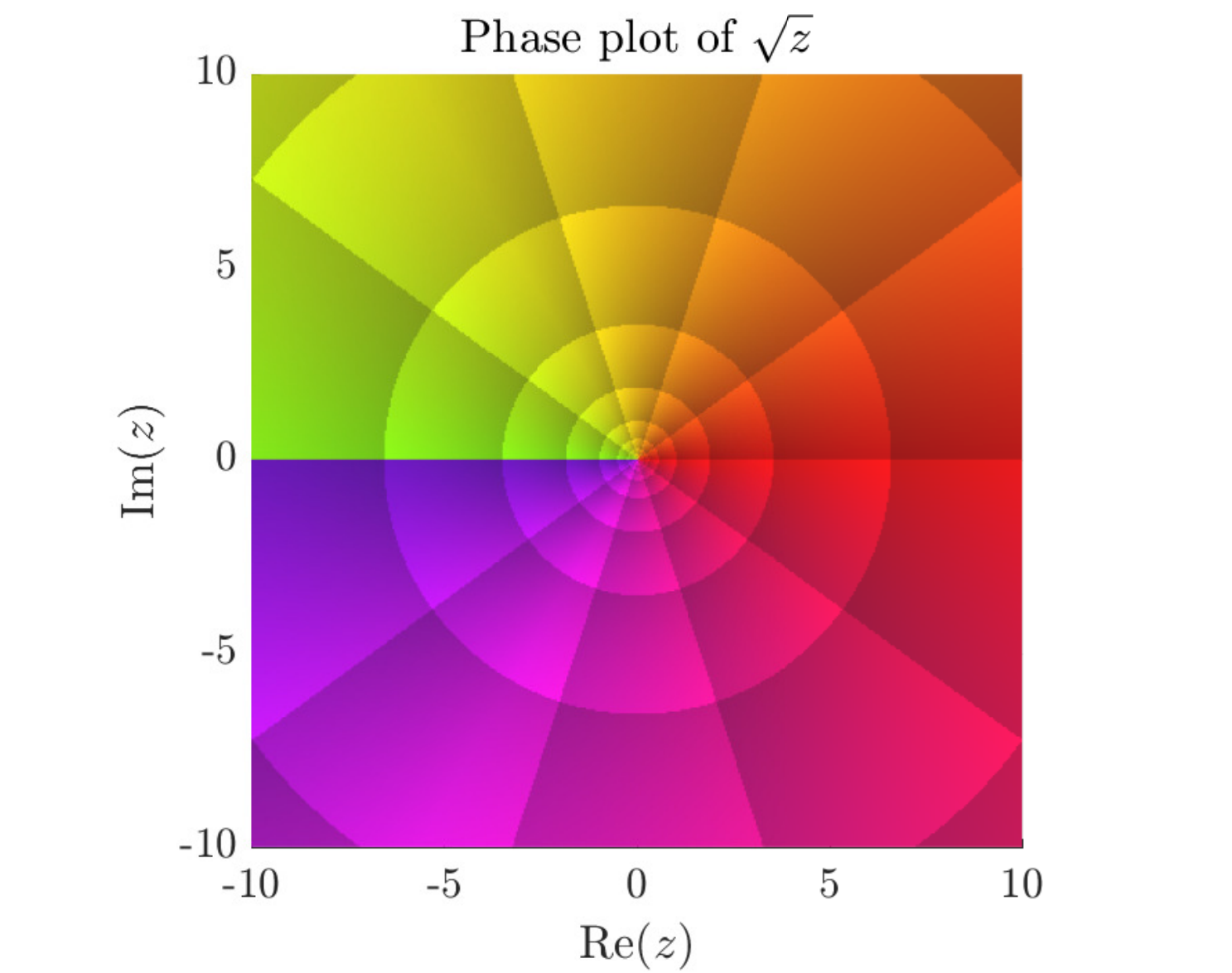}\includegraphics[width=0.3\textwidth]{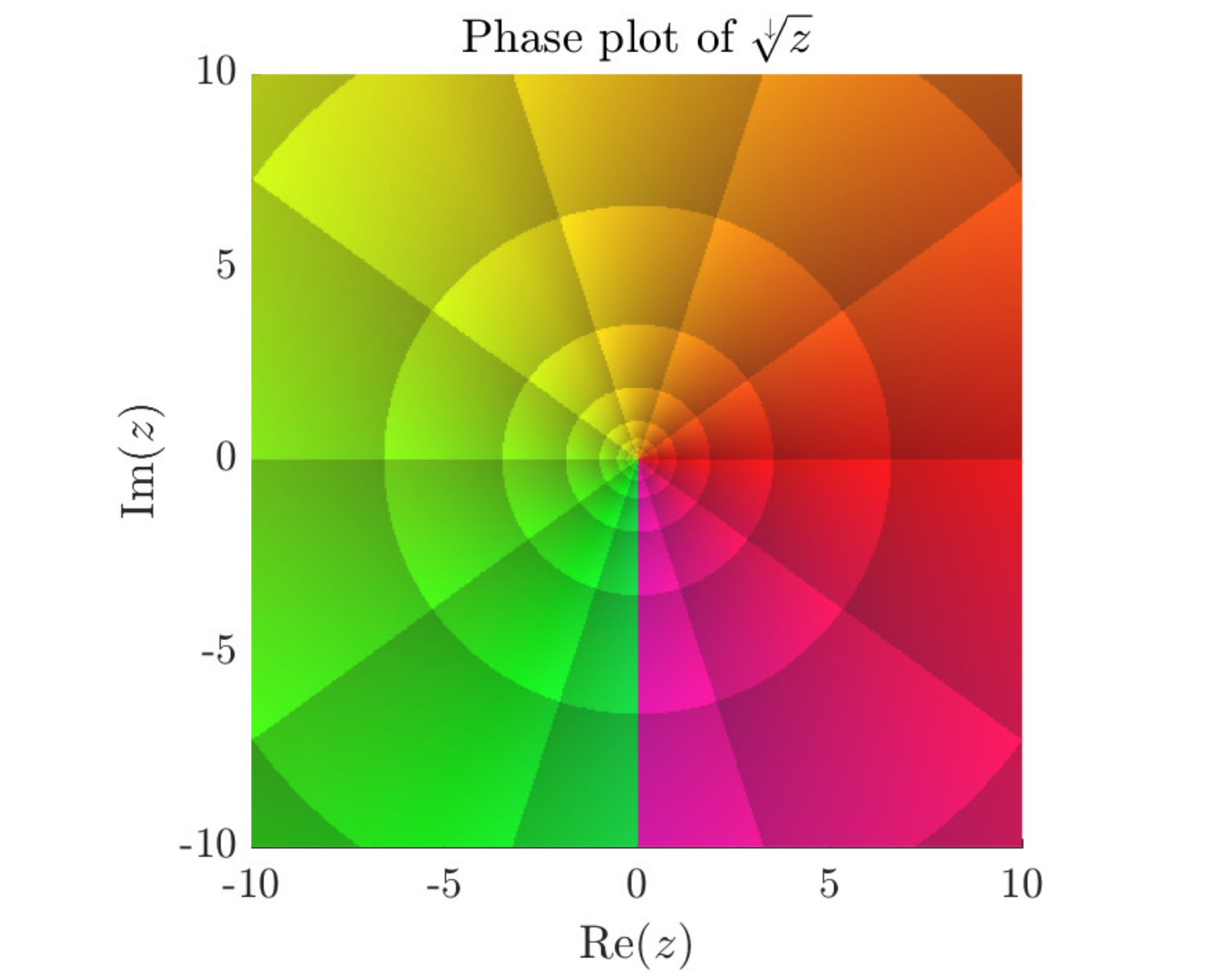}\includegraphics[width=0.3\textwidth]{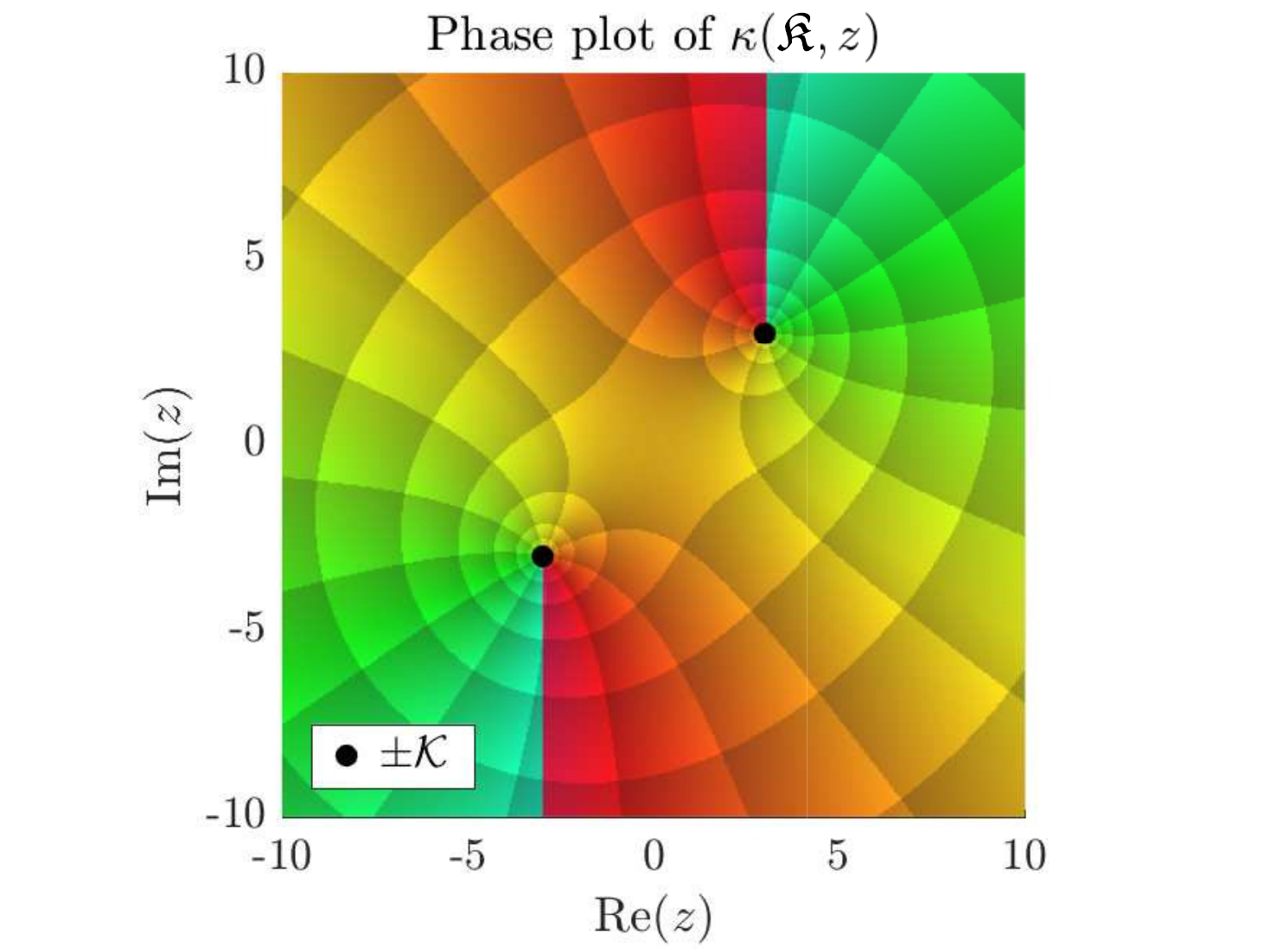}
		\caption{Phase portraits of the three functions $\sqrt{z}$,
			$\sqrt[\downarrow]{z}$, and $\kappa (\mathfrak{K}, z)$ for $\mathfrak{K}= 3
			+ 3 i$.}
		\label{fig:differentrootsandkappa}
	\end{figure}
	
	In the rest of the paper, we will sometimes use the bold notation
	$\tmmathbf{\alpha}$ to represent the two variables $(\alpha_1, \alpha_2)$. Let
	us now define the function $K (\tmmathbf{\alpha})$ as follows:
	\begin{eqnarray}
	K (\alpha_1, \alpha_2) & = & \frac{1}{\kappa (\kappa (k, \alpha_2),
		\alpha_1)},  \label{eq:defofK}
	\end{eqnarray}
	such that we have
	\begin{eqnarray*}
		& (K (\tmmathbf{\alpha}))^2 = \frac{1}{(\kappa (\kappa (k, \alpha_2),
			\alpha_1))^2} = \frac{1}{(\kappa (k, \alpha_2))^2 - \alpha_1^2} =
		\frac{1}{k^2 - \alpha_2^2 - \alpha_1^2} & 
	\end{eqnarray*}
	and define the function $\gamma (\tmmathbf{\alpha})$ as
	\begin{eqnarray}
	\gamma (\tmmathbf{\alpha}) = - i / K (\tmmathbf{\alpha}) & \text{ such that } & (\gamma (\tmmathbf{\alpha}))^2 = \alpha_1^2 +
	\alpha_2^2 - k^2 .  \label{eq:deffuncgamma}
	\end{eqnarray}
	Note that, by definition of $\kappa$, we have $1 / K (0, 0) = k$ and $\gamma (0, 0) = - ik$.
	
	\subsection{Double Fourier transform representation}\label{sec:doubleFourier}
	
	Let us now apply the dou\-ble Fourier transform (denoted by the operator
	$\mathfrak{F}$) in the ($x_1, x_2$) directions. Let us call $U (\alpha_1,
	\alpha_2, x_3)$ the double Fourier transform of $u (x_1, x_2, x_3)$, such that
	we have
	\begin{eqnarray*}
		U (\alpha_1, \alpha_2, x_3) &=\mathfrak{F} [u] = \int_{- \infty}^{\infty}
		\int_{- \infty}^{\infty} u (x_1, x_2, x_3) e^{i (\alpha_1 x_1 + \alpha_2
			x_2)} \, \mathd x_2 \mathd x_1, \\
		u (x_1, x_2, x_3) &=\mathfrak{F}^{- 1} [U] = \frac{1}{(2 \pi)^2}
		\int_{\mathcal{A}_1}^{} \int_{\mathcal{A}_2} U (\alpha_1, \alpha_2, x_3)
		e^{- i (\alpha_1 x_1 + \alpha_2 x_2)} \, \mathd \alpha_2 \mathd \alpha_1 .
	\end{eqnarray*}
	The contours of integration $\mathcal{A}_1$ and $\mathcal{A}_2$ in the inverse
	transform will  in general not completely lie on the real line, but will start at $-
	\infty$ and end at $+ \infty$. An exact description will be given in \COM{S}ection
	\ref{sub:givingA1A2}. Under this double Fourier transformation, the Helmholtz
	equation is changed into $(- \alpha_1^2 - \alpha_2^2) U + \frac{\partial^2 U}{\partial x_3^2} + k^2 U = 0$, which can be rewritten as
	\begin{eqnarray}
	\frac{\partial^2 U}{\partial x_3^2} - \gamma^2 (\tmmathbf{\alpha}) U = 0, &
	\text{ where, as already stated, } & \gamma^2 (\tmmathbf{\alpha}) = \alpha_1^2 + \alpha_2^2 - k^2.
	\label{eq:reducedhelmholtzfourier}
	\end{eqnarray}
	The contours $\mathcal{A}_1$ and $\mathcal{A}_2$
	will be chosen later such that $\text{Re}(\gamma
	(\tmmathbf{\alpha})) \geqslant 0$ when $\tmmathbf{\alpha} \in \mathcal{A}_1
	\times \mathcal{A}_2$. Hence in order not to have exponential growth as $x_3$
	tends to infinity, and because $x_3 \geqslant 0$, we must have
	\begin{eqnarray}
	U (\tmmathbf{\alpha}, x_3) & = & \COM{G} (\tmmathbf{\alpha}) e^{- \gamma
		(\tmmathbf{\alpha}) x_3}.  \label{eq:solutioninFourierspace}
	\end{eqnarray}
	Hence, we can write $u (\tmmathbf{x})$ using the inverse Fourier representation
	\begin{eqnarray}
	u (\tmmathbf{x}) & = & \frac{1}{(2 \pi)^2} \int_{\mathcal{A}_1}^{} \int_{\mathcal{A}_2}^{} \COM{G}
	(\alpha_1, \alpha_2) e^{- i (\alpha_1 x_1 + \alpha_2 x_2)} e^{- \gamma
		(\alpha_1, \alpha_2) x_3} \, \mathd \alpha_2 \mathd \alpha_1 . 
	\label{eq:ubeforeintroducingK}
	\end{eqnarray}
	We can write $f (x_1, x_2)$ in a similar fashion, using the symmetry of the
	solution (see Section \ref{sub:symmetry}):
	\begin{eqnarray}
	& f (x_1, x_2) = 2 \frac{\partial u}{\partial x_3} (x_1, x_2, 0^+) =
	\frac{-2}{(2 \pi)^2} \int_{\mathcal{A}_1}^{} \int_{\mathcal{A}_2}^{} \gamma
	(\tmmathbf{\alpha}) \COM{G} (\tmmathbf{\alpha}) e^{- i (\alpha_1 x_1 + \alpha_2
		x_2)} \mathd \alpha_2 \mathd \alpha_1 &  \label{eq:interFdebut}
	\end{eqnarray}
	Hence, upon introducing $F (\tmmathbf{\alpha})$ defined by
	\begin{eqnarray}
	F (\tmmathbf{\alpha}) & = & - 2 \gamma (\tmmathbf{\alpha}) \COM{G}(\tmmathbf{\alpha}),  \label{eq:introF}
	\end{eqnarray}
	the equation (\ref{eq:interFdebut}) becomes
	\begin{eqnarray*}
		f (x_1, x_2) & = & \frac{1}{(2 \pi)^2} \int_{\mathcal{A}_1}^{}
		\int_{\mathcal{A}_2}^{} F (\tmmathbf{\alpha}) e^{- i (\alpha_1 x_1 +
			\alpha_2 x_2)} \, \mathd \alpha_2 \mathd \alpha_1,
	\end{eqnarray*}
	which means that the function $F$ introduced in
	(\ref{eq:introF}) is in fact the double Fourier transform of $f$, i.e.
	\begin{eqnarray}
	F (\alpha_1, \alpha_2) & = & \int_{- \infty}^{\infty} \int_{-
		\infty}^{\infty} f (x_1, x_2) e^{i (\alpha_1 x_1 + \alpha_2 x_2)} \, \mathd
	x_2 \mathd x_1. \label{eq:fFT}
	\end{eqnarray}
	In what follows, it will be \COM{convenient} to use $K$
	instead of $\gamma$ \COM{and rewrite (\ref{eq:introF}) as
		\begin{align}
		G(\boldsymbol{\alpha})&=\tfrac{1}{2i}F(\boldsymbol{\alpha})K(\boldsymbol{\alpha}), \label{eq:revisionfunctionaleqraph}
		\end{align}
		which is the most important \textit{functional equation}\footnote{\COM{An alternative derivation, based on Green's identity, is
				given in \cite{Assier2018a}}.} of the problem. It relates the Fourier transform of $u$ and $\tfrac{\partial u}{\partial x_3}$ at $x_3=0^+$ and will be exploited to obtain the main result of the paper: equations (\ref{eq:genericWHeq1})--(\ref{eq:genericWHeq2}). Using (\ref{eq:revisionfunctionaleqraph}) in (\ref{eq:introF}) the wave field $u$ is given by }
	\begin{eqnarray}
	u (\tmmathbf{x}) & = & \frac{1}{(2 \pi)^2} \int_{\mathcal{A}_1}^{}
	\int_{\mathcal{A}_2}^{} \frac{F (\tmmathbf{\alpha}) K (\tmmathbf{\alpha})}{2
		i} e^{- i (\alpha_1 x_1 + \alpha_2 x_2)} e^{i \frac{x_3}{K
			(\tmmathbf{\alpha})}} \mathd \alpha_2 \mathd \alpha_1 .
	\label{eq:explicituFK}
	\end{eqnarray}
	
	\subsection{A small departure from the usual approach}
	
	As is usually the case when using the Wiener-Hopf technique, we could start by
	assuming that $k$ has a small positive imaginary part. Following this
	approach, it is possible to show that there \COM{exist} four real numbers $b_1,
	\delta_1, b_2$ and $\delta_2$, with $b_1 < \delta_1$ and $b_2 < \delta_2$,
	such that the function of interest $F (\tmmathbf{\alpha}) K
	(\tmmathbf{\alpha})$ is analytic on the tubular domain $\mathcal{D}^{\star}
	\subset \mathbb{C}^2$ defined by $\mathcal{D}^{\star} (b_1, b_2, \delta_1, \delta_2) = \mathcal{S} (b_1,
	\delta_1) \times \mathcal{S} (b_2, \delta_2)$, where for two real numbers $b < \delta$, the strip $\mathcal{S} (b, \delta)
	\subset \mathbb{C}$ is defined by $\mathcal{S} (b, \delta) = \{z \in \mathbb{C}, b < \text{Im} (z) <\delta \}$.
	In fact, it is possible to get an explicit expression for $\delta_{1, 2}$ and
	$b_{1, 2}$:
	\begin{eqnarray}
	\delta_1 = \text{Im} (k) | \cos (\varphi_0) |, &\text{ } \delta_2 = \text{Im} (k) | \sin
	(\varphi_0) |, &\text{ } b_{1, 2} = \max (- \delta_{1, 2}, \text{Im} (a_{1, 2})). 
	\end{eqnarray}
	However, if we want the solution for real $k$, the strips shrink to the real
	axes, and indented contours are needed in order to evaluate the inverse
	Fourier transforms. Our approach here, in the spirit of {\cite{Abrahams1997}},
	will be to start directly from such indented contours and avoid the limiting
	procedure discussion that would be required with the usual approach. We want
	to choose two contours $\mathcal{A}_1$ and $\mathcal{A}_2$ in the $\alpha_1$
	and $\alpha_2$ complex planes such that:
	\begin{itemize}[leftmargin=.25in]
		\item[(i)] For any $\alpha^{\star}_1 \in \mathcal{A}_1$, the functions $F
		(\alpha_1^{\star}, \cdot)$ and $K (\alpha_1^{\star}, \cdot)$ are analytic on
		$\mathcal{A}_2$.
		
		\item[(ii)] For any $\alpha^{\star}_2 \in \mathcal{A}_2$, the functions $F
		(\cdot^{}, \alpha_2^{\star})$ and $K (\cdot, \alpha_2^{\star})$ are analytic
		on $\mathcal{A}_1$.
		
		\item[(iii)] $\mathcal{A}_1$ and $\mathcal{A}_2$ are smooth contours starting at $-
		\infty$ and finishing at $+ \infty$. 
		
		\item[(iv)] For simplicity we would prefer that $\mathcal{A}_1$ be independent of
		$\alpha_2$ and $\mathcal{A}_2$ be independent of $\alpha_1$.
		
		\item[(v)] For any $\boldsymbol{\alpha} \in \mathcal{A}_1 \times \mathcal{A}_2$,
		$\text{Re}(\gamma(\boldsymbol{\alpha}))=\text{Im} (1 / K (\boldsymbol{\alpha})) \geqslant 0$.
	\end{itemize}
	\subsubsection{On fulfilling the requirements (i)--(v) for $K
		(\tmmathbf{\alpha})$}\label{sub:givingA1A2}
	
	In this subsection, we will show that there \RED{exist} contours $\mathcal{A}_1$
	and $\mathcal{A}_2$ that fulfil all the previous requirements (i)--(v)
	relative to the function $K (\tmmathbf{\alpha})$. Remember that $K
	(\tmmathbf{\alpha})$ is defined by $1 / \kappa (\kappa (k, \alpha_2),
	\alpha_1)$, and that by this definition (which breaks the symmetry between $\alpha_1$ and $\alpha_2$), $K$ does not behave in the same way
	in the $\alpha_1$ plane and in the $\alpha_2$ plane. In other words, even if
	by definition we have $K^2 (\alpha_1, \alpha_2) \equiv K^2 (\alpha_2,
	\alpha_1)$, we will not necessarily have $K (\alpha_1,
	\alpha_2) = K (\alpha_2, \alpha_1)$ for every $(\alpha_1, \alpha_2) \in
	\mathbb{C}^2$.
	
	To be more precise, for a fixed $\alpha_2^{\star}$ such that $\tmop{Im}
	(\kappa (k, \alpha_2^{\star})) \geqslant 0$, we expect the function $K
	(\alpha_1, \alpha_2^{\star})$ to simply have two branch points at $\pm \kappa
	(k, \alpha_2^{\star})$, with branch cuts extending vertically up and down, respectively, in the $\alpha_1$ complex plane, \COM{see Figure \ref{fig:proofofanalyticityK} (left)}. Hence, a suitable contour
	$\mathcal{A}_1$ would lie on the real line indented above $- \kappa (k, \alpha_2^{\star})$ and below
	$\kappa (k, \alpha_2^{\star})$ for any $\alpha_2^{\star} \in \mathcal{A}_2$.
	
	If we now fix an $\alpha_1^{\star}$ and consider the function $K
	(\alpha_1^{\star}, \alpha_2)$, we expect the analyticity structure to be a bit
	more complicated in the $\alpha_2$ plane. In particular, we expect to have
	potential problems at $\alpha_2 = \pm k$ due to the term $\kappa (k,
	\alpha_2)$, perhaps leading to a branch cut extending vertically upwards from $\pm
	k$. However, we also expect to have branch points where $\kappa (k, \alpha_2)
	= \pm \alpha_1^{\star}$, i.e., points where $\alpha_2 = \pm \kappa (k,
	\alpha_1^{\star})$, \COM{see Figure \ref{fig:proofofanalyticityK} (right)}. Hence, a suitable contour $\mathcal{A}_2$ would pass above
	$- k$ and $- \kappa (k, \alpha_1^{\star})$ and below $k$ and $\kappa (k,
	\alpha_1^{\star})$ for any $\alpha_1^{\star}\in\mathcal{A}_1$.
	
	If, as mentioned previously, it is possible to prove rigorously that some
	contours are valid in the case when $k$ has a small positive imaginary part,
	it is much harder to do so for real $k$. Instead, we will provide a
	{\tmem{visual proof}} that a given choice of $\mathcal{A}_1$ and
	$\mathcal{A}_2$ is suitable. Let us then consider the contours $\mathcal{A}_1$
	and $\mathcal{A}_2$ to be smoothly passing above $- k$ and below $k$ and also
	passing through the origins of their respective complex planes. A practical
	realisation of such contours can be obtained by the parametrisation $\mathcal{A}_1 (s_1) = s_1 + \frac{s_1}{a (s_1^4 + c)}$ and $\mathcal{A}_2 (s_2) = s_2 + \frac{s_2}{a (s_2^4 + c)}$, for $s_{1,2} \in \mathbb{R}$ and some complex constants $a$ and $c$. As such
	$\mathcal{A}_1$ and $\mathcal{A}_2$ satisfy (iii)-(iv).
	
	Given such a choice, it is possible to plot the loci of points $\pm \kappa
	(k, \mathcal{A}_2)$ in the $\alpha_1$ plane and the loci $\pm \kappa (k,
	\mathcal{A}_1)$ in the $\alpha_2$ plane. As long as our contours do not
	intersect these curves and do not intersect any resultant branch cuts, they should be
	valid. In fact, this can be seen in Figure \ref{fig:proofofanalyticityK},
	where the phase plots of $K (\alpha_1, \alpha_2^{\star})$ and $K
	(\alpha_1^{\star}, \alpha_2)$ are shown for different values of
	$\alpha_1^{\star} \in \mathcal{A}_1$ and $\alpha_2^{\star} \in \mathcal{A}_2$,
	together with the loci mentioned above.
	
	\begin{figure}[htbp]
		\centering
		\includegraphics[width=0.4\textwidth]{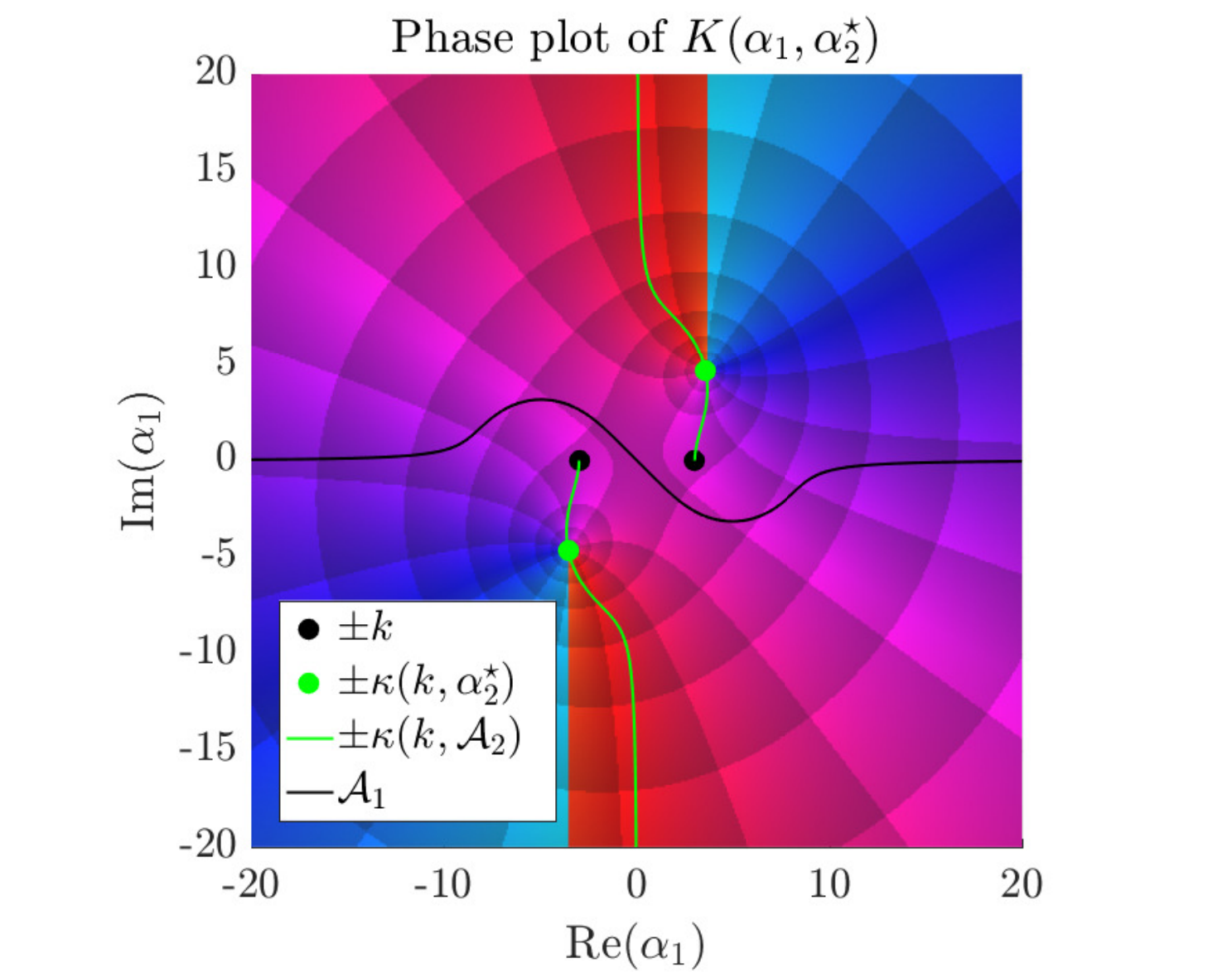}\includegraphics[width=0.4\textwidth]{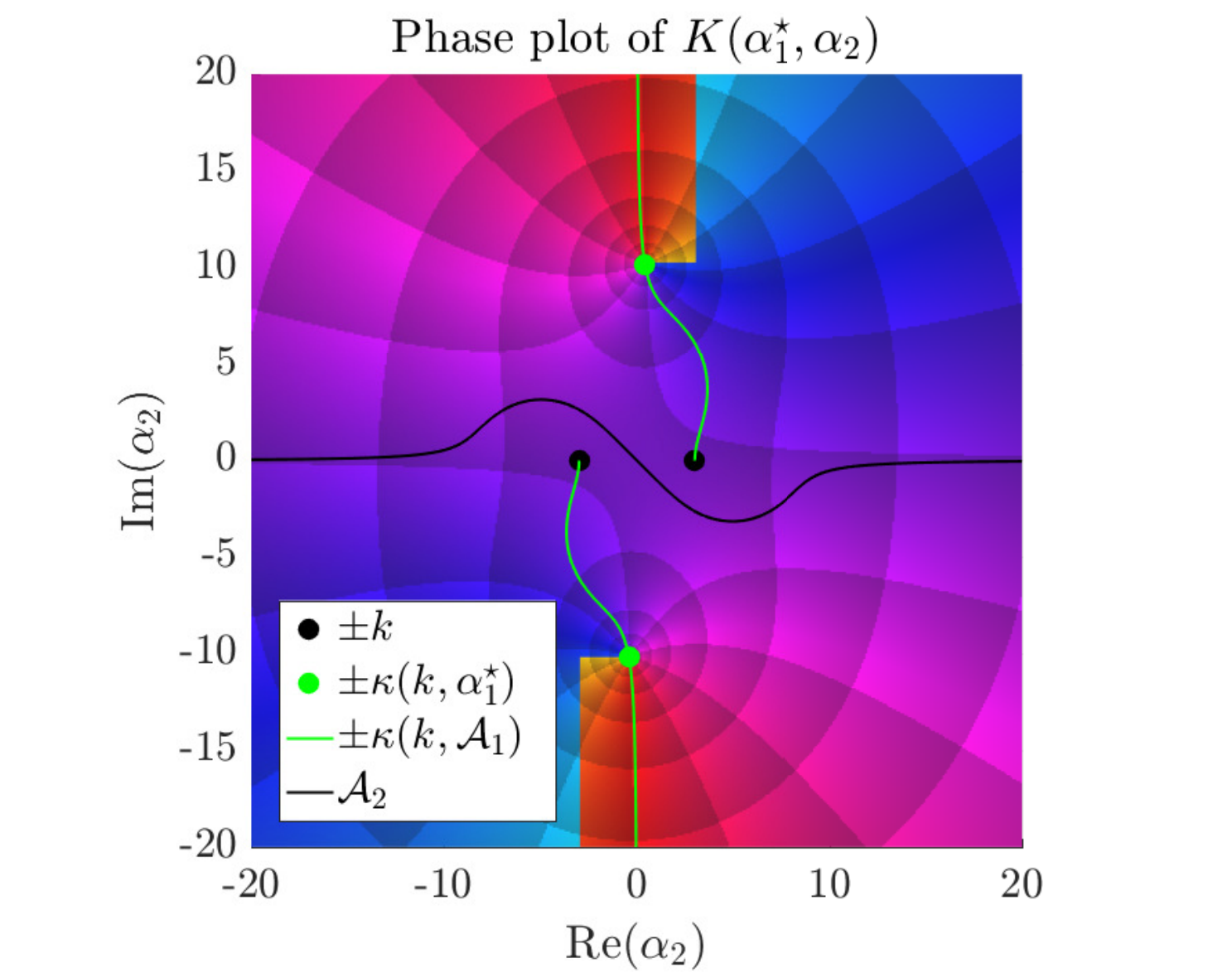}
		
		\includegraphics[width=0.4\textwidth]{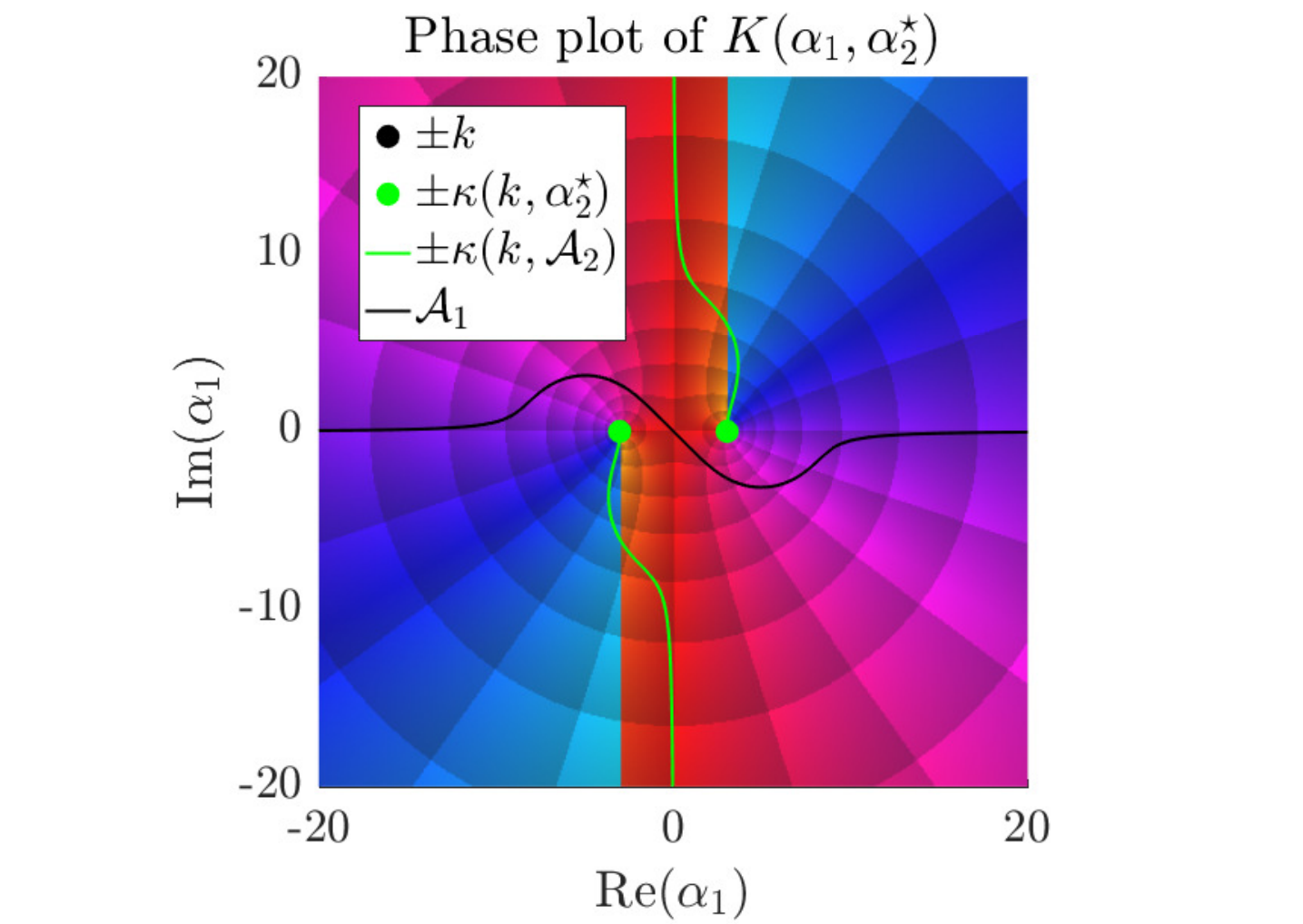}\includegraphics[width=0.4\textwidth]{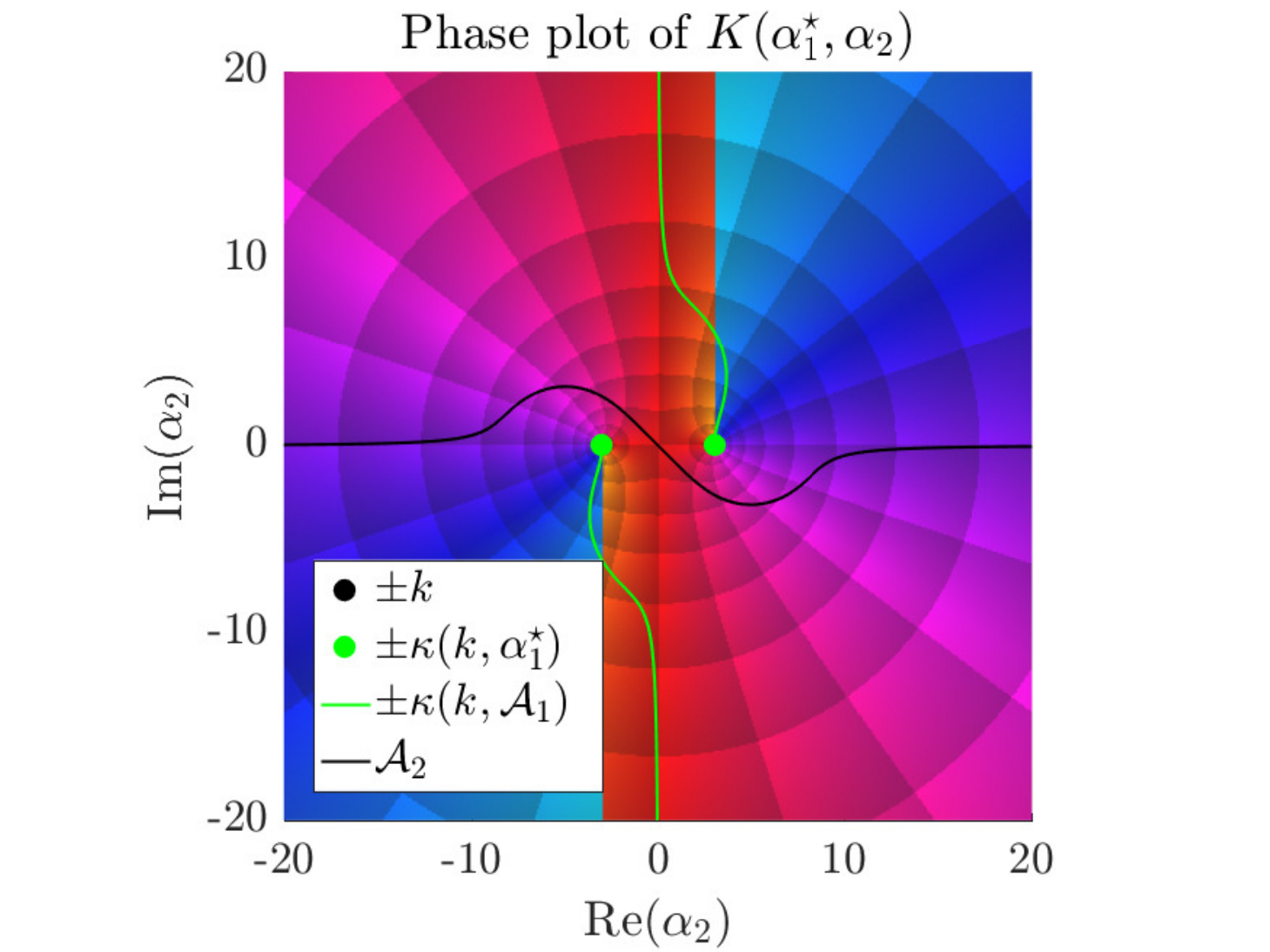}
		\caption{({\tmstrong{Visual proof of analyticity}}) Visualisation of $K$ in the $\alpha_1$ plane (left:
			$\alpha_2^{\star} =\mathcal{A}_2 (5)$(top) and $\alpha_2^{\star}
			=\mathcal{A}_2 (0)$ (bottom)) and in the $\alpha_2$ plane (right:
			$\alpha_2^{\star} =\mathcal{A}_1 (10)$(top) and $\alpha_1^{\star}
			=\mathcal{A}_1 (0)$ (bottom)). Here and in Figure \ref{fig:proofsigncompatibility} we chose $k=3$, $a=0.0012+0.0006 i$ and $c=1000i$.}
		\label{fig:proofofanalyticityK}
	\end{figure}
	
	As one can infer from Figure \ref{fig:proofofanalyticityK}, the contours
	$\mathcal{A}_1$ and $\mathcal{A}_2$, chosen suitably, avoid the singularities of $K$. In other words, for any $\alpha_2^{\star} \in \mathcal{A}_2$,
	the function $K (\alpha_1, \alpha_2^{\star})$ is analytic on $\mathcal{A}_1$,
	while for any $\alpha_1^{\star} \in \mathcal{A}_1$, the function $K
	(\alpha_1^{\star}, \alpha_2)$ is analytic on $\mathcal{A}_2$. Hence, as far as
	$K$ is concerned, this choice satisfies the conditions (i)--(iv). We still
	need to check that the condition (v) is satisfied. Again, here we will use a
	{\tmem{visual approach}}. The phase portrait of $\tmop{Im} (1 / K (\alpha_1,
	\alpha_2^{\star}))$ and $\tmop{Im} (1 / K (\alpha_1^{\star}, \alpha_2))$ for
	different values of $\alpha_1^{\star} \in \mathcal{A}_1$ and $\alpha_2^{\star}
	\in \mathcal{A}_2$ are displayed in Figure \ref{fig:proofsigncompatibility}.
	The regions where $\tmop{Im} (1 / K) > 0$ appear in red, while those where
	$\tmop{Im} (1 / K) < 0$ appear in blue.
	
	\begin{figure}[htbp]
		\centering
		\includegraphics[width=0.4\textwidth]{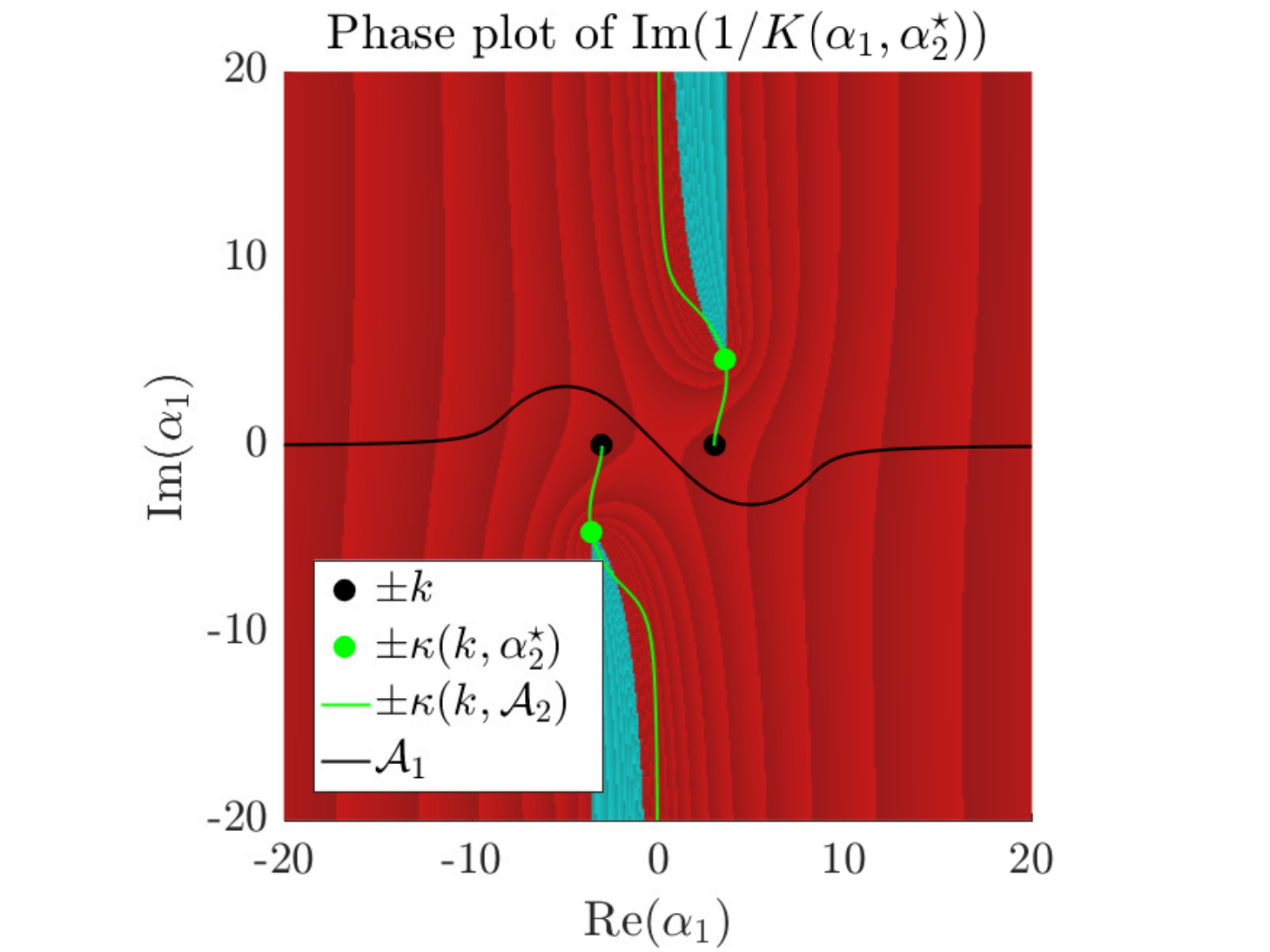}\includegraphics[width=0.4\textwidth]{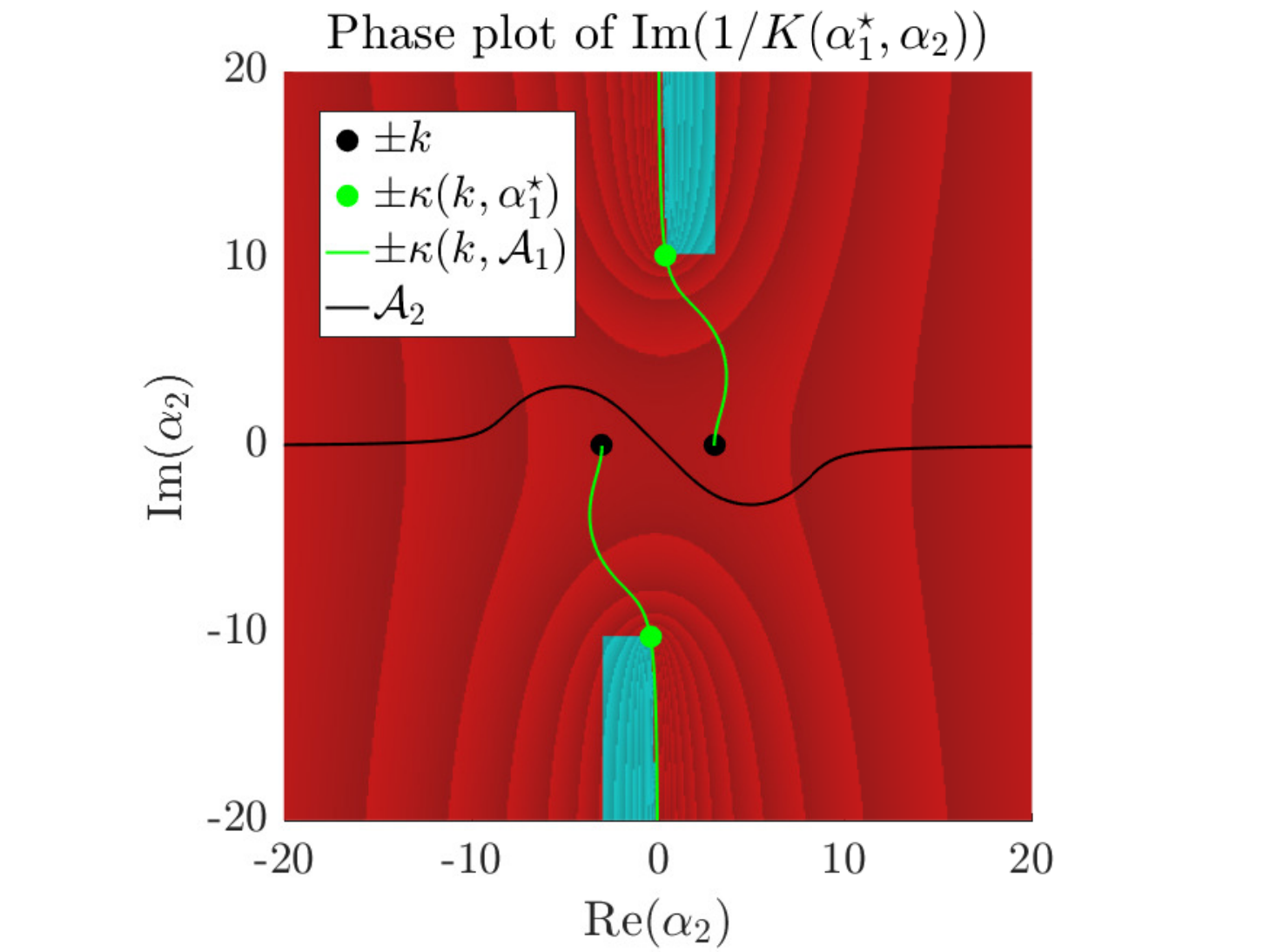}
		
		\includegraphics[width=0.4\textwidth]{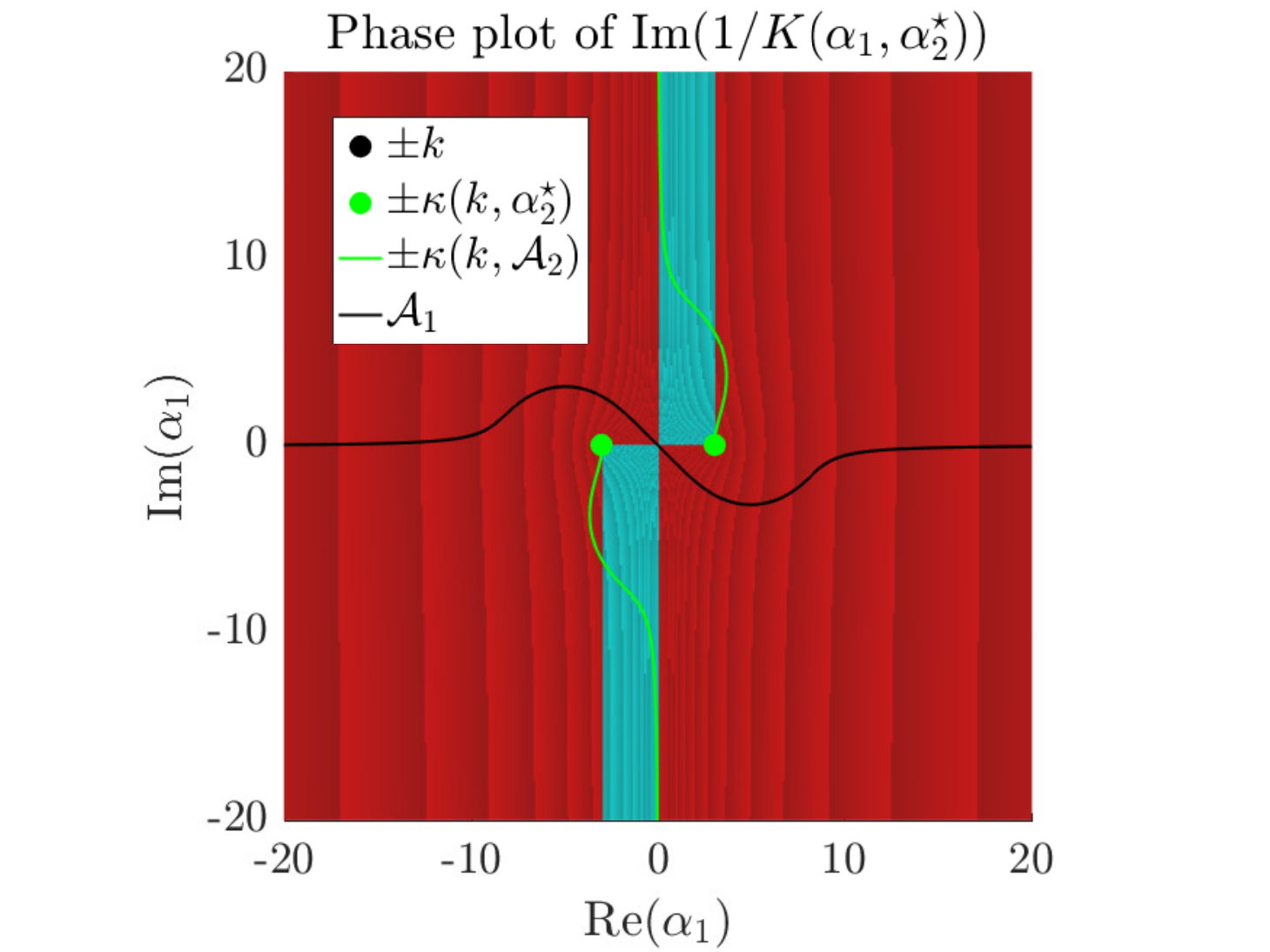}\includegraphics[width=0.4\textwidth]{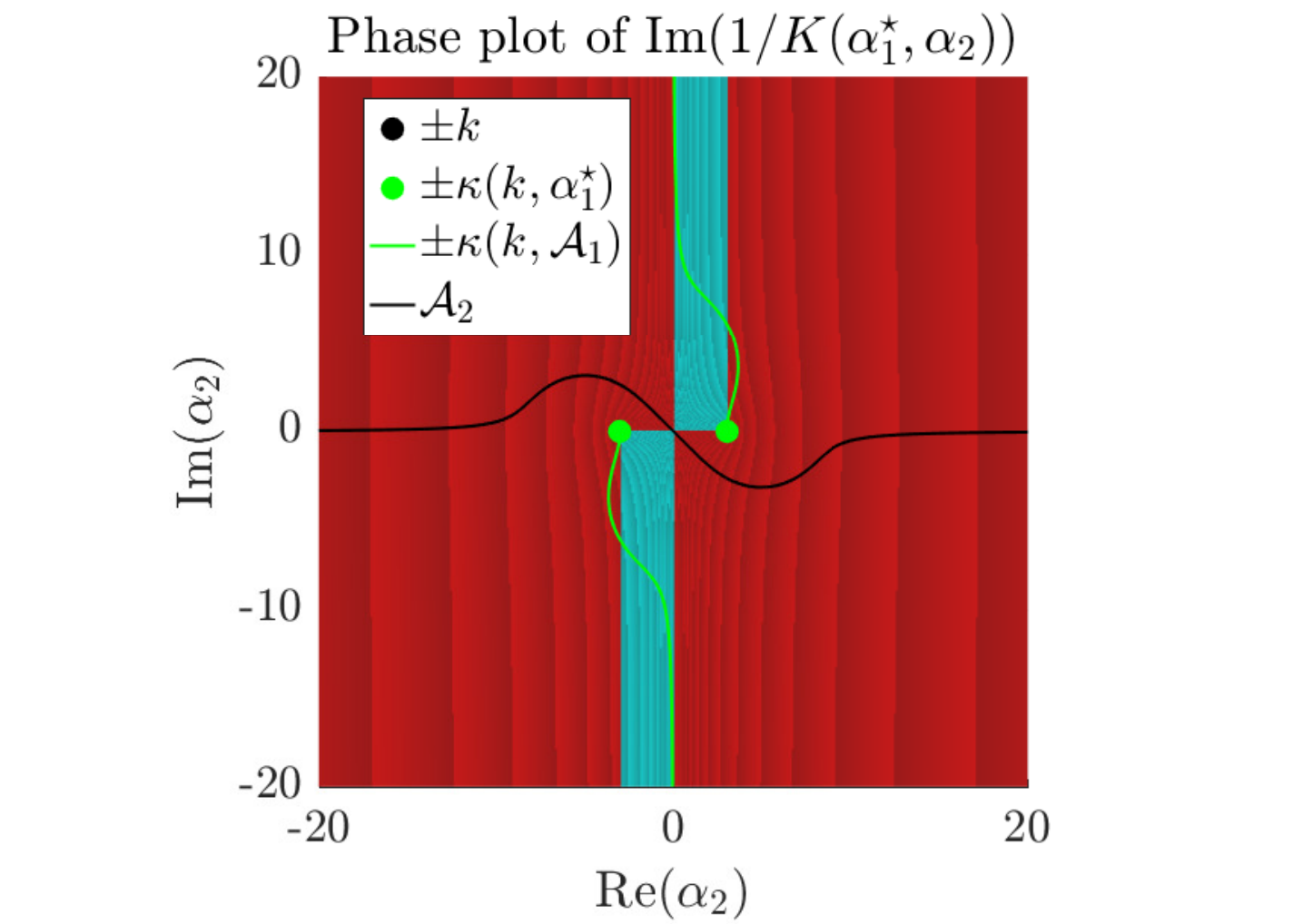}
		\caption{({\tmstrong{Visual proof of sign
					compatibility}}) Visualisation of $\tmop{Im} (1 / K)$ in the $\alpha_1$ plane (left:
			$\alpha_2^{\star} =\mathcal{A}_2 (5)$ (top) and $\alpha_2^{\star}
			=\mathcal{A}_2 (0)$ (bottom)) and in the $\alpha_2$ plane (right:
			$\alpha_1^{\star} =\mathcal{A}_1 (10)$ (top) and $\alpha_1^{\star}
			=\mathcal{A}_1 (0)$ (bottom)). The region where $\tmop{Im} (1 / K) \geqslant
			0$ appears in red on the plots. }
		\label{fig:proofsigncompatibility}
	\end{figure}
	
	As one can infer from Figure \ref{fig:proofsigncompatibility}, for any
	$\boldsymbol{\alpha} \in \mathcal{A}_1 \times \mathcal{A}_2$, we have
	$\tmop{Im} (1 / K (\boldsymbol{\alpha})) \geqslant 0$, as required in order for
	(v) to be satisfied. Note also that it only becomes zero when both $\alpha_1$
	and $\alpha_2$ are zero. It also shows that if $\mathcal{A}_1$ is chosen as
	above, $\mathcal{A}_2$ is forced to pass through the origin, and vice-versa.
	
	\subsubsection{On fulfilling the requirements (i)--(ii) for $F
		(\tmmathbf{\alpha})$}\label{sub:onFproperties}
	
	Remember that $F$ is defined in (\ref{eq:fFT}), and so using the condition
	(\ref{eq:finQ234}), it reduces to
	\begin{eqnarray}
	F (\alpha_1, \alpha_2) & = & \int_0^{\infty} \int_0^{\infty} f (x_1, x_2)
	e^{i (\alpha_1 x_1 + \alpha_2 x_2)} \, \mathd x_2 \mathd x_1 . \label{eq:defFwith0integrals}
	\end{eqnarray}
	In order to understand the analyticity property of $F$, we need to use the
	following lemma.
	
	\begin{lemma}
		\label{lem:quarter++}Let $\phi (x_1, x_2)$ be a function of the two real
		variables $x_1$ and $x_2$ and let $\gamma_1, \gamma_2 \in \mathbb{R}$ \RED{be} such
		that $| \phi (x_1, x_2) | \leqslant A_1 \exp (\gamma_1 x_1 + \gamma_2 x_2)$
		as $| x_1 | \rightarrow \infty$ and $| x_2 | \rightarrow \infty$ and $(x_1,
		x_2) \in Q_1$. Then the function $\Phi (\alpha_1, \alpha_2)$ defined by
		\begin{eqnarray*}
			\Phi (\alpha_1, \alpha_2) & = & \int_0^{\infty} \int_0^{\infty} \phi
			(x_1, x_2) e^{i (\alpha_1 x_1 + \alpha_2 x_2)} \, \mathd x_2 \mathd x_1
		\end{eqnarray*}
		can be interpreted as a function of the complex variable $\tmmathbf{\alpha}\in\mathbb{C}^2$, and as such, it is analytic in
		$\tmop{UHP} (\gamma_1) \times \tmop{UHP} (\gamma_2)$ considered an open
		subset of $\mathbb{C}^2$, where \COM{the upper-half plane} $\tmop{UHP} (\gamma_{1, 2})$ is the region
		in the $\alpha_{1, 2}$ complex plane lying above the horizontal line
		$\tmop{Im} (\alpha_{1, 2}) = \gamma_{1, 2}$. 
	\end{lemma}
	
	In our case, because of the estimate (\ref{eq:finQ1}), we can show that there
	exists $M > 0$, such that $| f (x_1, x_2) | \leqslant M \exp (\text{Im} (a_1) x_1 +
	\text{Im} (a_2) x_2)$ as $x_1, x_2 \rightarrow \infty$ within $Q_1$, where $a_{1,2}$ are related to the incident wave direction as defined below (\ref{eq:incidentwave}). Moreover, since
	$k$ is considered real, $\tmop{Im} (a_{1, 2}) = 0$. Hence, in the notation of
	Lemma \ref{lem:quarter++}, we have $\gamma_{1, 2} = 0$ and we can conclude
	that $F$ is analytic on $\tmop{UHP} (0) \times \tmop{UHP} (0)$, i.e, for
	$\tmop{Im} (\alpha_{1, 2}) > 0$.
	
	However, this does not mean that $F$ cannot be analytically continued onto a
	bigger domain. This realisation is important since the contours
	$\mathcal{A}_1$ and $\mathcal{A}_2$ defined in \RED{Section} \ref{sub:givingA1A2} do not lie
	within $\tmop{UHP} (0) \times \tmop{UHP} (0)$ since they both drop under their
	respective real axes.
	
	Hence, let us try to infer \textit{a priori}\footnote{Note that this particular aspect
		is studied more rigorously in {\cite{Assier2018a}}.} a bit more about the
	behaviour of $F$ outside $\tmop{UHP} (0) \times \tmop{UHP} (0)$. First of all,
	the estimate (\ref{eq:finQ1}), giving the behaviour of $f (x_1, x_2)$ at
	infinity gives us some information about the behaviour of $F
	(\tmmathbf{\alpha})$ within a finite part of the complex planes. Namely, we
	can expect that $F (\alpha_1, \alpha_2)$ will have a simple pole in the
	$\alpha_1$ plane at $\alpha_1 = a_1$ and a simple pole in the $\alpha_2$ plane
	at $\alpha_2 = a_2$. It also seems reasonable to expect that other possible
	singular behaviour\COM{s} would occur in the lower-half planes, e.g. branch points at $- k$ and maybe also on $- \kappa (k, \mathcal{A}_{1, 2})$ and at $- \kappa (k,
	a_{1, 2})$, once $\mathcal{A}_{1, 2}$ have been specified.
	
	Therefore, if $a_1$ and $a_2$ are negative, the contours $\mathcal{A}_1$ and
	$\mathcal{A}_2$ will be appropriate, since they are passing above the poles
	and the possible singular parts of $F$.
	
	\begin{remark}\label{note:introtildecontour}
		The situation is different if $a_{1,2}$ is positive, as then the contours $\mathcal{A}_{1,2}$ shown in Figure \ref{fig:proofofanalyticityK} will pass below the pole. A simple way to overcome \RED{what is} a technical difficulty is to allow $a_{1,2}$ to have a small imaginary part $\epsilon<0$ say, when $\text{Re}(a_{1,2})>0$. Then one can choose the contour $\mathcal{A}_{1,2}$ to lie sufficiently close to the real line that it passes above the pole, and the pole itself is located so that its residue will yield the correct behaviour for (\ref{eq:solutioninFourierspace}). Once the solution has been obtained, by continuity it should remain valid as $\epsilon\rightarrow0$.
	\end{remark}
	In what follows, in particular when drawing explanatory diagrams, unless
	stated otherwise, we will assume that $a_1$ and $a_2$ are both negative. We
	will make sure to provide accurate ways of dealing with the case $a_{1, 2} >
	0$ when necessary.
	
	\subsection{Set notations} \label{sec:setnotations}
	
	Let us start by introducing notations to describe useful sets in the
	$\alpha_1$ and $\alpha_2$ planes. We define the lower-half planes
	$\tmop{LHP}_1$ and $\tmop{LHP}_2$ and upper-half planes $\tmop{UHP}_1$ and
	$\tmop{UHP}_2$ as follows:
	\begin{eqnarray*}
		\tmop{LHP}_1 = \left\{ \alpha_1 \in \mathbb{C}, \text{ s.t. } 
		\alpha_1 \text{ lies \COM{below} } \mathcal{A}_1 \right\}, & \text{ } &
		\tmop{LHP}_2 = \left\{ \alpha_2 \in \mathbb{C}, \text{ s.t. } 
		\alpha_2 \text{ lies \COM{below} } \mathcal{A}_2 \right\},\\
		\tmop{UHP}_1 = \left\{ \alpha_1 \in \mathbb{C}, \text{ s.t. } 
		\alpha_1 \text{ lies \COM{above} } \mathcal{A}_1 \right\}, & \text{ } &
		\tmop{UHP}_2 = \left\{ \alpha_2 \in \mathbb{C}, \text{ s.t. } 
		\alpha_2 \text{ lies \COM{above} } \mathcal{A}_2 \right\}.
	\end{eqnarray*}
	Note that these sets are defined to be inclusive of the contour
	$\mathcal{A}_1$ and $\mathcal{A}_2$ in the sense that $\mathcal{A}_1 = \tmop{LHP}_1 \cap \tmop{UHP}_1$ and $ \mathcal{A}_2
	= \tmop{LHP}_2 \cap \tmop{UHP}_2$. The four types of sets introduced so far are illustrated in Figure \ref{fig:setsinC}.
	
	\begin{figure}[htbp]
		\centering
		\includegraphics[width=0.3\textwidth]{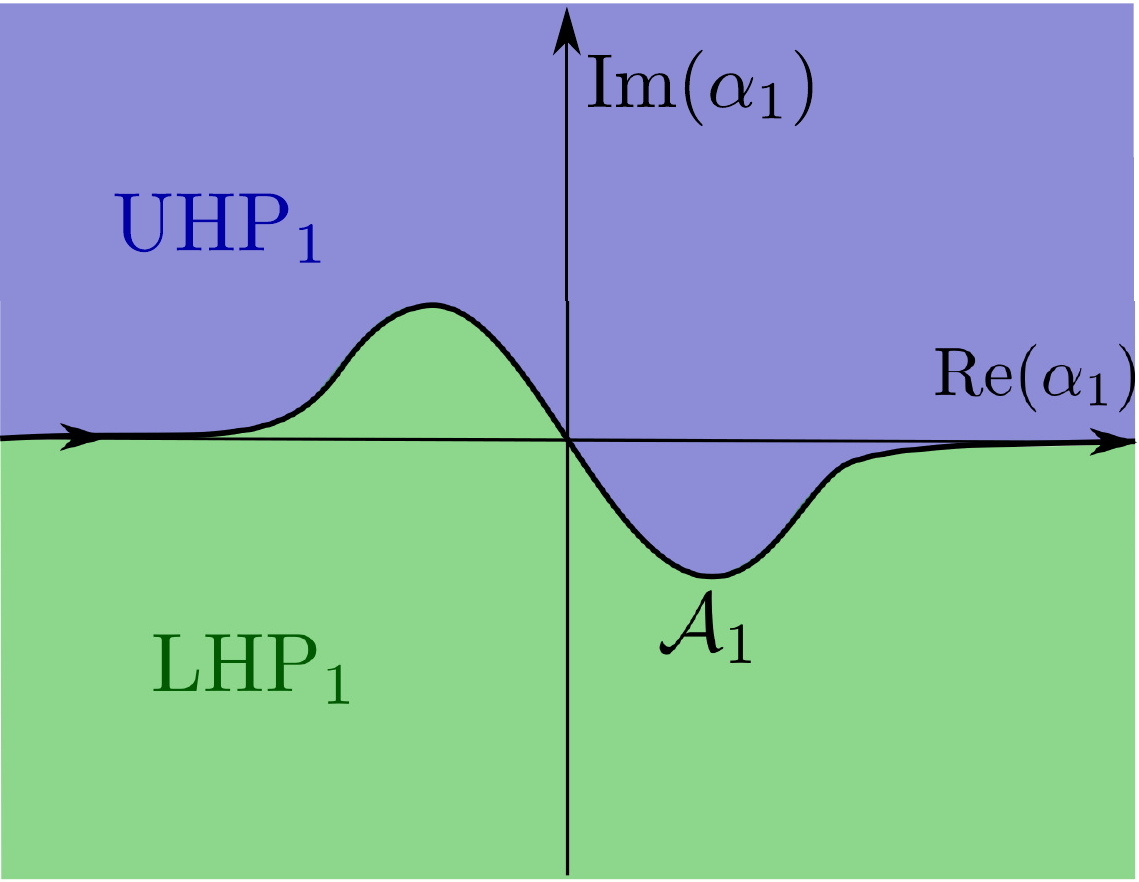}\qquad\includegraphics[width=0.3\textwidth]{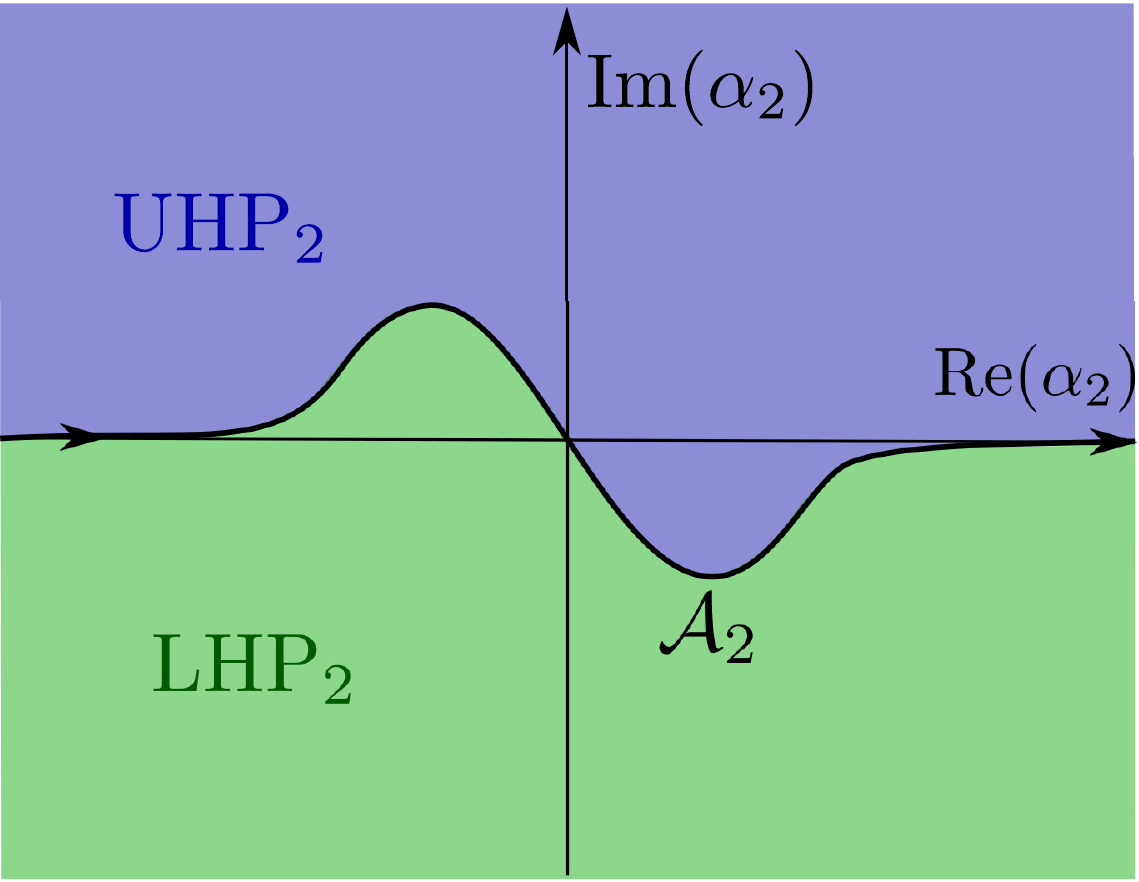}
		\caption{Diagrammatic description of the lower and upper-half planes used
			throughout this study.}
		\label{fig:setsinC}
	\end{figure}
	
	Let us now define a few different $\mathbb{C}^2$ sets derived from various
	products of the $\mathbb{C}$ spaces described above. We start with the set
	$\mathcal{D}=\mathcal{A}_1 \times \mathcal{A}_2$ where all of the functions we will deal with are well-behaved. It is also useful to define the $\mathbb{C}^2$ sets $\mathcal{D}_{+ +}= \tmop{UHP}_1 \times \tmop{UHP}_2$, $\mathcal{D}_{- +}= \tmop{LHP}_1 \times \tmop{UHP}_2$, $\mathcal{D}_{- -}= \tmop{LHP}_1 \times \tmop{LHP}_2$ and $\mathcal{D}_{+ -}= \tmop{UHP}_1 \times \tmop{LHP}_2$. Finally, let us introduce the sets $\mathcal{D}_{+ \circ}=\mathcal{D}_{+ +} \cap \mathcal{D}_{+ -} =  \tmop{UHP}_1 \times \mathcal{A}_2$ and $\mathcal{D}_{-\circ}=\mathcal{D}_{- -} \cap \mathcal{D}_{- +}= \tmop{LHP}_1 \times \mathcal{A}_2$.
	
	With the above points regarding analyticity now clarified, we can return to $F(\boldsymbol{\alpha})$ given in (\ref{eq:defFwith0integrals}) at the beginning of this subsection. It is clear that $F$ is analytic on $\mathcal{D}_{++}$ and hence we can rewrite it as
	\begin{eqnarray}
	F(\alpha_1,\alpha_2)=2 i F_{++}(\alpha_1,\alpha_2).
	\label{eq:defF++intermsofF}
	\end{eqnarray}
	
	\section{On the four-part factorisation of $K$}\label{sec:factK}
	
	Let us consider again the function $K (\tmmathbf{\alpha})$ defined by
	(\ref{eq:defofK}). We have shown in Section \ref{sub:givingA1A2} that $K
	(\tmmathbf{\alpha})$ is analytic on the product of contours
	$\mathcal{D}=\mathcal{A}_1 \times \mathcal{A}_2$. In this section, our aim is
	to show that there \COM{exist} four functions $K_{+ +} (\tmmathbf{\alpha})$, $K_{+
		-} (\tmmathbf{\alpha})$, $K_{- +} (\tmmathbf{\alpha})$ and $K_{- -}
	(\tmmathbf{\alpha})$, analytic on $\mathcal{D}_{+ +}$, $\mathcal{D}_{+ -}$,
	$\mathcal{D}_{- +}$ and $\mathcal{D}_{- -}$ respectively, such that for
	$\tmmathbf{\alpha} \in \mathcal{D}$, we have
	\begin{eqnarray*}
		K (\tmmathbf{\alpha}) & = & K_{+ +} (\tmmathbf{\alpha}) K_{+ -}
		(\tmmathbf{\alpha}) K_{- +} (\tmmathbf{\alpha}) K_{- -} (\tmmathbf{\alpha}).
	\end{eqnarray*}
	
	\subsection{Factorisation in the $\alpha_1$-plane} 
	
	\RED{Because of the definitions \COM{(\ref{eq:defkapparev2}) and (\ref{eq:defofK}) of $\kappa$ and $K$, we have}:
		\begin{align}
		K (\tmmathbf{\alpha}) & =  1 / \kappa (\kappa (k, \alpha_2), \alpha_1) = 1
		/ \left( \sqrt[\downarrow]{\kappa (k, \alpha_2) - \alpha_1}
		\sqrt[\downarrow]{\kappa (k, \alpha_2) + \alpha_1} \right),
		\label{eqrev:defofKandkapa}
		\end{align}
		\COM{and} one can see that for any $\tmmathbf{\alpha} \in \mathcal{D}$, it is possible to write
		\begin{align*}
		K (\tmmathbf{\alpha}) & =  K_{- \circ} (\tmmathbf{\alpha}) K_{+ \circ}
		(\tmmathbf{\alpha}),
		\end{align*}
		such that for a given $\alpha_2 \in \mathcal{A}_2$, $K_{- \circ} (\alpha_1,
		\alpha_2)$ is analytic (as a function of $\alpha_1$) in $\tmop{LHP}_1$ and
		$K_{+ \circ} (\alpha_1, \alpha_2)$ is analytic (as a function of $\alpha_1$)
		in $\tmop{UHP}_1$. 
		Exact expressions for $K_{- \circ}$ and $K_{+ \circ}$ follow from (\ref{eqrev:defofKandkapa}):}
	\begin{eqnarray}
	K_{- \circ} (\tmmathbf{\alpha}) = 1 / \sqrt[\downarrow]{\kappa (k, \alpha_2)
		- \alpha_1} & \text{ and } & K_{+ \circ} (\tmmathbf{\alpha}) = 1 /
	\sqrt[\downarrow]{\kappa (k, \alpha_2) + \alpha_1}.
	\label{eq:explicitKminusnote}
	\end{eqnarray}
	Indeed, for a given $\alpha_2 \in \mathcal{A}_2$, the only branch point of
	$K_{- \circ} (\tmmathbf{\alpha})$ is at $\alpha_1 = \kappa (k, \alpha_2)$,
	which is strictly within $\tmop{UHP}_1$ so that $K_{- \circ}
	(\tmmathbf{\alpha})$ is a {\tmem{minus function}} when considered as a
	function of $\alpha_1$, i.e., it is analytic in $\tmop{LHP}_1$. Similarly, the
	only branch point of $K_{+ \circ} (\tmmathbf{\alpha})$ is at $\alpha_1 = -
	\kappa (k, \alpha_2)$, which is strictly within $\tmop{LHP}_1$ so that $K_{+
		\circ} (\tmmathbf{\alpha})$ is a {\tmem{plus function}} when considered as a
	function of $\alpha_1$, i.e. it is analytic in $\tmop{UHP}_1$. This
	factorisation is illustrated in Figure \ref{fig:firstKfactorisation}. 
	
	\begin{figure}[htbp]
		\centering
		\raisebox{-.5\height}{\includegraphics[width=0.27\textwidth]{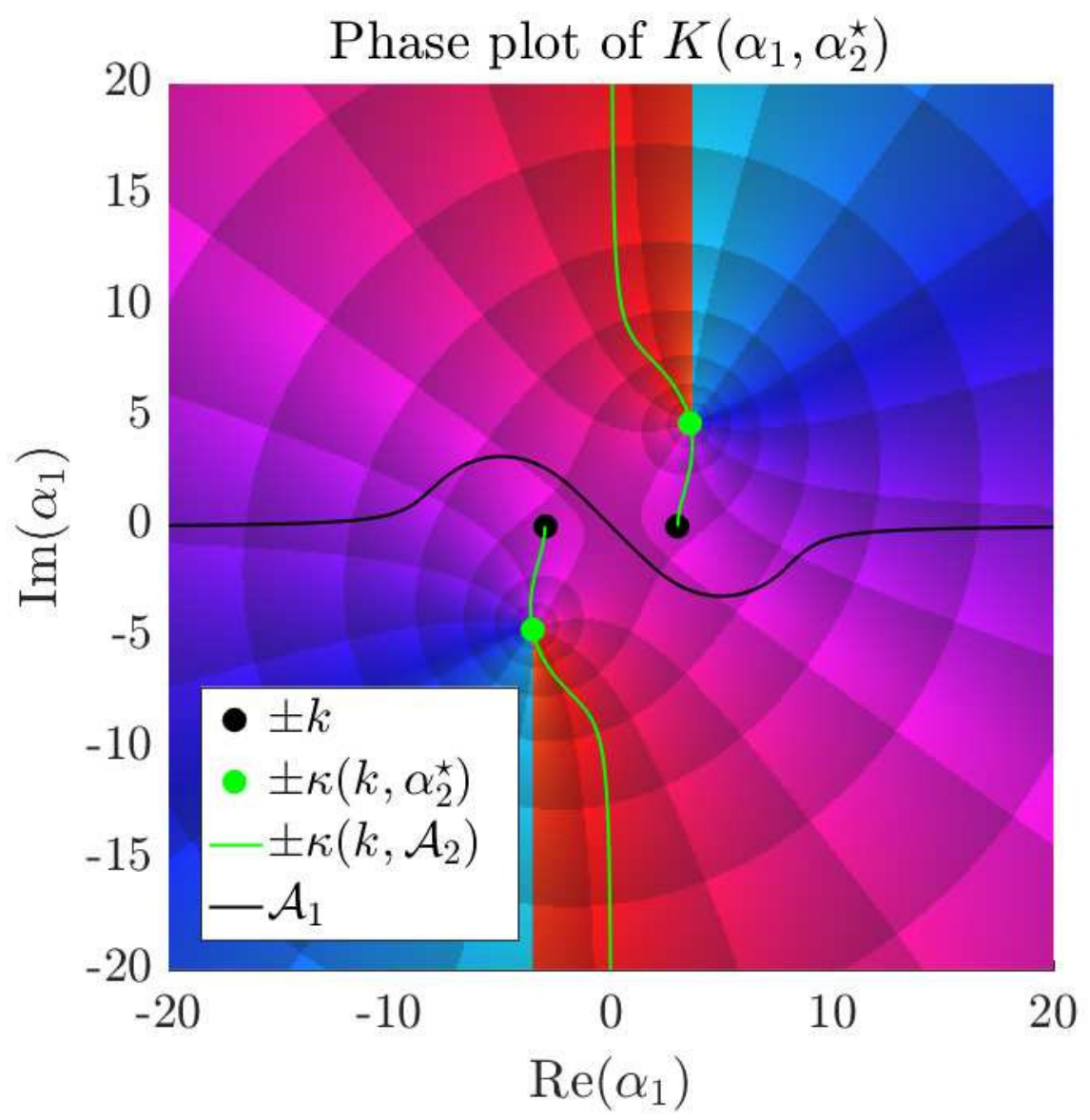}}$=$ \raisebox{-.5\height}{\includegraphics[width=0.27\textwidth]{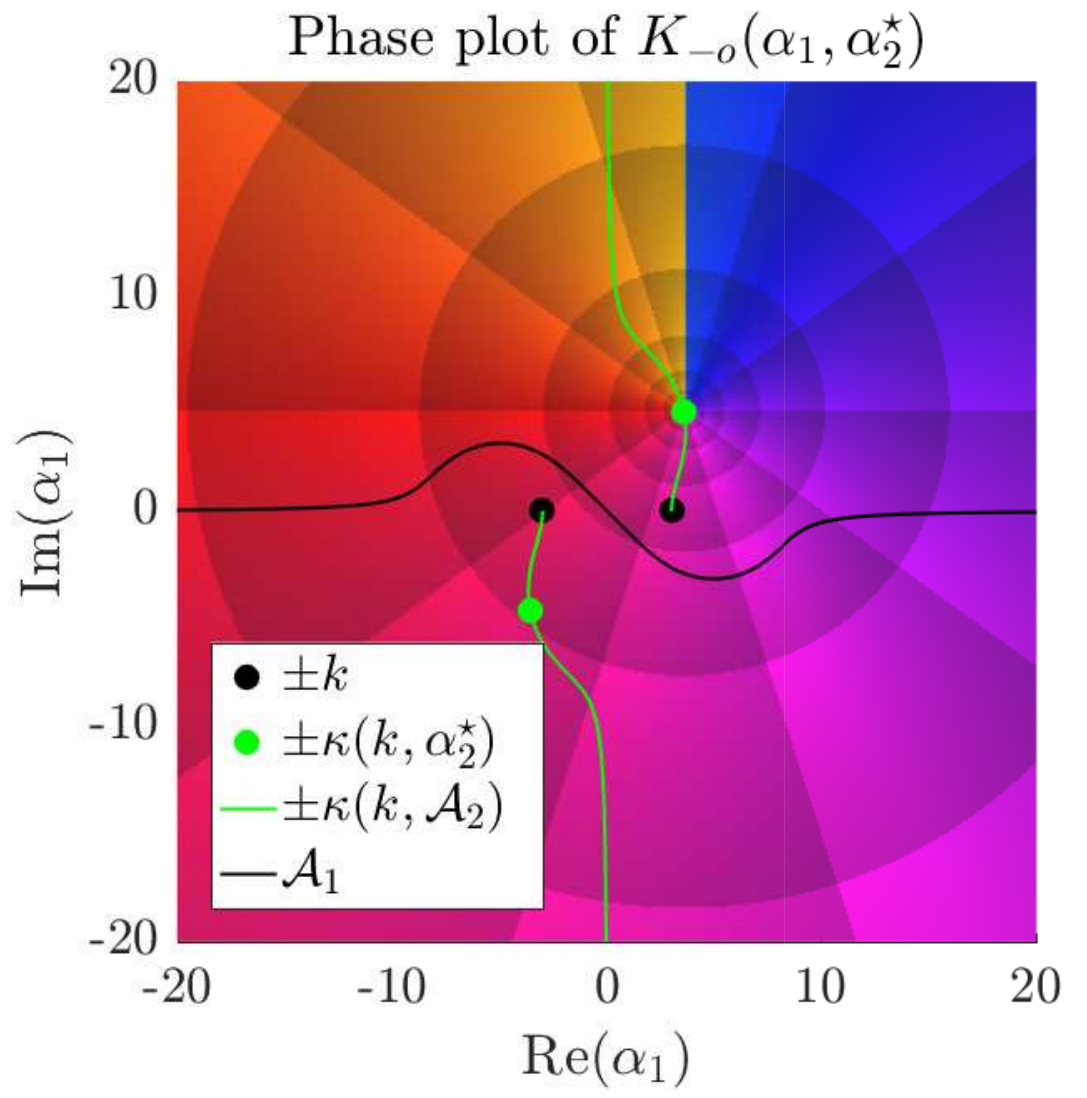}}$\times$ \raisebox{-.5\height}{\includegraphics[width=0.27\textwidth]{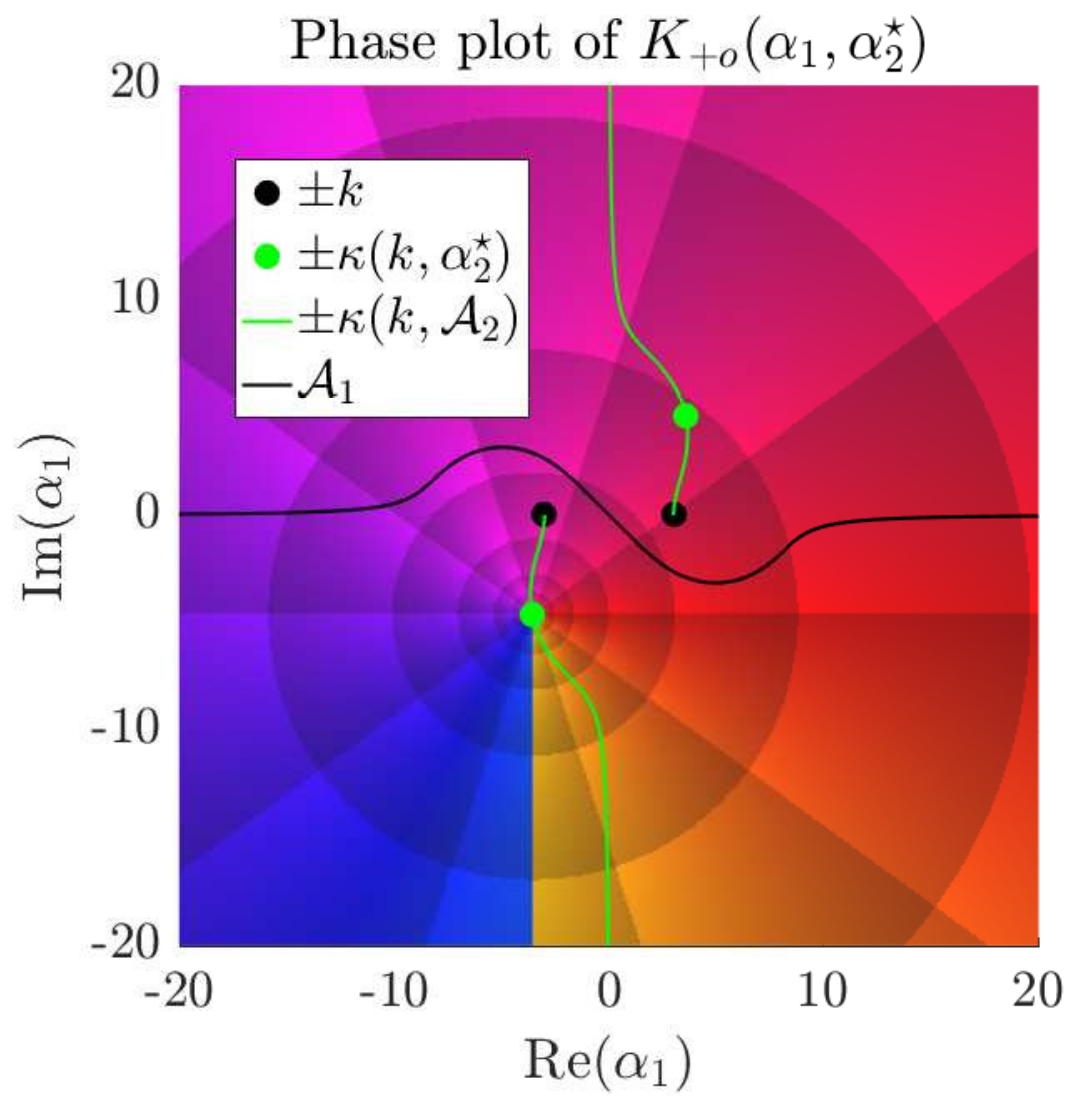}}
		\raisebox{-.5\height}{\includegraphics[width=0.28\textwidth]{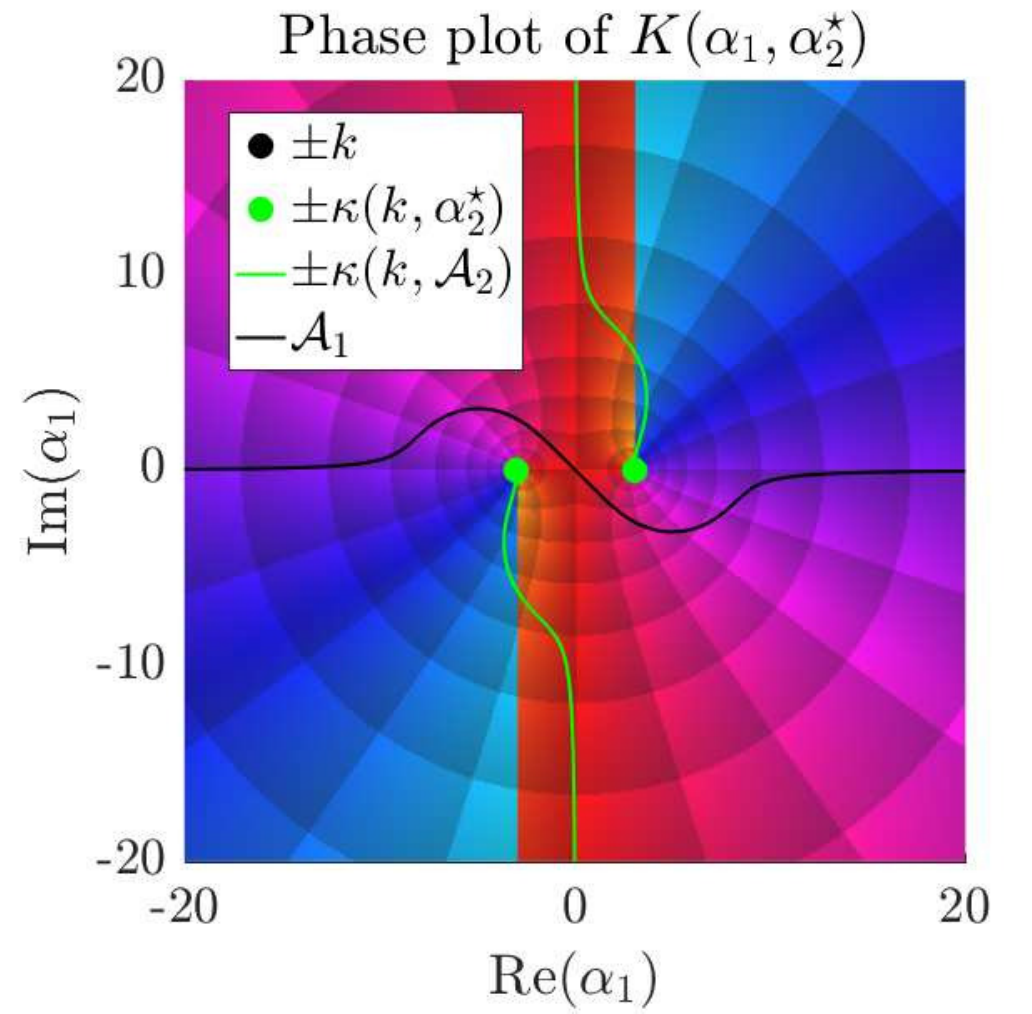}}$=$\raisebox{-.5\height}{\includegraphics[width=0.28\textwidth]{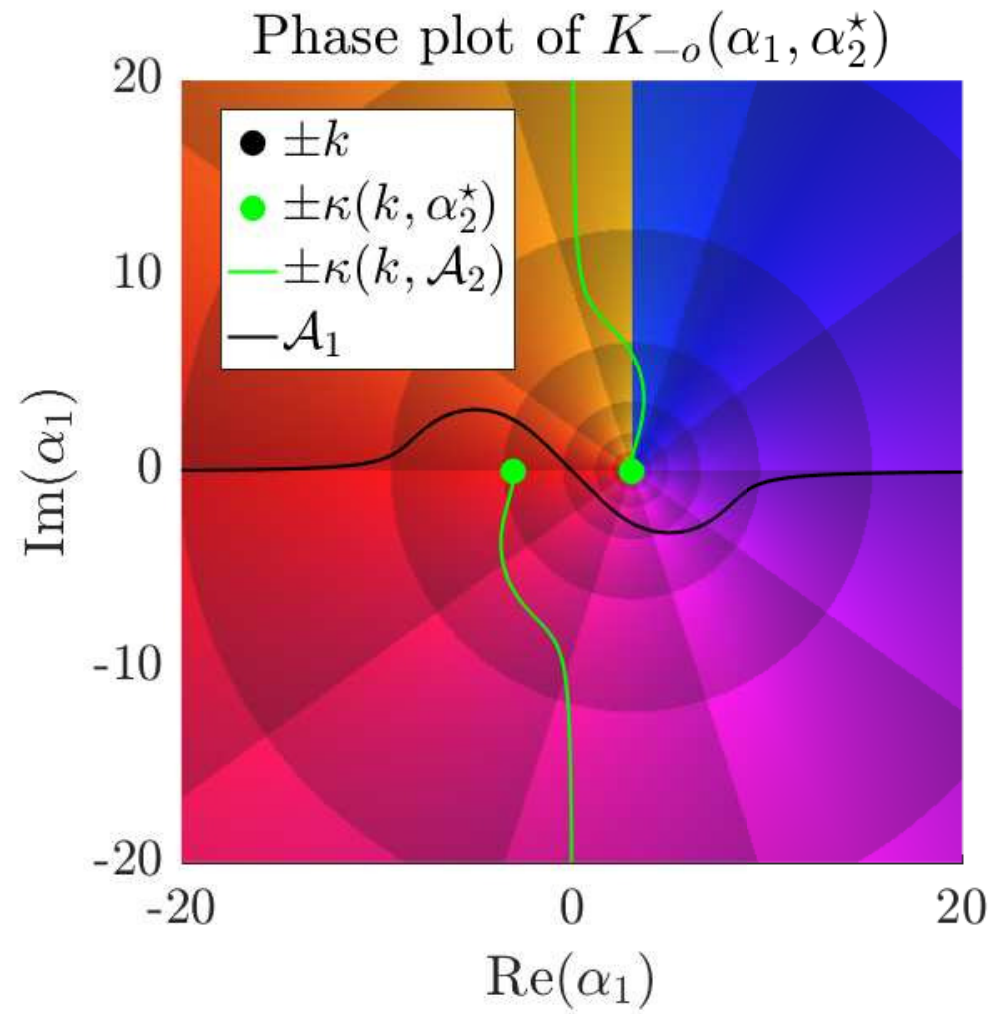}}$\times$ \raisebox{-.5\height}{\includegraphics[width=0.28\textwidth]{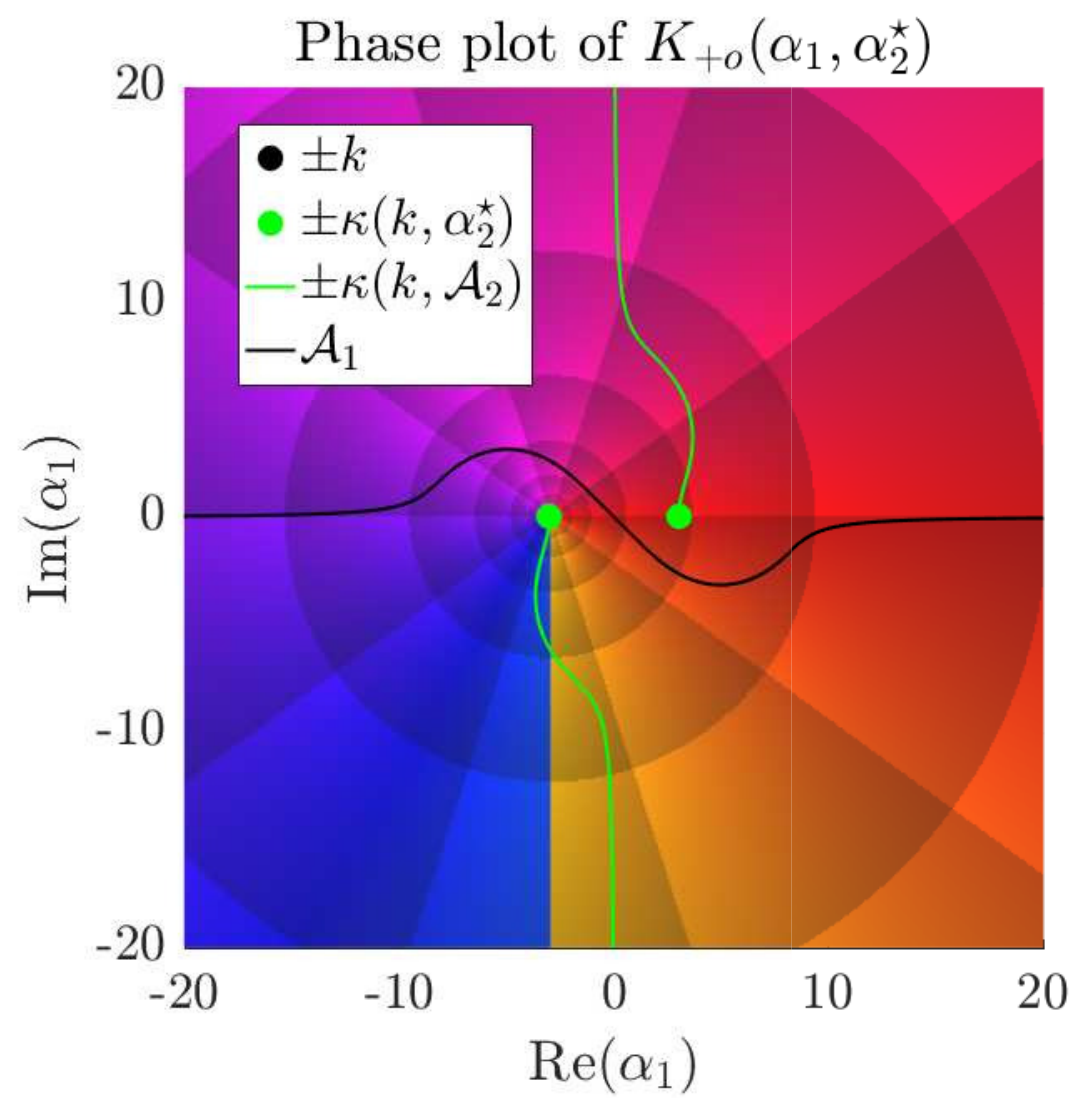}}
		\caption{Plots of the functions $K (\alpha_1, \alpha_2^{\star})$, $K_{-
				\circ} (\alpha_1, \alpha_2^{\star})$ and $K_{+ \circ} (\alpha_1,
			\alpha_2^{\star})$ in the $\alpha_1$ complex plane for $\alpha_2^{\star}
			=\mathcal{A}_2 (5)$ (top) and $\alpha_2^{\star} =\mathcal{A}_2 (0)$
			(bottom).}
		\label{fig:firstKfactorisation}
	\end{figure}
	
	It must be stressed that these functions do not have any useful analyticity properties when
	viewed as functions of $\alpha_2$, with branch cuts passing through both
	$\tmop{UHP}_2$ and $\tmop{LHP}_2$ as $\alpha_1$ moves along $\mathcal{A}_1$. This can be seen in Figure
	\ref{fig:firstKfactinalpha2plane}.
	
	\begin{figure}[htbp]
		\centering
		\raisebox{-.5\height}{\includegraphics[width=0.27\textwidth]{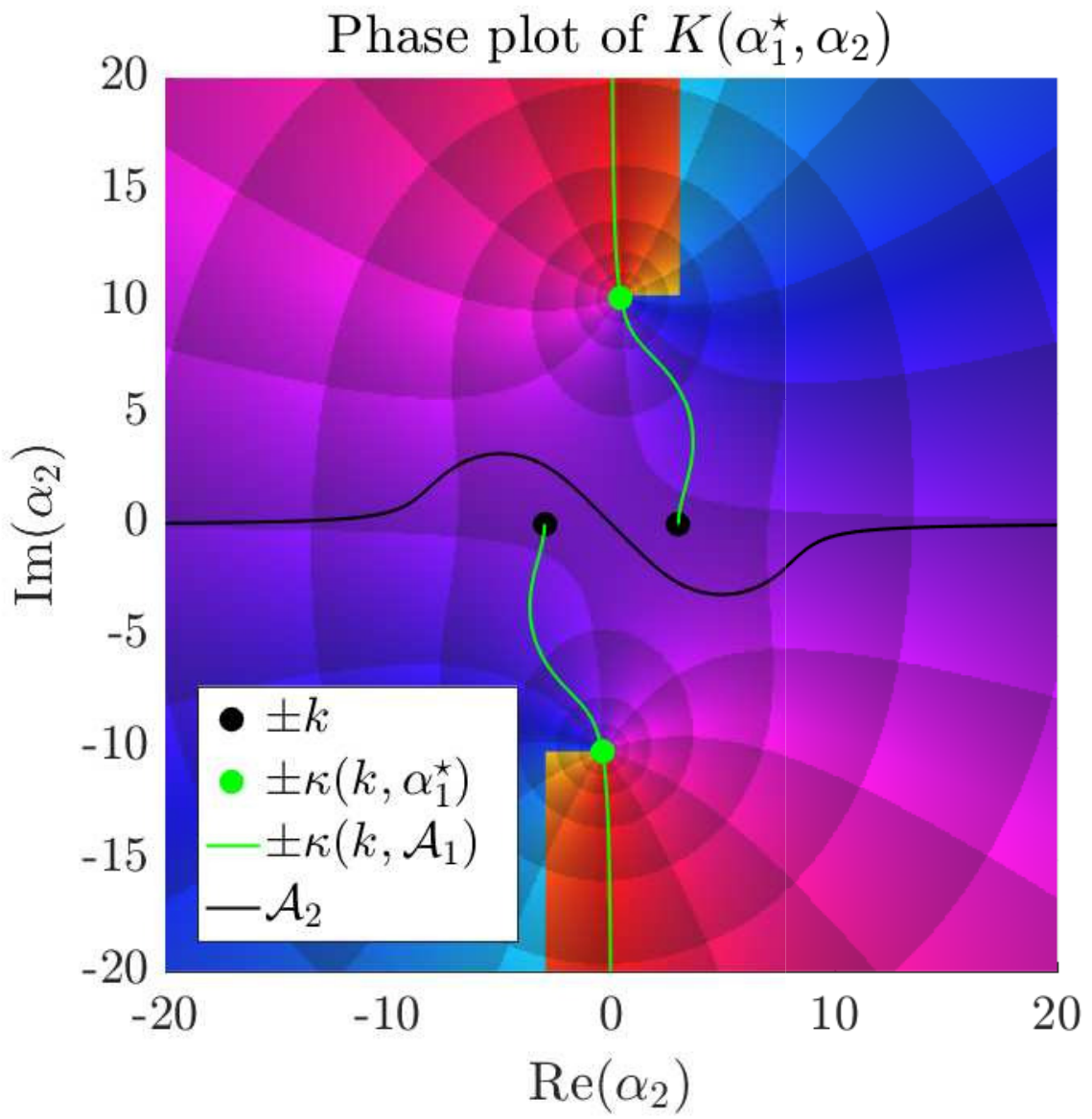}}$=$ \raisebox{-.5\height}{\includegraphics[width=0.27\textwidth]{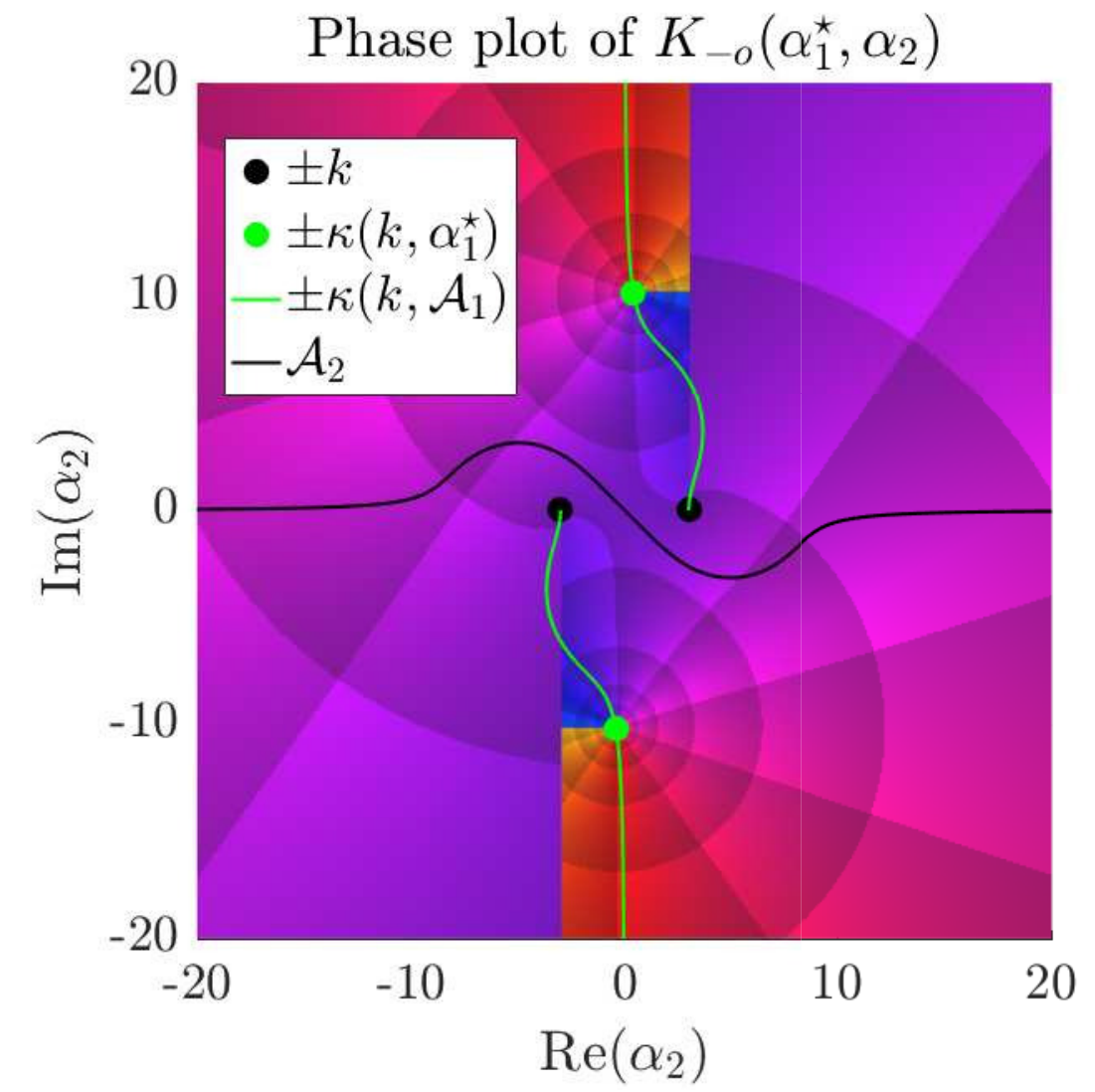}}$\times$ \raisebox{-.5\height}{\includegraphics[width=0.27\textwidth]{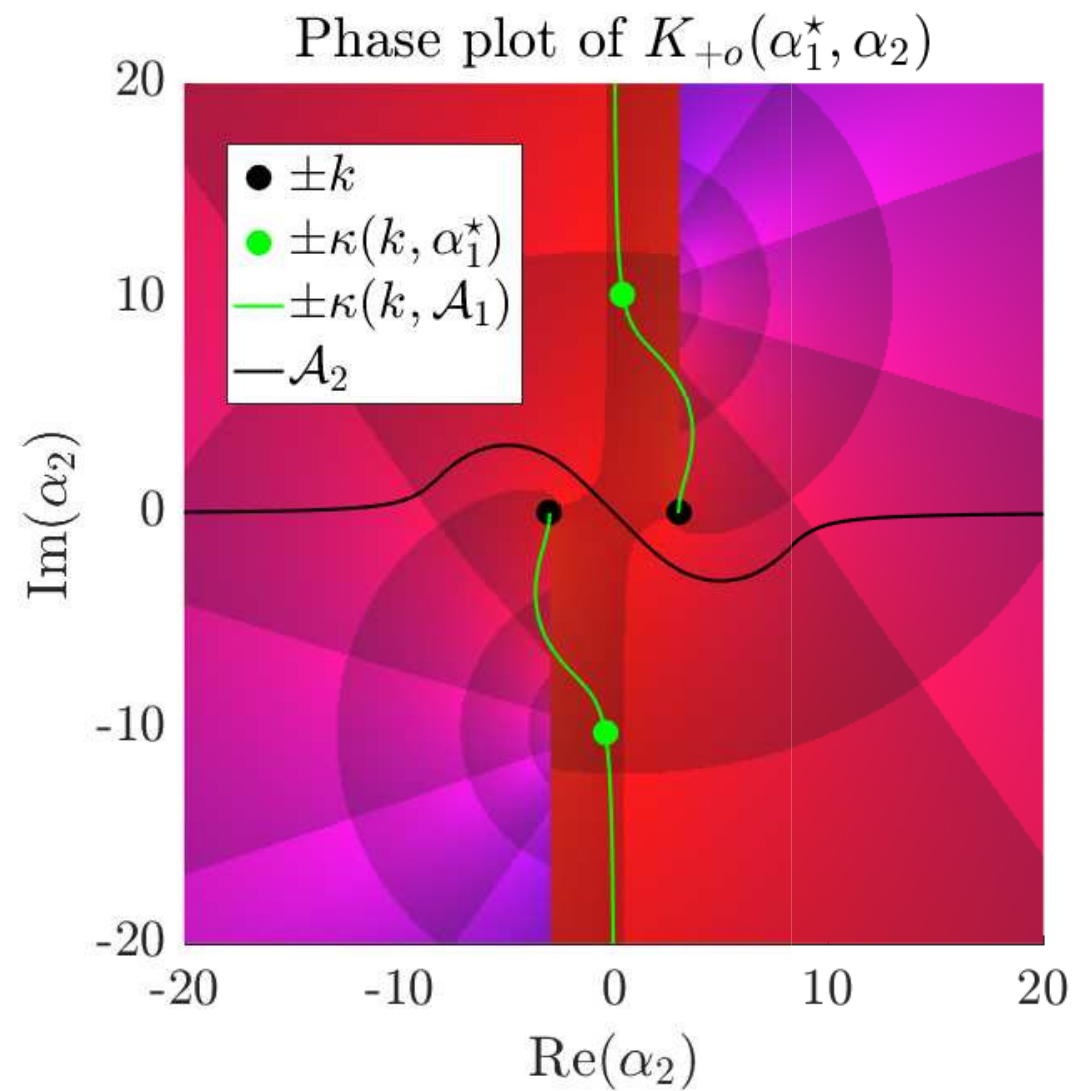}}\\
		
		\raisebox{-.5\height}{\includegraphics[width=0.27\textwidth]{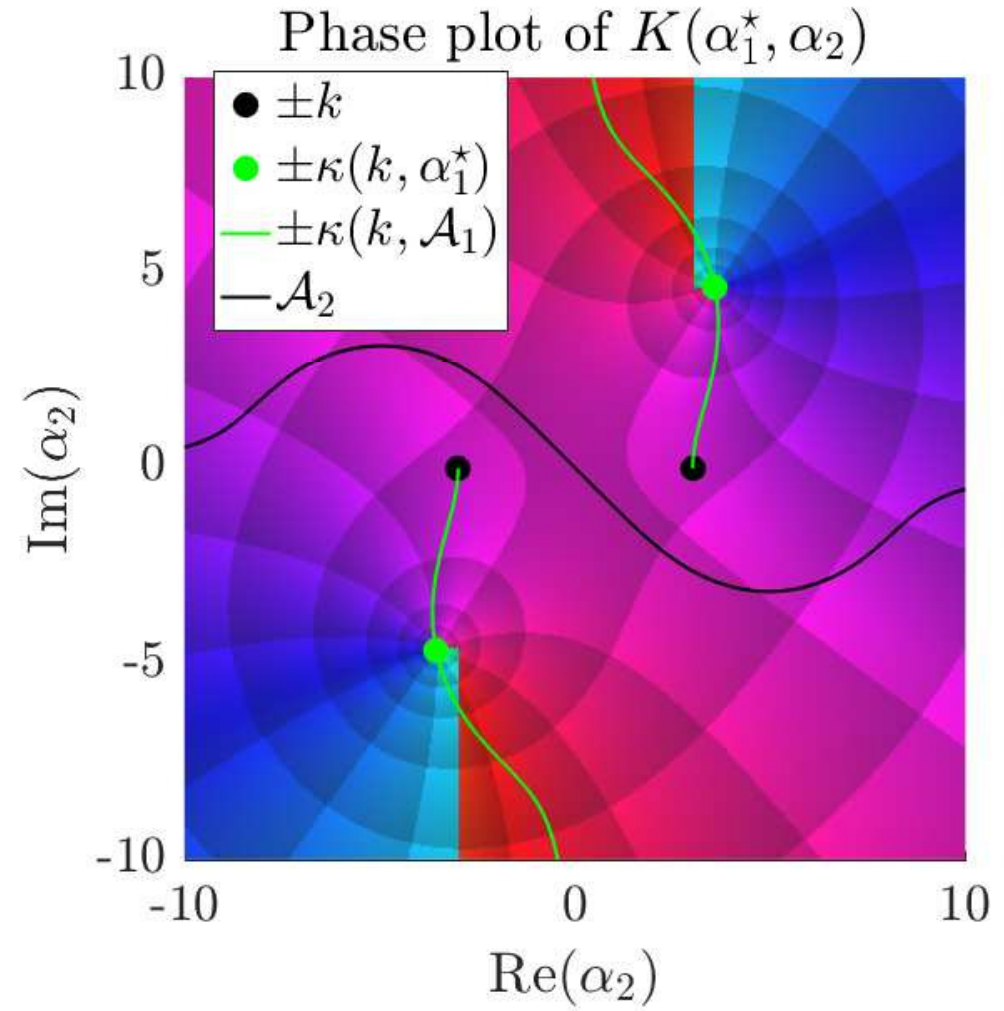}}$=$ \raisebox{-.5\height}{\includegraphics[width=0.27\textwidth]{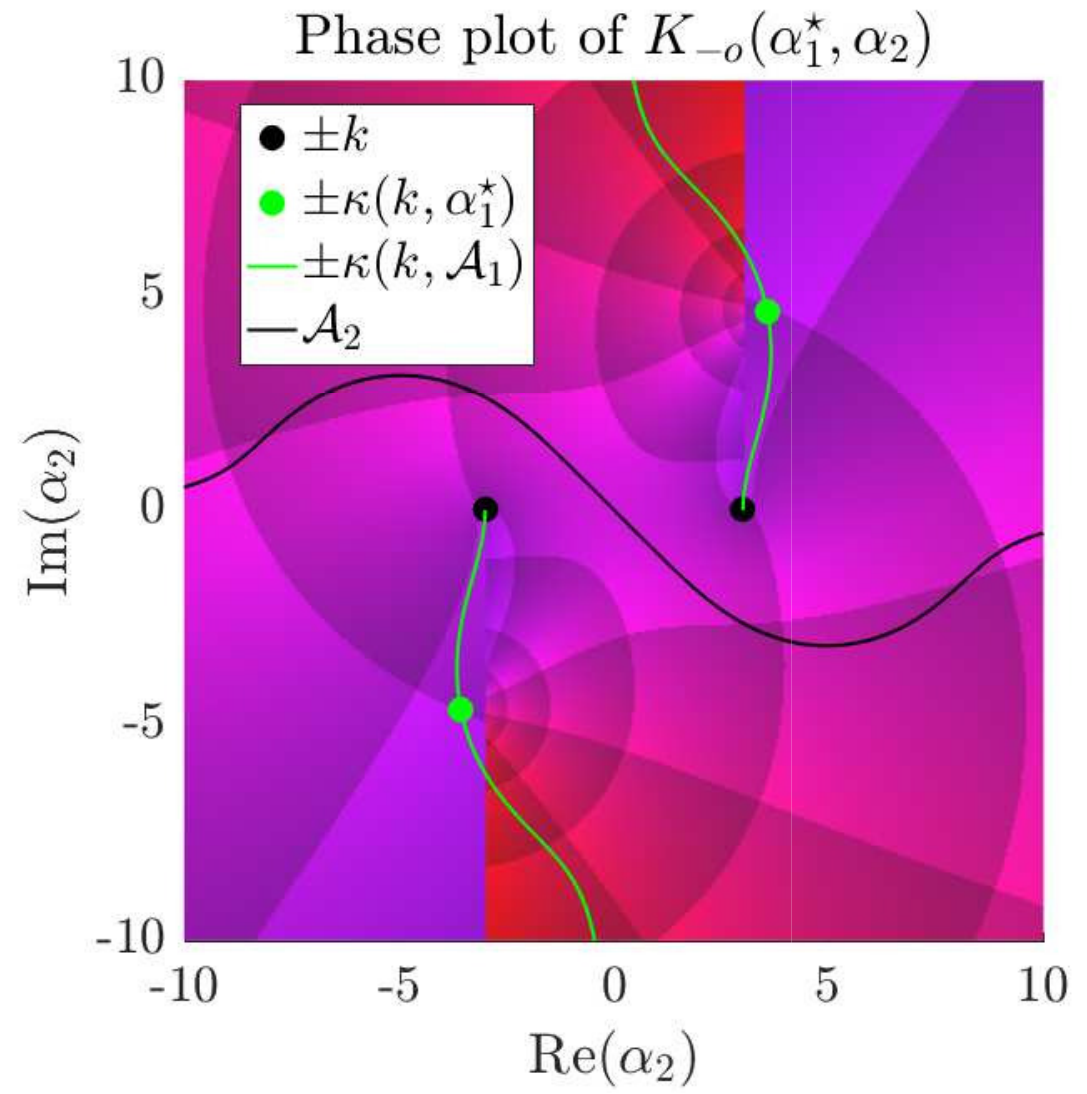}}$\times$ \raisebox{-.5\height}{\includegraphics[width=0.27\textwidth]{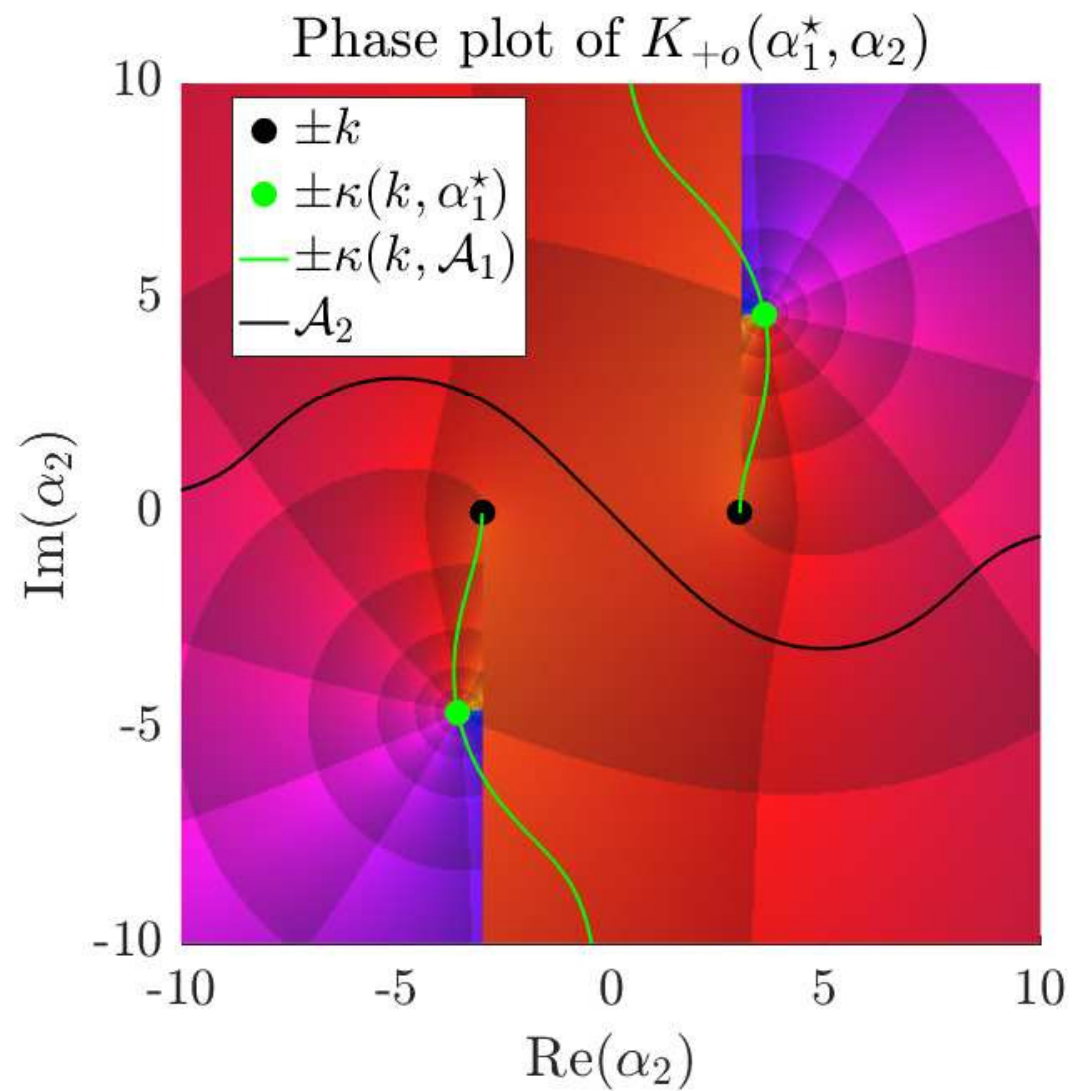}}
		\caption{Plots of the functions $K (\alpha_1^{\star}, \alpha_2)$, $K_{-
				\circ} (\alpha_1^{\star}, \alpha_2^{})$ and $K_{+ \circ} (\alpha_1^{\star},
			\alpha_2^{})$ in the $\alpha_2$ complex plane for $\alpha_1^{\star}
			=\mathcal{A}_1 (10)$ (top) and $\alpha_1^{\star} =\mathcal{A}_1 (5)$
			(bottom).}
		\label{fig:firstKfactinalpha2plane}
	\end{figure}
	
	It is also possible to introduce the functions $K_{\circ -}$ and $K_{\circ +}$
	defined as follows:
	\begin{eqnarray}
	K_{\circ -} (\alpha_1, \alpha_2) = 1/\sqrt[\downarrow]{\kappa (k, \alpha_1) -
		\alpha_2} & \text{ and } & K_{\circ +} (\alpha_1, \alpha_2) =
	1/\sqrt[\downarrow]{\kappa (k, \alpha_1) + \alpha_2}, 
	\label{eq:introothertypefactor}
	\end{eqnarray}
	which will prove useful in \COM{Section} \ref{sec:factK-oK+o}. 
	
	\subsection{Factorisation in the $\alpha_2$-plane}
	
	\subsubsection{Cauchy's formula and its application to factorisation problems}
	
	Let us state two useful results in complex analysis, that we will need in
	this section. The results are classic, and hence, the proofs are omitted. Please refer to e.g. \cite{Noble} for more details. Note that these are valid for a generic complex plane, and since in what we have done so far $\mathcal{A}_1$ and $\mathcal{A}_2$ are the same, we will just denote it by $\mathcal{A}$ in what follows. Similarly, we will use $\text{UHP}$ and $\text{LHP}$ without subscript.
	
	\begin{lemma}
		\label{th:Cauchysumdecomp}{\tmem{{\tmstrong{[Cauchy's formula and
						sum-split]}}}} Let $\Phi$ be a function analytic on a
		(potentially curved) strip $\mathcal{S} \subset \mathbb{C}$ containing
		$\mathcal{A}$, such that we have $\Phi (\alpha) = \Phi_+ (\alpha) +
		\Phi_- (\alpha)$ on $\mathcal{A}$ with $\Phi_+$ analytic on
		$\tmop{UHP}$ and $\Phi_-$ analytic on $\tmop{LHP}$. \COM{C}onsider
		$\mathcal{A}^b_{\varepsilon}$ and $\mathcal{A}^a_{\varepsilon}$ to be the
		contours oriented from left to right defined by $\mathcal{A}^b_{\varepsilon} =\mathcal{A} - i \varepsilon$ and $\mathcal{A}^a_{\varepsilon} =\mathcal{A} + i \varepsilon$, where $\varepsilon > 0$ is any number such that these contours lie within
		$\mathcal{S}$ and the superscripts $a$ and $b$ stand for ``above'' and
		``below'' respectively, as illustrated in Figure \ref{fig:Cauchysundecompcontours}. Let $\alpha \in
		\mathcal{A}$, then, provided that \COM{$\Phi(z)=\mathcal{O}(1/|z|^\lambda)$ for some $\lambda>0$} 
		as $|z|\rightarrow \infty$ within $\mathcal{S}$, the following formulae hold
		\begin{eqnarray*}
			\Phi_+ (\alpha) = \frac{1}{2 i \pi} \int_{\mathcal{A}^b_{
					\varepsilon}} \frac{\Phi (z)}{z - \alpha} \, \mathd z & \text{ and } &
			\Phi_- (\alpha) = \frac{- 1}{2 i \pi} \int_{\mathcal{A}^a_{
					\varepsilon}} \frac{\Phi (z)}{z - \alpha} \, \mathd z \COM{,}
		\end{eqnarray*}
		and can be used to analytically continue $\Phi_+$ ($\Phi_-$) from $\mathcal{A}$ onto $\tmop{UHP}$ ($\tmop{LHP}$).
	\end{lemma}
	\begin{figure}[htbp]
		\centering
		\includegraphics[width=0.45\textwidth]{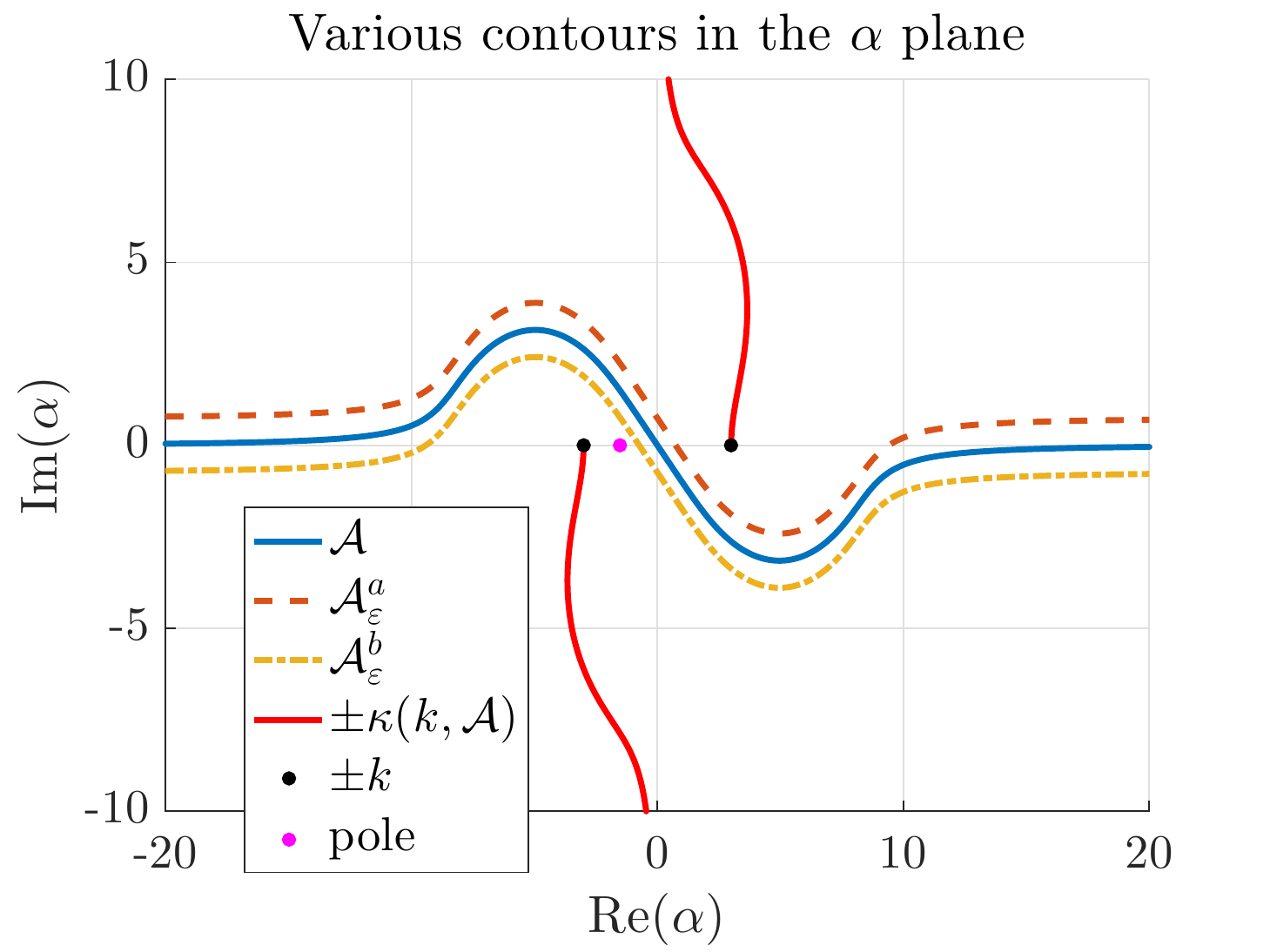}
		\caption{Diagrammatic illustrations of the contours introduced in Lemma
			\ref{th:Cauchysumdecomp}}
		\label{fig:Cauchysundecompcontours}
	\end{figure}
	\begin{corollary}
		\label{cor:CauchyFac}{\tmem{{\tmstrong{[Cauchy's formula and
						factorisation]}}}} Let $\Psi$ be a function analytic on a (potentially
		curved) strip $\mathcal{S} \subset \mathbb{C}$ containing $\mathcal{A}$,
		such that we have $\Psi (\alpha) = \Psi_+ (\alpha) \Psi_- (\alpha)$ on
		$\mathcal{A}$ with $\Psi_+$ analytic on $\tmop{UHP}$ and $\Psi_-$
		analytic on $\tmop{LHP}$. Let $\alpha \in \mathcal{A}$, then, provided that $\Psi(z)\rightarrow 1$ as $|z| \rightarrow \infty$ within $\mathcal{S}$, the following formulae hold
		\begin{eqnarray*}
			\Psi_+ (\alpha) = \exp \left\{ \frac{1}{2 i \pi}
			\int_{\mathcal{A}^b_{\varepsilon}} \frac{\log (\Psi (z))}{z - \alpha}
			\, \mathd z \right\}&\text{ and }&
			\Psi_- (\alpha)  = \exp \left\{ \frac{- 1}{2 i \pi}
			\int_{\mathcal{A}^a_{\varepsilon}} \frac{\log (\Psi (z))}{z - \alpha}
			\, \mathd z \right\},
		\end{eqnarray*}
		where $\mathcal{A}_\varepsilon^{a,b}$ are defined as in Lemma \ref{th:Cauchysumdecomp}, and can be used to analytically continue $\Psi_+$ from $\mathcal{A}$ onto $\tmop{UHP}$ and
		$\Psi_-$ from $\mathcal{A}$ onto $\tmop{LHP}$.
	\end{corollary}
	\subsubsection{Factorisation of $K_{- \circ}$ and $K_{+
			\circ}$}\label{sec:factK-oK+o}
	
	It does not seem possible to find an {\tmem{explicit}} factorisation of these
	functions. Nevertheless, a direct application of Cauchy's formulae does lead
	to a formal factorisation of $K_{- \circ}$ and $K_{+ \circ}$ in the $\alpha_2$
	plane. However, the resulting expressions can be quite slow to evaluate
	numerically. In \RED{Appendix} \ref{app:factoK-o}, we perform some manipulations of
	the integrals in order to obtain forms that are rapid to compute; these are employed in (\ref{eq:integralK-+})-(\ref{eq:integralK+-}). 
	$K_{- \circ}$ can be factorised as $K_{- \circ} (\tmmathbf{\alpha}) = K_{- +} (\tmmathbf{\alpha}) K_{- -}
	(\tmmathbf{\alpha})$, and $K_{+ \circ}$ can be factorised as $K_{+ \circ} (\tmmathbf{\alpha})= K_{+ +} (\tmmathbf{\alpha}) K_{+ -}
	(\tmmathbf{\alpha})$,
	where we have
	\begin{eqnarray}
	K_{- +} (\boldsymbol{\alpha}) & = & \frac{1}{\sqrt[\downarrow]{\sqrt[\downarrow]{k + \alpha_2}}}\exp \left\{ \frac{- 1}{4 i \pi}
	\int_{\mathcal{A}_{\varepsilon}^b} \frac{\overset{\swarrow}{\log}
		\left( 1 - \frac{\alpha_1}{\kappa (k, z)} \right)}{z - \alpha_2} \, \mathd z
	\right\} \text{ for } \boldsymbol{\alpha}\in\mathcal{D_{-+}},
	\label{eq:integralK-+} \\
	K_{- -} (\boldsymbol{\alpha}) & = & \frac{1}{\sqrt[\downarrow]{\sqrt[\downarrow]{k - \alpha_2}}}\exp \left\{ \frac{1}{4 i \pi}
	\int_{\mathcal{A}_{\varepsilon}^a} \frac{\overset{\swarrow}{\log}
		\left( 1 - \frac{\alpha_1}{\kappa (k, z)} \right)}{z - \alpha_2} \, \mathd z
	\right\} \text{ for } \boldsymbol{\alpha}\in\mathcal{D_{--}},
	\label{eq:integralK--} \\
	K_{+ +} (\boldsymbol{\alpha}) & = &\frac{1}{\sqrt[\downarrow]{\sqrt[\downarrow]{k + \alpha_2}}}\exp \left\{ \frac{- 1}{4 i \pi}
	\int_{\mathcal{A}_{\varepsilon}^b} \frac{\overset{\swarrow}{\log}
		\left( 1 + \frac{\alpha_1}{\kappa (k, \COM{z})} \right)}{z - \alpha_2} \,
	\mathd z \right\} \text{ for } \boldsymbol{\alpha}\in\mathcal{D_{++}},
	\label{eq:integralK++} \\
	K_{+ -} (\boldsymbol{\alpha}) & = & \frac{1}{\sqrt[\downarrow]{\sqrt[\downarrow]{k - \alpha_2}}}\exp \left\{ \frac{1}{4 i \pi}
	\int_{\mathcal{A}_{\varepsilon}^a} \frac{\overset{\swarrow}{\log}
		\left( 1 + \frac{\alpha_1}{\kappa (k, \COM{z})} \right)}{z - \alpha_2} \,
	\mathd z \right\} \text{ for } \boldsymbol{\alpha}\in\mathcal{D_{+-}}. 
	\label{eq:integralK+-}
	\end{eqnarray}
	These formulae allow for a fast evaluation of the four components of the
	factorisation of $K$, allowing us to gain a good {\tmem{visual}} understanding
	of the singularity structure of $K_{- +}$, $K_{- -}$, $K_{+ +}$ and $K_{+ -}$,
	as illustrated in Figures \ref{fig:visualfactK-o} and \ref{fig:visualfactK+o}.
	To give an idea of the speed, for each plot we need to evaluate the functions
	160,000 times and it takes about 14 seconds to run on a standard laptop.
	
	Another method (\RED{see e.g.\ }{\cite{Albertsen1997}}), involving the Dilog function, has also been used to evaluate these factors. \RED{Both methods are very fast to evaluate, though, upon implementing them both in Matlab, it transpires that ours leads to a faster evaluation of $K_{++}$ say. Moreover, our formula (\ref{eq:integralK++}) giving $K_{++}$ is more compact than that involving the Dilog function.}
	
	\begin{figure}[htbp]
		\centering
		\raisebox{-.5\height}{\includegraphics[width=0.3\textwidth]{pictures/K-o_alpha2_plane_scaled-eps-converted-to.pdf}}$=$\raisebox{-.5\height}{\includegraphics[width=0.3\textwidth]{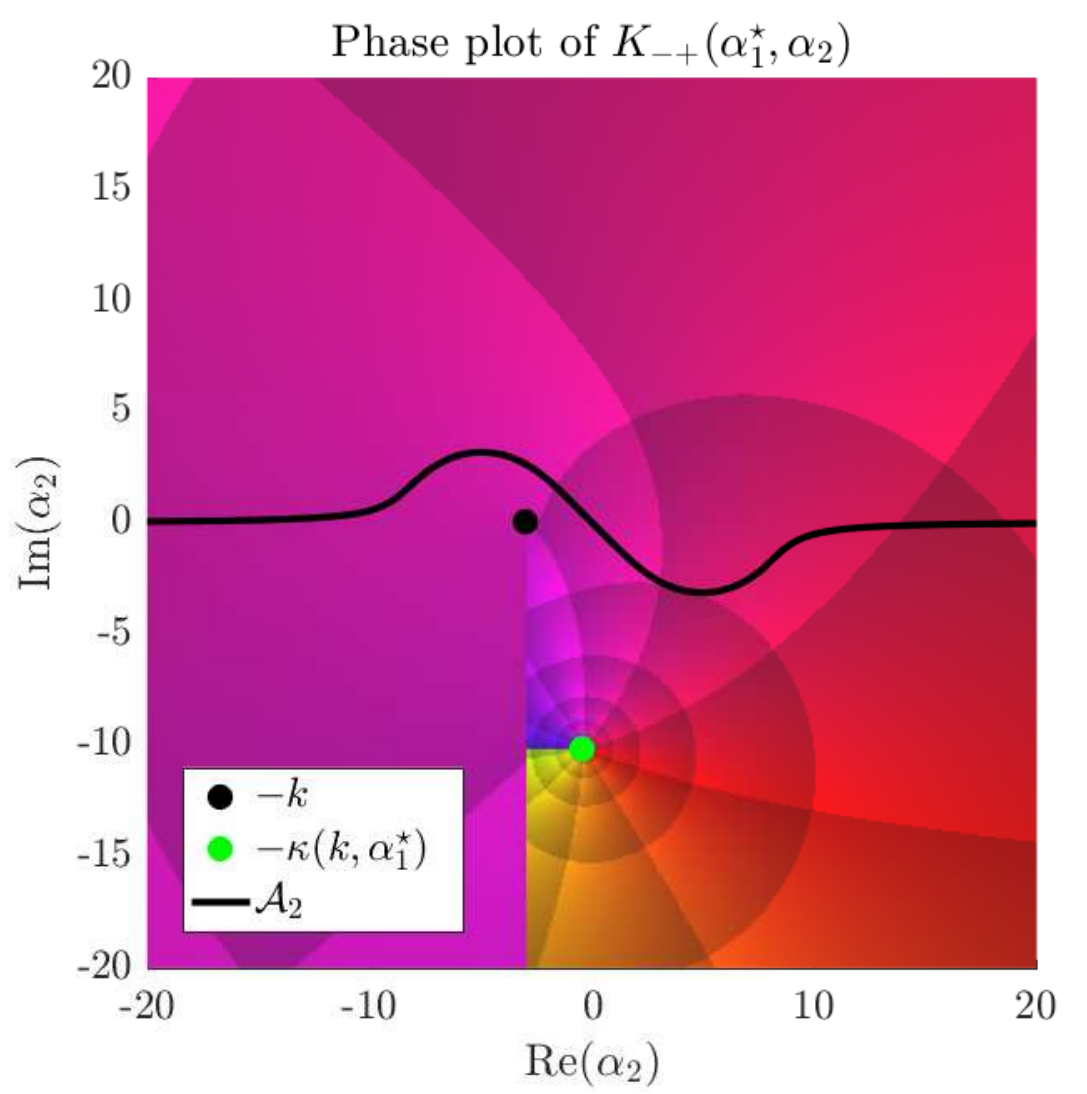}}$\times$\raisebox{-.5\height}{\includegraphics[width=0.3\textwidth]{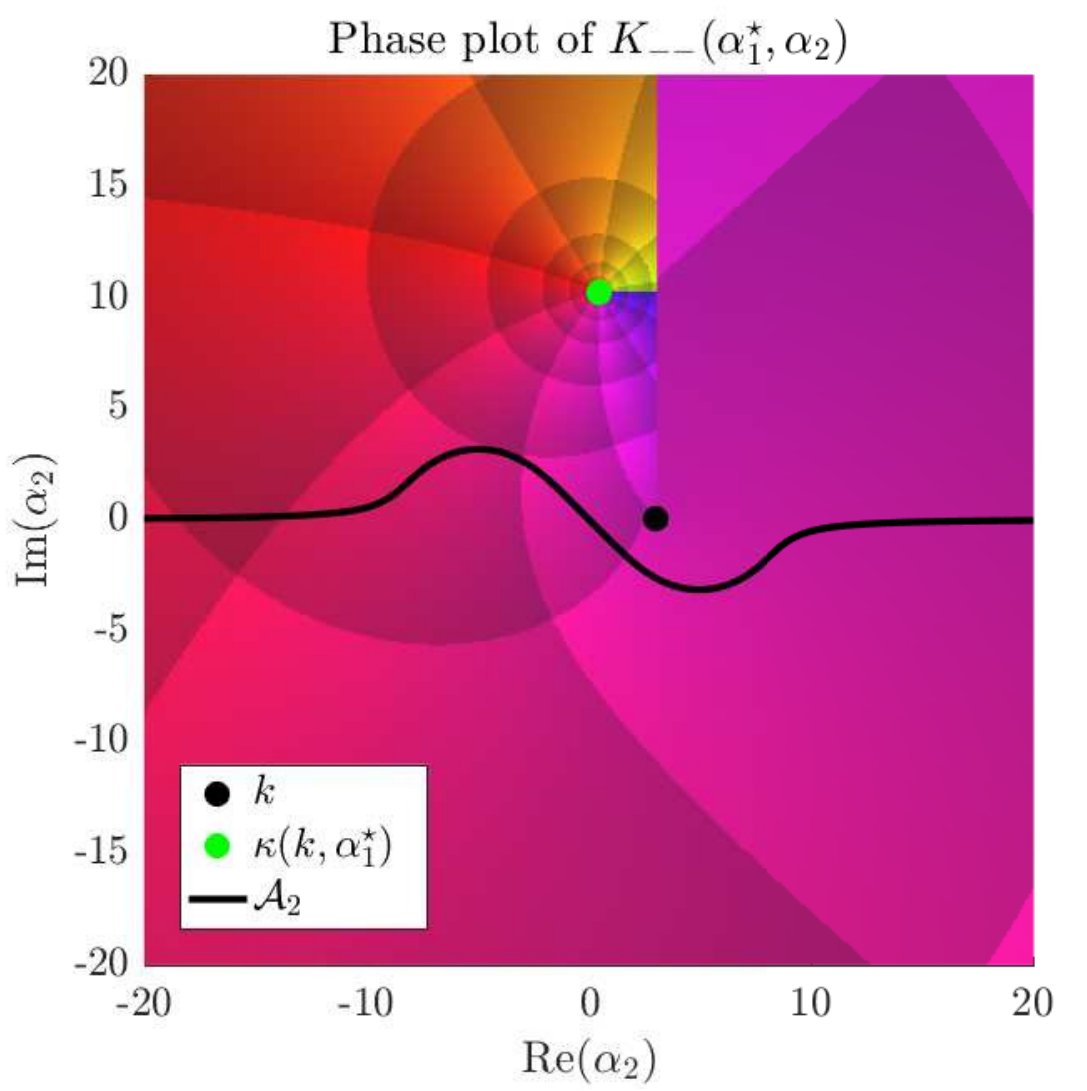}}
		\caption{Plots of the functions $K_{- \circ} (\alpha_1^{\star}, \alpha_2)$,
			$K_{- +} (\alpha_1^{\star}, \alpha_2^{})$ and $K_{- -} (\alpha_1^{\star},
			\alpha_2^{})$ in the $\alpha_2$ complex plane for $\alpha_1^{\star}
			=\mathcal{A}_1 (10)$. In its region of analyticity, $\tmop{UHP}_2$, $K_{-
				+}$ has been obtained via (\ref{eq:integralK-+}), while in $\tmop{LHP}_2$,
			it has been obtained by analytical continuation using $K_{- +} = K_{- \circ}
			/ K_{- -}$. A similar strategy has been used to plot $K_{- -}$.}
		\label{fig:visualfactK-o}
	\end{figure}
	
	\begin{figure}[htbp]
		\centering
		\raisebox{-.5\height}{\includegraphics[width=0.3\textwidth]{pictures/K+o_alpha2_plane_scaled-eps-converted-to.pdf}}$=$ \raisebox{-.5\height}{\includegraphics[width=0.3\textwidth]{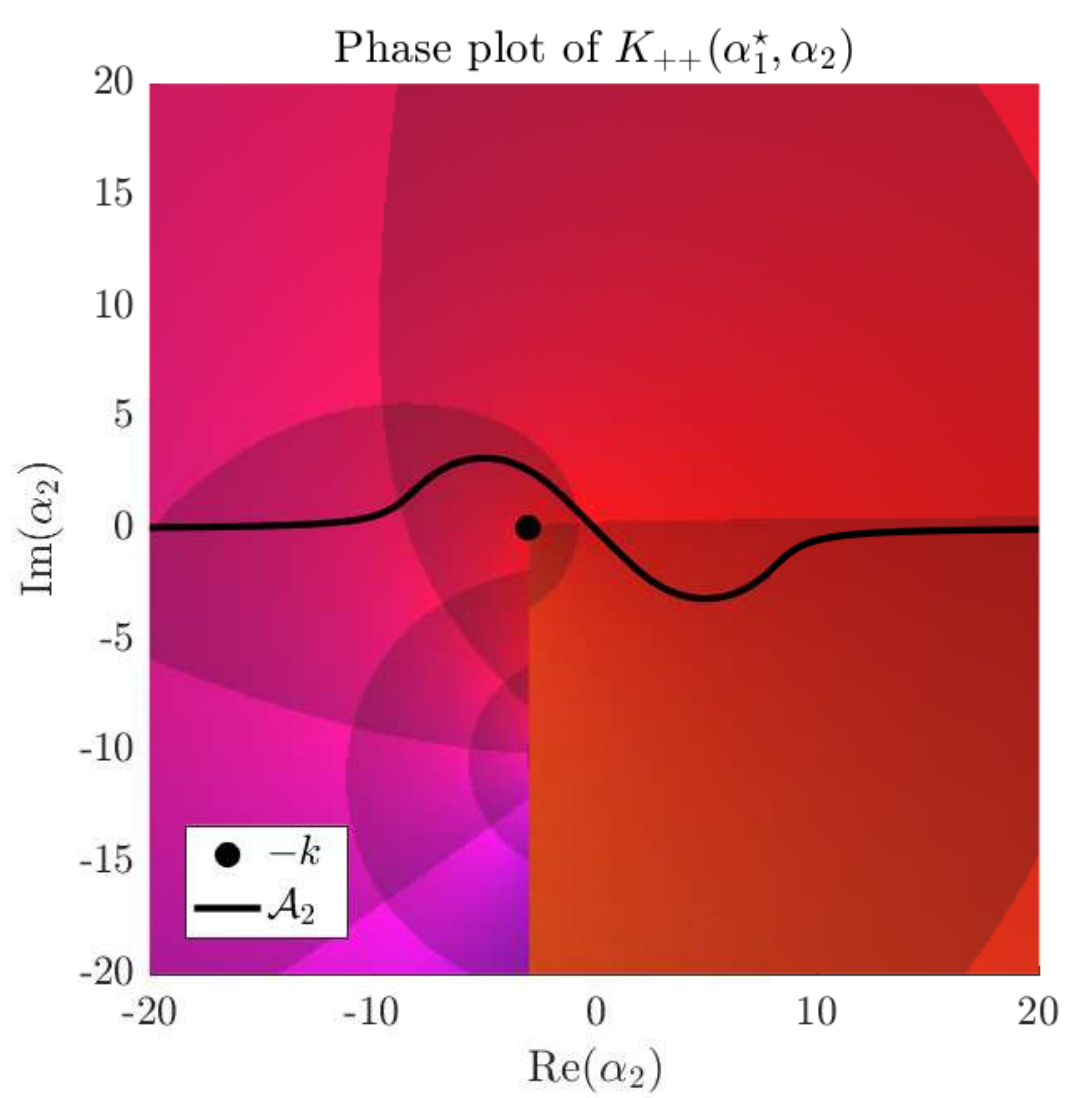}}$\times$ \raisebox{-.5\height}{\includegraphics[width=0.3\textwidth]{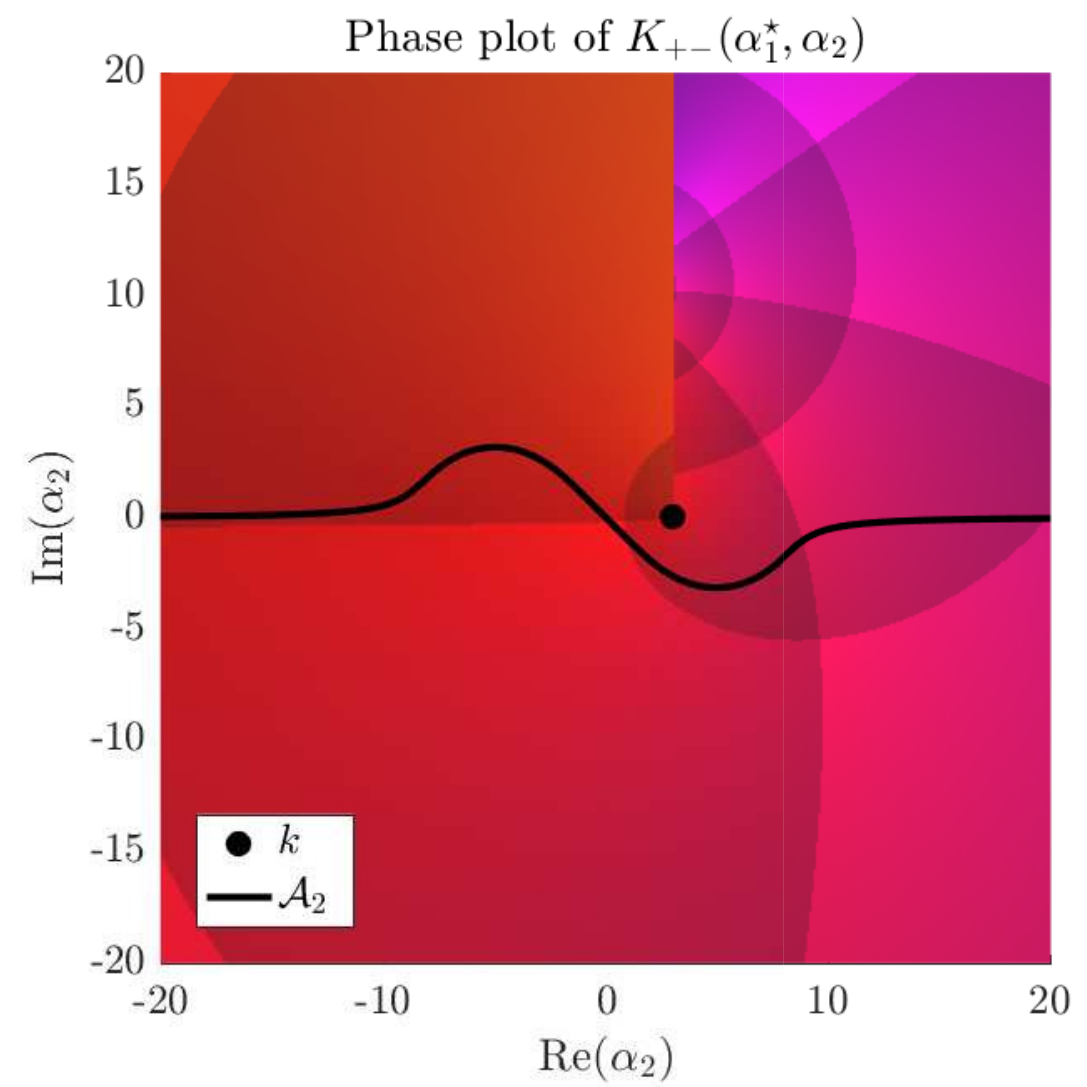}}
		\caption{Plots of the functions $K_{+ \circ} (\alpha_1^{\star}, \alpha_2)$,
			$K_{+ +} (\alpha_1^{\star}, \alpha_2^{})$ and $K_{+ -} (\alpha_1^{\star},
			\alpha_2^{})$ in the $\alpha_2$ complex plane for $\alpha_1^{\star}
			=\mathcal{A}_1 (10)$. In its region of analyticity, $\tmop{UHP}_2$, $K_{+
				+}$ has been obtained via (\ref{eq:integralK++}), while in $\tmop{LHP}_2$,
			it has been obtained by analytical continuation using $K_{+ +} = K_{+ \circ}
			/ K_{+ -}$. A similar strategy has been used to plot $K_{+ -}$.}
		\label{fig:visualfactK+o}
	\end{figure}
	
	On Figures \ref{fig:visualfactK-o} and \ref{fig:visualfactK+o}, $\alpha_1^\star$ has been chosen on $\mathcal{A}_1$ for illustration purpose, but it could have been chosen anywhere in $\text{LHP}_1$ for Figure \ref{fig:visualfactK-o} and anywhere in $\text{UHP}_1$ for Figure \ref{fig:visualfactK+o}. We chose to visualise this factorisation in the $\alpha_2$ plane, but it is also possible to visualise it in the $\alpha_1$ plane for a given $\alpha_2^\star$ on $\mathcal{A}_2$. In this case, in order to analytically continue the factors past their natural domain of analyticity, one should use the functions $K_{\circ \pm}$ introduced in (\ref{eq:introothertypefactor}).
	
	\section{The (generic) Wiener-Hopf system in
		$\mathbb{C}^2$}\label{sec:successiveWH}
	
	\subsection{Quadruple sum-split}
	
	\COM{Using the function $F_{+ +}$ defined in (\ref{eq:defF++intermsofF}), the functional equation (\ref{eq:revisionfunctionaleqraph}) can be rewritten $G(\tmmathbf{\alpha})=F_{+ +} (\tmmathbf{\alpha}) K (\tmmathbf{\alpha})$ and, as seen in Section \ref{sec:Fourier}, $G$ is analytic on $\mathcal{D}=\mathcal{A}_1\times\mathcal{A}_2$. Hence, we can\footnote{A more rigorous approach to obtain this would be to refer to Bochner's theorem {\cite{Bochner1938}}.} write its additive
		decomposition}
	\begin{align}
	F_{+ +} (\tmmathbf{\alpha}) K (\tmmathbf{\alpha}) \COM{=G(\boldsymbol{\alpha})} =  G_{+ +}
	(\tmmathbf{\alpha}) + G_{- +} (\tmmathbf{\alpha}) + G_{- -}
	(\tmmathbf{\alpha}) + G_{\noplus + -} (\tmmathbf{\alpha})_{}, 
	\label{eq:wh2}
	\end{align}
	where $G_{+ +} (\tmmathbf{\alpha})$, $G_{- +} (\tmmathbf{\alpha})$, $G_{- -} (\tmmathbf{\alpha})$  and  $G_{ + -} (\tmmathbf{\alpha})$ are analytic on $\mathcal{D}_{+ +}$,
	$\mathcal{D}_{- +}$, $\mathcal{D}_{- -}$ and $\mathcal{D}_{+ -}$ respectively.
	\COM{Note that by definition of $G(\boldsymbol{\alpha})$, see (\ref{eq:solutioninFourierspace}), we have $G(\boldsymbol{\alpha})=\mathfrak{F} [u
		(x_1, x_2, 0)] (\tmmathbf{\alpha})$},
	where $\mathfrak{F}$ is the double Fourier transform operator as defined in
	\COM{Section} \ref{sec:doubleFourier}. \COM{Therefore, upon defining} the functions $u_j$, $j = 1
	\ldots 4$, by
	\begin{eqnarray*}
		u_j (x_1, x_2) = u (x_1, x_2, 0) H_j (x_1, x_2), & \text{ where } & H_j (x_1,
		x_2) = \left\{ \begin{array}{l}
			1 \text{ if } (x_1, x_2) \in Q_j, \\
			0 \text{ otherwise} \COM{,}
		\end{array} \right.
	\end{eqnarray*}
	\COM{i}t is then possible to define the additive terms as quarter-range Fourier
	transform\COM{s:}
	\begin{eqnarray}
	G_{+ +} (\tmmathbf{\alpha}) =\mathfrak{F} [u_1 (x_1, x_2)]
	(\tmmathbf{\alpha}),&\text{ }& G_{- +} (\tmmathbf{\alpha}) =\mathfrak{F} [u_2
	(x_1, x_2)] (\tmmathbf{\alpha}), \nonumber\\
	G_{- -} (\tmmathbf{\alpha}) =\mathfrak{F} [u_3 (x_1, x_2)]
	(\tmmathbf{\alpha}), &\text{ }& G_{+ -} (\tmmathbf{\alpha}) =\mathfrak{F} [u_4
	(x_1, x_2)] (\tmmathbf{\alpha}) .  \label{eq:G+-definition}
	\end{eqnarray}
	We can also define the auxiliary functions $G_{+ \circ} = G_{+ +} + G_{+ -}$ and $G_{- \circ} = G_{- +} + G_{-
		-}$ that are analytic on $\mathcal{D}_{+ \circ}$ and $\mathcal{D}_{- \circ}$
	respectively.
	
	\subsection{On the function $G_{+ +}$}
	
	Because we impose the Dirichlet condition (\ref{eq:Dirichletcondition}), it
	follows that we have
	\begin{eqnarray*}
		u_1 (x_1, x_2) & = & - u_{\tmop{in}} (x_1, x_2, 0) H_1 (x_1, x_2) = - e^{- i
			(a_1 x_1 + a_2 x_2)} H_1 (x_1, x_2)
	\end{eqnarray*}
	and so, since $G_{+ +}$ is defined on $\mathcal{D}_{+ +}$ by $G_{+ +} (\tmmathbf{\alpha})=\mathfrak{F} [u_1 (x_1, x_2)](\tmmathbf{\alpha})$, we obtain
	\begin{eqnarray}
	G_{+ +} (\tmmathbf{\alpha}) & = & \frac{1}{(\alpha_1 - a_1) (\alpha_2 -
		a_2)} \text{ } . \label{eq:exactfunctionG++}
	\end{eqnarray}
	
	Note that each pole must lie in its respective lower-half plane whether $a_{1,2}$ is positive or negative in order to ensure that $G_{++}$ is analytic in $\mathcal{D}_{++}$. As discussed in Remark \ref{note:introtildecontour}, when $a_{1,2}$ is positive, we allow it to have a small imaginary part, $\epsilon<0$, which places it below $\mathcal{A}_{1,2}$, and then later allow $\epsilon\rightarrow 0$.
	
	
	Hence, at the moment, we have four unknown functions, namely $F_{+ +}$, $G_{+
		-}$, $G_{- +}$ and $G_{- -}$. In the following two subsections, we will show
	how (\ref{eq:wh2}) can be reduced to four equations, involving our four
	unknowns.
	
	\subsection{A first split\COM{,} in the $\alpha_1$ plane}
	
	Let us start by rewriting\footnote{For brevity we will only specify the argument of the functions involved if it is not $\tmmathbf{\alpha}$\COM{.}} (\ref{eq:wh2}) as follows:
	\begin{eqnarray*}
		F_{+ +} K_{+ \circ} K_{- \circ} & = & G_{+ +} + G_{- \circ} + G_{+ -} \, .
	\end{eqnarray*}
	Upon dividing by $K_{- \circ}$, we obtain
	\begin{eqnarray}
	F_{+ +} K_{+ \circ} & = & G_{+ +}/K_{- \circ} + G_{-
		\circ}/K_{- \circ} + G_{+ -}/K_{- \circ}. 
	\label{eq:otherwaywh1}
	\end{eqnarray}
	Now, formally, using for example Lemma \ref{th:Cauchysumdecomp}, it is
	possible to perform a sum-split in the $\alpha_1$-plane of the terms $G_{+ +} / K_{- \circ}$ and $G_{+
		-} / K_{- \circ}$ by writing
	\begin{eqnarray*}
		\frac{G_{+ +}}{K_{- \circ}} = \left[ \frac{G_{+ +}}{K_{- \circ}}
		\right]_{+ \circ} \!\!+ \left[ \frac{G_{+ +}}{K_{- \circ}} \right]_{- \circ} & \!\!\text{ and }&
		\frac{G_{+ -}}{K_{- \circ}} = \left[ \frac{G_{+ -}}{K_{- \circ}}
		\right]_{+ \circ} \!\!+ \left[ \frac{G_{+ -}}{K_{- \circ}} \right]_{- \circ},
	\end{eqnarray*}
	where the operators $\left[ \ \right]_{- \circ}$ and $\left[
	\ \right]_{+ \circ}$ represent respectively the $\alpha_1$-minus part and
	$\alpha_1$-plus part of a given function that is analytic on $\mathcal{A}_1$
	when considered a function of $\alpha_1$. 
	With this split, (\ref{eq:otherwaywh1}) may be rearranged as
	\begin{eqnarray}
	F_{+ +} K_{+ \circ} - \left[ \frac{G_{+ +}}{K_{- \circ}} \right]_{+ \circ} \!\!-
	\left[ \frac{G_{+ -}}{K_{- \circ}} \right]_{+ \circ} & \!\!= & \frac{G_{-
			\circ}}{K_{- \circ}} + \left[ \frac{G_{+ +}}{K_{- \circ}} \right]_{- \circ}
	\!\!+ \left[ \frac{G_{+ -}}{K_{- \circ}} \right]_{- \circ} \!\!. 
	\label{eq:firstsplitinter}
	\end{eqnarray}
	Because of the simplicity of $G_{+ +}$ (see (\ref{eq:exactfunctionG++})), the
	sum-split of $G_{+ +} / K_{- \circ}$ can be \COM{achieved explicitly} via the pole removal
	technique:
	\begin{eqnarray*}
		\left[ \frac{G_{+ +}}{K_{- \circ}} \right]_{+ \circ} \!\!= \frac{G_{+ +}}{K_{- \circ} (a_1, \alpha_2)} & \text{ and } & \left[
		\frac{G_{+ +}}{K_{- \circ}} \right]_{- \circ} \!\!= G_{+ +} 
		\left( \frac{1}{K_{- \circ}} - \frac{1}{K_{- \circ}
			(a_1, \alpha_2)} \right) .
	\end{eqnarray*}
	Now, by construction, the left-hand side (LHS) of (\ref{eq:firstsplitinter}) is analytic in $\mathcal{D}_{+
		\circ}$, while the right-hand side (RHS) of (\ref{eq:firstsplitinter}) is analytic in
	$\mathcal{D}_{- \circ}$. \COM{Hence it is possible to use (\ref{eq:firstsplitinter}) to construct a function $E_{1\circ}$ that
		is analytic on $\mathbb{C} \times \mathcal{A}_2$ and defined by}
	\begin{eqnarray}
	\label{eq:EalphaLiouville1}
	E_{1\circ} & = & \left\{ \begin{array}{rc}
	F_{+ +} K_{+ \circ} - \frac{G_{+ +}}{K_{-
			\circ} (a_1, \alpha_2)} - \left[ \frac{G_{+ -}}{K_{- \circ}} \right]_{+
		\circ} & \!\!\!\!\!\!\text{ if }\tmmathbf{\alpha} \in \mathcal{D}_{+ \circ},\\
	\frac{G_{- \circ}}{K_{- \circ}} + G_{+ +} 
	\left( \frac{1}{K_{- \circ}} - \frac{1}{K_{- \circ}
		(a_1, \alpha_2)} \right) + \left[ \frac{G_{+ -}}{K_{- \circ}} \right]_{-
		\circ} & \!\!\!\!\!\!\text{ if }\tmmathbf{\alpha} \in \mathcal{D}_{- \circ}\COM{.} 
	\end{array} \right.  
	\end{eqnarray}
	\COM{Moreover, it}
	can be shown that \RED{it tends to zero as $|\alpha_1|\rightarrow \infty$} (see \RED{Appendix} \ref{app:Liouville1}), and \COM{so} we can apply Liouville's theorem \COM{in the $\alpha_1$ plane} to get $E_{1\circ}\equiv 0$; hence
	\begin{eqnarray}
	\small{F_{+ +} K_{+ \circ} - \frac{G_{+ +}}{K_{- \circ}
			(a_1, \alpha_2)} - \left[ \frac{G_{+ -}}{K_{- \circ}} \right]_{+ \circ}}\!\! & =
	& 0, \label{eq:split1eq1}\\
	\small{\frac{G_{- \circ}}{K_{- \circ}} + G_{+ +} \left(
		\frac{1}{K_{- \circ}} - \frac{1}{K_{- \circ} (a_1,
			\alpha_2)} \right) + \left[ \frac{G_{+ -}}{K_{- \circ}} \right]_{- \circ}}\!\! &
	= & 0.  \label{eq:split1eq2}
	\end{eqnarray}
	
	\subsection{A second split\COM{,} in the $\alpha_2$
		plane}\label{sub:secondsplitgeneric}
	
	Upon multiplying equation (\ref{eq:split1eq1}) by $K_{- +} (a_1, \alpha_2) / K_{+ -}$, it becomes
	\begin{eqnarray}
	F_{+ +} K_{+ +} K_{- +} (a_1,
	\alpha_2) & = & \frac{G_{+ +}}{K_{- -} (a_1, \alpha_2)
		K_{+ -}} + \frac{K_{- +} (a_1, \alpha_2)}{K_{+ -}} \left[ \frac{G_{+ -}}{K_{- \circ}} \right]_{+ \circ}\!\! . 
	\label{eq:secondsplitinter1}
	\end{eqnarray}
	The LHS is a $+ +$ function, and, once again, formally, \COM{using Lemma \ref{th:Cauchysumdecomp}}, each of the two terms
	on the RHS of (\ref{eq:secondsplitinter1}) ha\COM{s} a sum-split decomposition in
	the $\alpha_2$ plane \COM{(the associated operators being denoted $\left[ \ \right]_{\circ -}$ and $\left[
		\ \right]_{\circ +}$)}, such that we can rewrite (\ref{eq:secondsplitinter1}) as
	\begin{eqnarray}
	F_{+ +} K_{+ +} K_{- +} (a_1,
	\alpha_2) - \left[ \frac{G_{+ +}}{K_{- -} (a_1,
		\alpha_2) K_{+ -}} \right]_{\COM{\circ} +} \!\!- \left[ \frac{K_{- +}
		(a_1, \alpha_2)}{K_{+ -}} \left[ \frac{G_{+ -}}{K_{-
			\circ}} \right]_{+ \circ} \right]_{\COM{\circ} +}\!\!\!\!\!\! &&  \label{eq:secondsplitinter2}\\
	= \left[ \frac{G_{+ +}}{K_{- -} (a_1, \alpha_2) K_{+
			-}} \right]_{\COM{\circ} -} \!\!+ \left[ \frac{K_{- +} (a_1,
		\alpha_2)}{K_{+ -}} \left[ \frac{G_{+ -}}{K_{- \circ}}
	\right]_{+ \circ} \right]_{\COM{\circ} -}\!\!. &&  \nonumber
	\end{eqnarray}
	Again, because of the form of $G_{+ +}$, the related split can be
	performed explicitly by pole removal to get
	\begin{eqnarray*}
		\left[ \frac{G_{+ +}}{K_{- -} (a_1, \alpha_2) K_{+ -}
		} \right]_{\COM{\circ} -} &\!\! = & G_{+ +}
		\left( \frac{1}{K_{- -} (a_1, \alpha_2) K_{+ -}} -
		\frac{1}{K_{- -} (a_1, a_2) K_{+ -} (\alpha_1, a_2)} \right),\\
		\left[ \frac{G_{+ +}}{K_{- -} (a_1, \alpha_2) K_{+ -}
		} \right]_{\COM{\circ} +} &\!\! = & \small{\frac{G_{+ +}}{K_{- -} (a_1, a_2) K_{+ -} (\alpha_1, a_2)}} \text{ } .
	\end{eqnarray*}
	Now\COM{, by inspection, it is clear that} the LHS of (\ref{eq:secondsplitinter2}) is analytic on $\mathcal{D}_{+
		+}$, while its RHS is analytic on $\mathcal{D}_{+ -}$. Hence, it is possible
	to \COM{construct} a function $E_{+ 2}$ that is analytic on $\tmop{UHP}_1 \times
	\mathbb{C}$ and defined by
	\scriptsize
	\begin{eqnarray}
	E_{+ 2} & = & \!\! \left\{ \begin{array}{rc}
	\tiny{\!\!\!\! F_{+ +}  K_{+ +} K_{- +}
		(a_1, \alpha_2) - \frac{G_{+ +}}{K_{- -} (a_1, a_2)
			K_{+ -} (\alpha_1, a_2)} - \left[ \frac{K_{- +} (a_1, \alpha_2)}{K_{+ -}} \left[ \frac{G_{+ -}}{K_{- \circ}} \right]_{+ \circ}
		\right]_{\COM{\circ} +}} & \!\!\!\!\!\!\!\text{if } \tmmathbf{\alpha} \in \mathcal{D}_{+ +},\\
	\tiny{\!\!\!\! G_{+ +} \left( \frac{1}{K_{- -} (a_1, \alpha_2)
			K_{+ -}} - \frac{1}{K_{- -} (a_1, a_2) K_{+ -}
			(\alpha_1, a_2)} \right) + \left[ \frac{K_{- +} (a_1, \alpha_2)}{K_{+ -}
		} \left[ \frac{G_{+ -}}{K_{- \circ}} \right]_{+ \circ}
		\right]_{\COM{\circ} -}} & \!\!\!\!\!\!\!\text{if } \tmmathbf{\alpha} \in \mathcal{D}_{+ -}\COM{.}
	\end{array} \right.  \label{eq:Eplus2alphaLiouville}
	\end{eqnarray}
	\normalsize
	\COM{One of the aims of this work is to provide a \textit{constructive path} towards Radlow's ansatz. In order to do so, we wish to apply Liouville's theorem in the $\alpha_2$ plane and, for this, we need to examine the right-hand sides of (\ref{eq:Eplus2alphaLiouville}) as $|\alpha_2|\rightarrow\infty$ in their respective half-planes of analyticity. Using the proof given in Appendix \ref{app:Liouville2}, we can show that $E_{+2}\equiv0$ and hence obtain the two main equations of the paper:}
	\begin{align}
	F_{+ +} & = \small{\frac{G_{+ +}}{K_{+ +} K_{- +} (a_1, \alpha_2)
			K_{- -} (a_1, a_2) K_{+ -} (\alpha_1, a_2)}} \label{eq:genericWHeq1}\\
	& + \small{\frac{1}{K_{+ +} K_{- +} (a_1, \alpha_2)}
		\left[ \frac{K_{- +} (a_1, \alpha_2)}{K_{+ -}} \left[
		\frac{G_{+ -}}{K_{- \circ}} \right]_{+ \circ} \right]_{\COM{\circ} +}}\!\!,
	\nonumber\\
	0 & = G_{+ +} \left( \frac{1}{K_{- -}
		(a_1, \alpha_2) K_{+ -}} - \frac{1}{K_{- -} (a_1, a_2)
		K_{+ -} (\alpha_1, a_2)} \right) \label{eq:genericWHeq2}\\
	& + \small{\left[ \frac{K_{- +} (a_1, \alpha_2)}{K_{+ -}} \left[ \frac{G_{+ -}}{K_{- \circ}} \right]_{+ \circ}
		\right]_{\COM{\circ} -}}\!\!. \nonumber 
	\end{align}
	
	Remember that in order to recover \COM{the} physical field everywhere via
	(\ref{eq:explicituFK}), the unknown of interest is the function $F_{+ +}
	(\tmmathbf{\alpha})$. We can at this stage make two important remarks
	regarding (\ref{eq:genericWHeq1}). Firstly, provided that we know the function
	$G_{+ -}$, then $F_{+ +}$ can in theory be recovered. Secondly, it is
	important to note that the first term on the RHS of
	(\ref{eq:genericWHeq1}) is exactly Radlow's ansatz published in
	{\cite{Radlow1965}}. The main issue with Radlow's solution was that the
	resulting physical field did not behave as expected near the tip of the
	quarter-plane (Radlow's ansatz predicts a behaviour of $\mathcal{O} (r^{1 /
		4})$, while the correct behaviour is $\mathcal{O} (r^{\nu_1 - 1 / 2})$, where
	$\nu_1$ is related to the first eigenvalue of the Laplace-Beltrami operator).
	As such the benefit of this equation is dual. On \COM{the} one hand, it is clear that
	(\ref{eq:genericWHeq1}) \COM{indicates} the error in Radlow's analysis, since a term
	is missing from his ansatz. On the other hand, we provide here a constructive procedure
	showing how this ansatz is obtained, which can be enlightening in view of
	the fact that no derivation was provided in Radlow's original work. Indeed, it was the fact that Radlow merely stated a solution in \cite{Radlow1965}, that has partially led to difficulties in establishing and quantifying the error to-date.
	
	In addition, we also know that the correct physical behaviour of the solution
	should be enforced by the term involving $G_{+ -}$. Equation
	(\ref{eq:genericWHeq2}), that we will refer to as a {\tmem{compatibility
			equation}}, is very interesting in that respect. Firstly, it does not appear in
	Radlow's work, nor in any subsequent work to our knowledge. Secondly, if it
	can somehow be inverted (which \COM{in practice} is a very difficult thing to do), it provides a
	way to obtain $G_{+ -}$. \COM{Even though} it is not possible to do this exactly (as
	the authors believe is the case), it provides a way of testing any approximation to $G_{+
		-}$. Hence, we believe that the compatibility equation (\ref{eq:genericWHeq2})
	is key to solving the problem at hand. We will not go through this route in
	this paper, but it will be the basis of a future article.
	
	Before going further, note also that (\ref{eq:split1eq2}) has not been used so
	far. It is possible to employ it to obtain two more equations involving $G_{-
		-}$, $G_{+ -}$ and $G_{- +}$, by introducing similarly a function $E_{- 2}$
	entire in the $\alpha_2$ complex plane (which is again zero by application of Liouville's theorem). However, we do not believe that these
	will provide further information on the solution, and so are extraneous. Moreover, nowhere in this section did we use the
	definition of $\mathcal{A}_{1, 2}$ explicitly; hence the results obtained
	remain valid when $a_1$ or $a_2$ are positive.
	
	To summarise, in order to solve our problem and find $F_{+ +}$, we need to
	gain some information about $G_{+ -}$ and find an approximation that will be
	compatible both with the physics of the problem and with the compatibility
	equation (\ref{eq:genericWHeq2}). \COM{A possible approximation scheme for $G_{+-}$, involving and explicit canonical integral, is suggested in \cite{interestingintegral}.}
	\COM{ However, f}or the purpose of this paper, let us assume that we know $F_{+ +}$ and let us
	try to find out what can be inferred about the diffraction coefficient.
	
	\subsection{Link with diffraction coefficient}
	
	Classically, (see e.g. {\cite{Assier2012,Assier2012b}}) the
	Dirichlet corner diffraction coefficient $f^d (\theta, \varphi, \theta_0,
	\varphi_0)$ is defined by
	\begin{eqnarray}
	u_{\tmop{sph}} & \underset{kr \rightarrow \infty}{\approx} & 2 \pi
	\frac{e^{ikr}}{kr} f^d (\theta, \varphi ; \theta_0, \varphi_0), 
	\label{eq:originaldiffcoeffart}
	\end{eqnarray}
	where $u_{\tmop{sph}}$ represents the spherical wave emanating from the tip.
	Assuming that $F_{+ +}$ is known, using complexified spherical coordinates, one can apply a double steepest-descent analysis as $kr\rightarrow\infty$ \cite{Bleistein2012,Abrahamssteepest1994} to obtain the following relationship between the diffraction coefficient
	and $F_{+ +}$:
	\begin{eqnarray}
	f^d (\theta, \varphi ; \theta_0, \varphi_0) & = & \frac{kF_{+ +} (- k \cos
		(\varphi) \sin (\theta), - k \sin (\varphi) \sin (\theta) )}{4 \pi^2 i} \cdot
	\label{eq:reltwodiffcoeffart}
	\end{eqnarray}
	We believe that this formula should remain valid everywhere. We may of course get other far-field contributions (edge-diffracted
	waves, reflected wave, etc.) that will result from crossing poles when deforming the
	various contours to their steepest-descent paths. However, the $1 / kr$
	component can only be the one given in (\ref{eq:reltwodiffcoeffart}). In
	particular, it should have the same singular regions as those obtained
	(explicitly) with the embedding procedure, but most importantly, this
	formula should be valid in the regions that the embedding formulae cannot (yet)
	reach. We can easily observe that the polar singularity structure is similar. In fact we
	have seen in {\cite{Assier2012}} that if we write $\xi = \cos (\varphi) \sin
	(\theta)$, $\xi_0 = \cos (\varphi_0) \sin (\theta_0)$, $\eta = \sin (\varphi)
	\sin (\theta)$ and $\eta_0 = \sin (\varphi_0) \sin (\theta_0)$, the
	diffraction coefficient had simple poles when $\xi = - \xi_0$ and $\eta = -
	\eta_0$. Upon noticing that in (\ref{eq:reltwodiffcoeffart}) we evaluate $F_{++}$ at $(\alpha_1,\alpha_2)=(-k\xi,-k\eta)$, realising that $(a_1,a_2)=(k\xi_0,k\eta_0)$, and remembering that $F_{++}$ has poles at $\alpha_{1,2}=a_{1,2}$, we recover the expected polar singularities.
	
	\RED{Note\footnote{\COM{The authors are most grateful to the anonymous reviewer for this and many other suggestions, which have significantly improved the clarity and focus of the article.}} that (\ref{eq:reltwodiffcoeffart}) implies that the diffraction coefficient does not depend on $k$. To see this, let $v(\boldsymbol{x})$ be the scattered field of the Dirichlet quarter-plane problem for $k=1$. One can show directly that, for $k>0$, the solution $u$ of our problem summarised in Section \ref{sec:formulation-summary} is given by $u(\boldsymbol{x})=v(k\boldsymbol{x})$. Using the basic definition of the double Fourier transform, and the fact that $\tfrac{\partial u}{\partial x_3}(x_1,x_2,0^+)=0$ on $Q_2 \cup Q_3 \cup Q_4$, we can show that  $ik F_{+ +} (k\tmmathbf{\alpha}) =\mathfrak{F} [ \tfrac{\partial
			v}{\partial x_3} (x_1, x_2, 0^+) ] (\tmmathbf{\alpha})$, which is clearly independent of $k$.}
	
	
	\RED{Another} interesting feature to be considered is that we know
	{\cite{Assier2012}} that the diffraction coefficient should in fact be
	purely imaginary (at least where the MSF are valid). However, it is not
	obvious that the RHS of (\ref{eq:reltwodiffcoeffart}) is indeed purely imaginary.
	
	
	One issue with the formula (\ref{eq:reltwodiffcoeffart}) is that the function
	$F_{+ +}$ is evaluated on the real interval $(- k, k)$ in both complex planes.
	However, it is clear from the above analysis that the segment $(- k, 0)$ does
	not lie in $\tmop{UHP}_1$ or $\tmop{UHP}_2$. Hence, we are forced to evaluate a
	$+ +$ function {\tmem{outside}} $\mathcal{D}_{+ +}$. This problem can be dealt with by means of analytically continuing $F_{+ +}$ within that
	region.
	
	\section{Comparison between Radlow's ansatz and MSF}\label{sec:results}
	
	In this section, we compare the diffraction coefficient obtained by the MSF to
	that obtained by using Radlow's {\tmem{erroneous}} ansatz. The MSF is now an
	established method known to be correct within a certain domain of the observer
	space. The idea of comparing both method\COM{s} is mainly due to serendipity. Whilst
	testing a method to evaluate the effect of $G_{+ -}$ on the diffraction
	coefficient, we once accidentally set $G_{+ -} = 0$, which is equivalent to
	using Radlow's ansatz exactly. To our surprise, this led to a \RED{very good} 
	agreement with the MSF results, where these formulae were valid. We decided to
	explore the incidence space, and so far we could not find any incident angle
	leading to \RED{an obvious disagreement} 
	between the two methods. Here we present four
	distinct incidence\COM{s} (we keep $\theta_0 = \pi / 4$, and choose four different
	$\varphi_0$), corresponding to different signs for $a_{1, 2}$. The chosen
	incidence\COM{s} are summarised in Figure \ref{fig:incidencediagram}.
	
	\begin{figure}[htbp]
		\centering
		\includegraphics[width=0.75\textwidth]{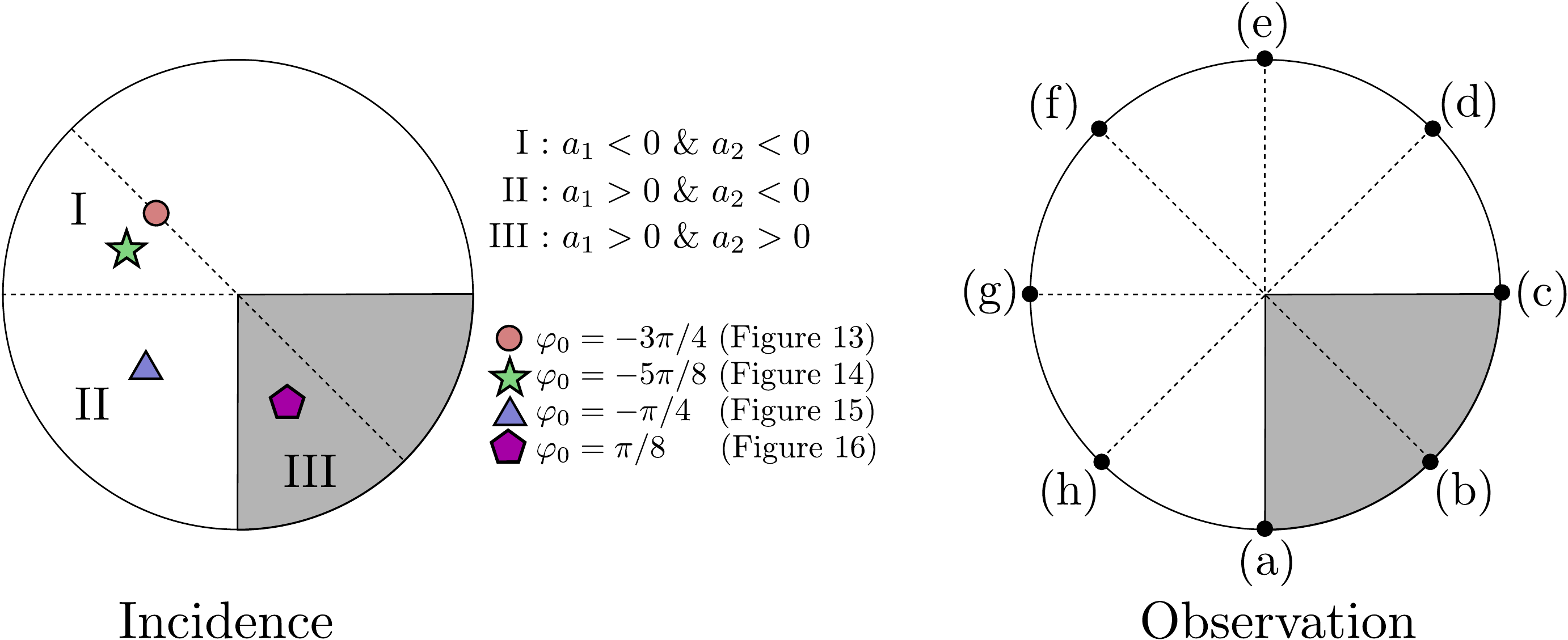}
		\caption{Left: Illustration of the incident angles used in the presentation of the
			results. We have ensured that each region corresponding to a different sign
			combination of $a_1$ and $a_2$ was considered. Right: Illustration of the 8 arcs of observation used in the presentation of the
			results.}
		\label{fig:incidencediagram}
	\end{figure}
	
	For each incidence, we pick 8 arcs surrounding the quarter-plane on which we
	evaluate the diffraction coefficient. That is we pick 8 values of $\varphi$
	between 0 and $2 \pi$, and for each value of $\varphi$, we evaluate the
	coefficient for $\theta \in [0, \pi / 2]$. The results are presented in
	Figures \ref{fig:run0}--\ref{fig:run4}, showing \RED{very good} 
	agreement between
	the two methods. 
	When the diffraction coefficient does not have any singularities, as in Figures \ref{fig:run0}(e)(f)(g), \ref{fig:run2}(e)(f)(g), \ref{fig:run3}(a)(g)(h) and \ref{fig:run4}(b), it means that the only far-field component in the observation region is the spherical wave emanating from the tip, this is the so-called oasis zone. The diffraction coefficient becomes singular at the boundaries of existence of the edge-diffracted fields. Another important point to mention is the validity of this ansatz in the
	region where the MSF are not valid (see \cite{Assier2012} for discussion) due to double diffraction of the field (Figures \ref{fig:run3}(c) and \ref{fig:run4}(a)(c)). Passing the limit of validity, we notice
	that the diffraction coefficient given by Radlow's ansatz, which is purely imaginary everywhere else,
	becomes purely real. Mathematically this corresponds to saddle points
	going through a branch point during the steepest-descent procedure. Having no
	data to compare to in this region, it remains to be seen if this yields the correct physical solution.

\newcommand\mywidth{0.35}
\begin{figure}[htbp]
\centering
(a)\raisebox{-.1\height}{\includegraphics[width=\mywidth\textwidth]{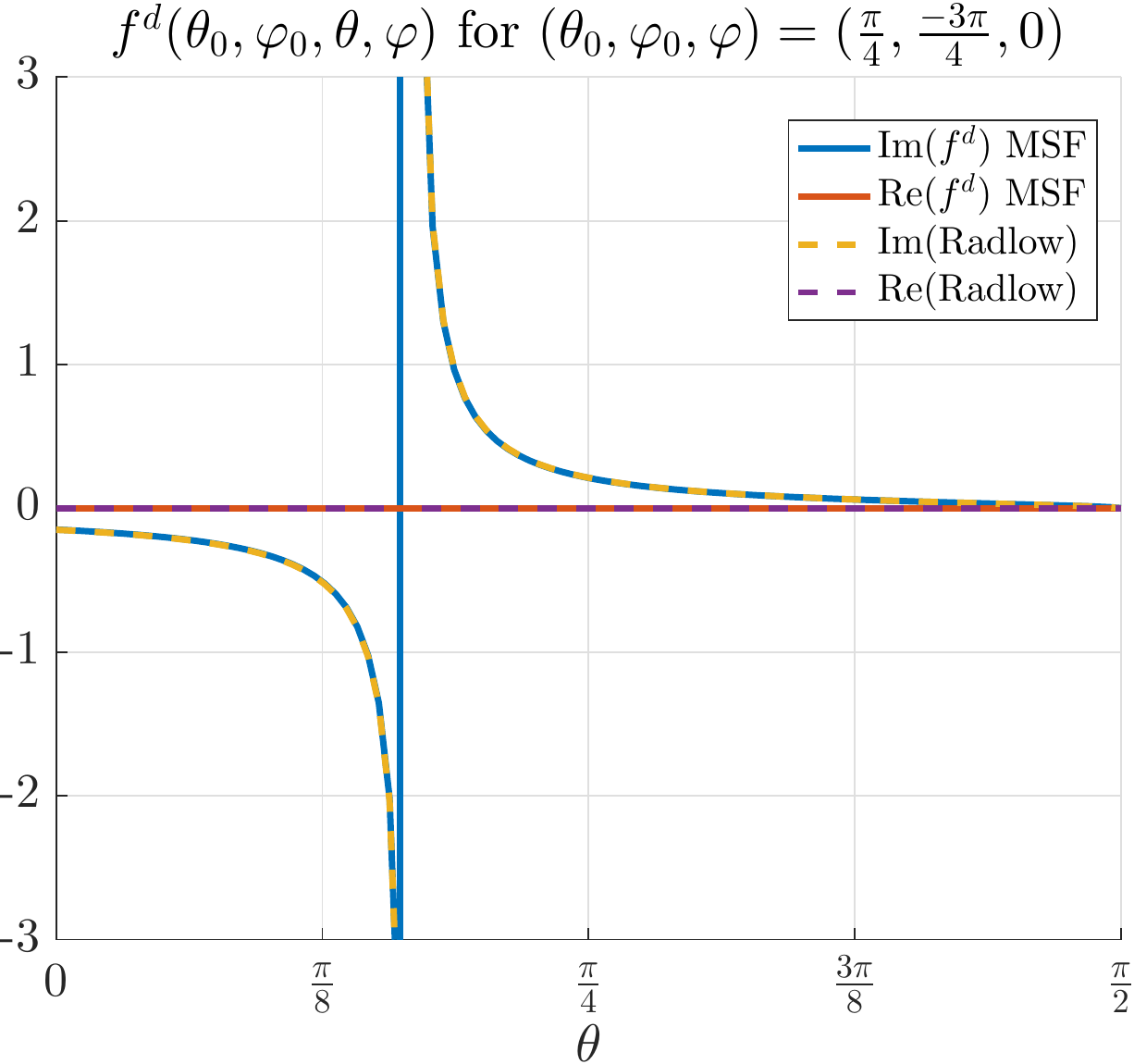}}\quad(b)\raisebox{-.1\height}{\includegraphics[width=\mywidth\textwidth]{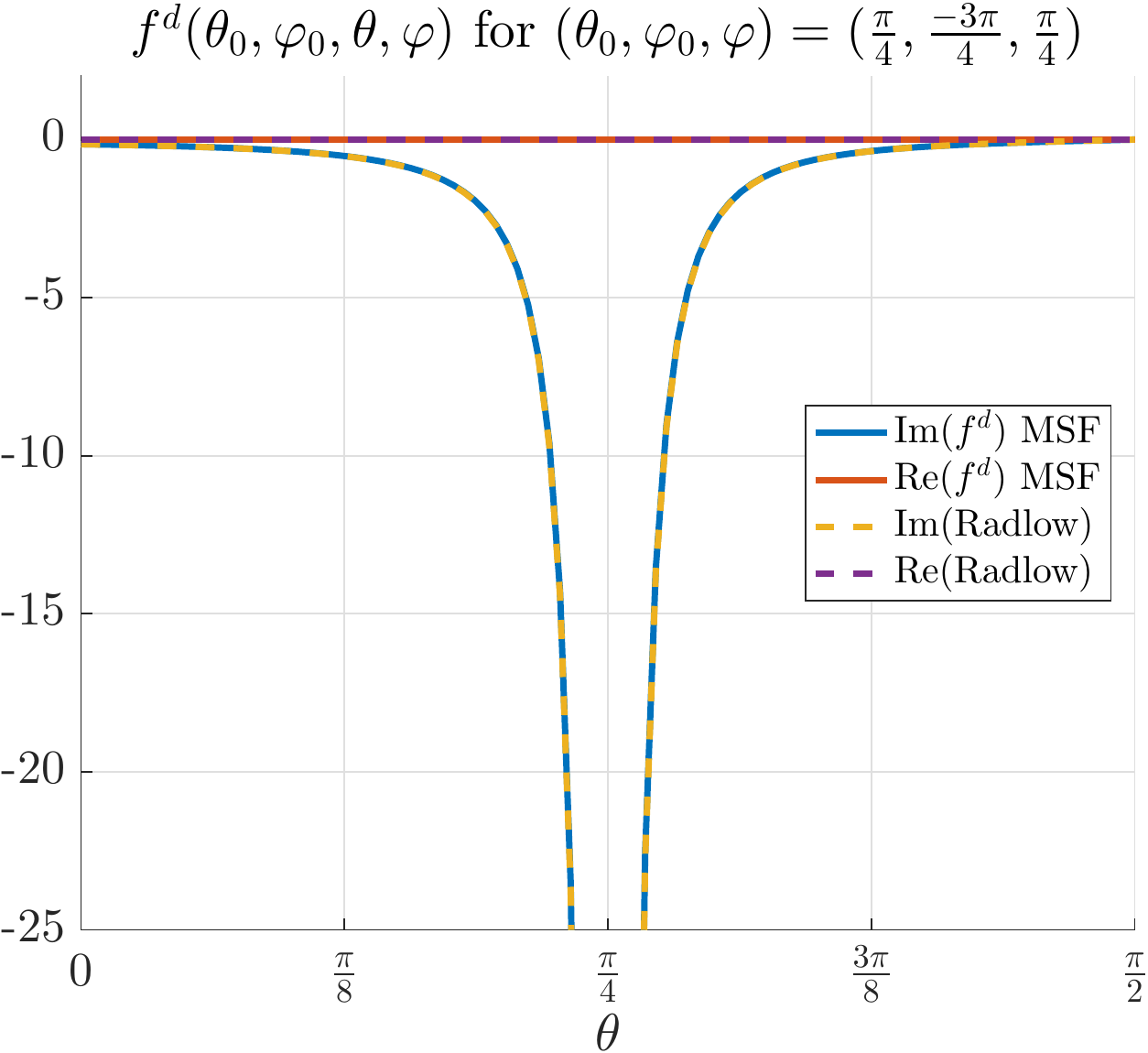}}
  
(c)\raisebox{-.1\height}{\includegraphics[width=\mywidth\textwidth]{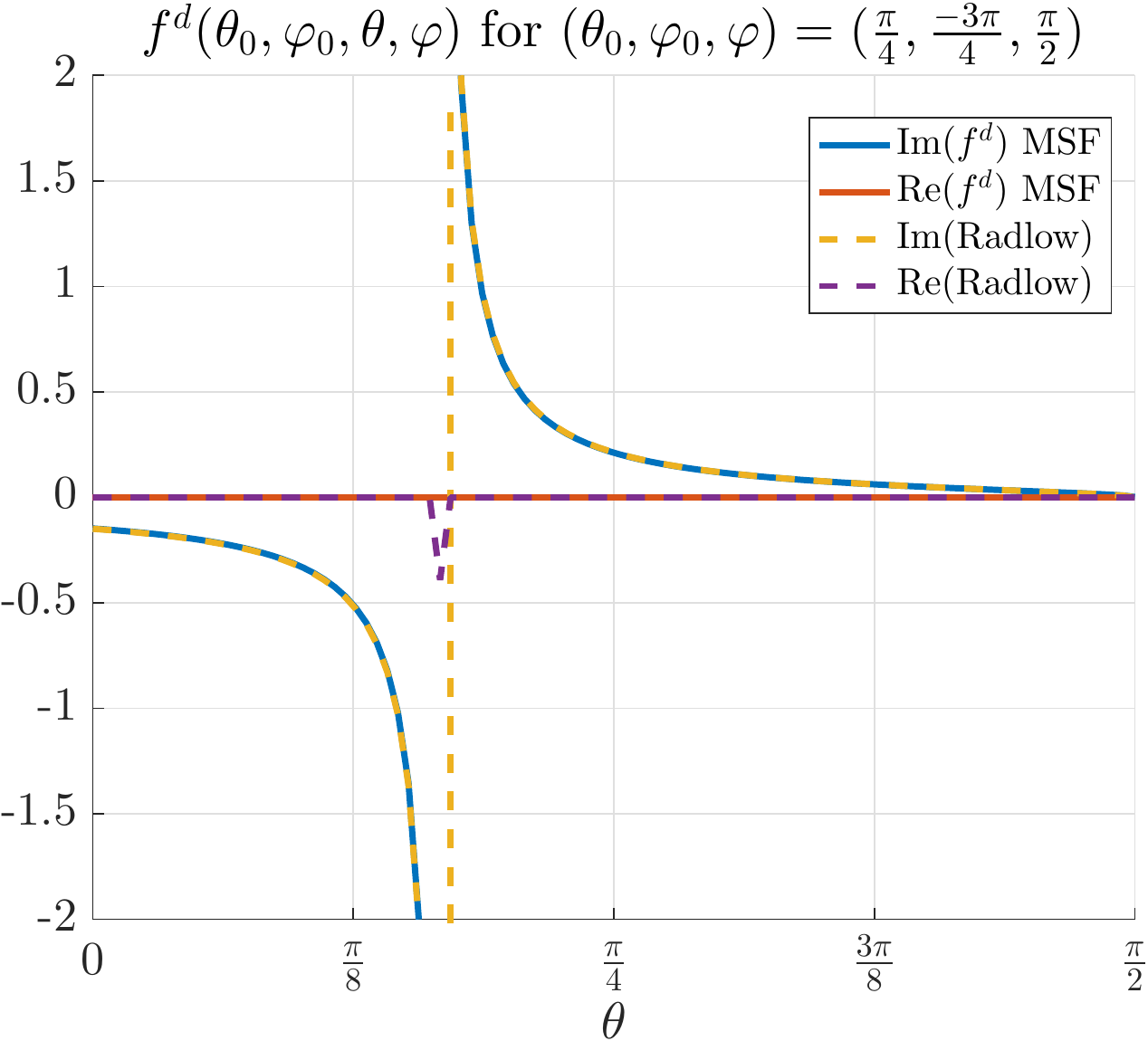}}\quad(d)\raisebox{-.1\height}{\includegraphics[width=\mywidth\textwidth]{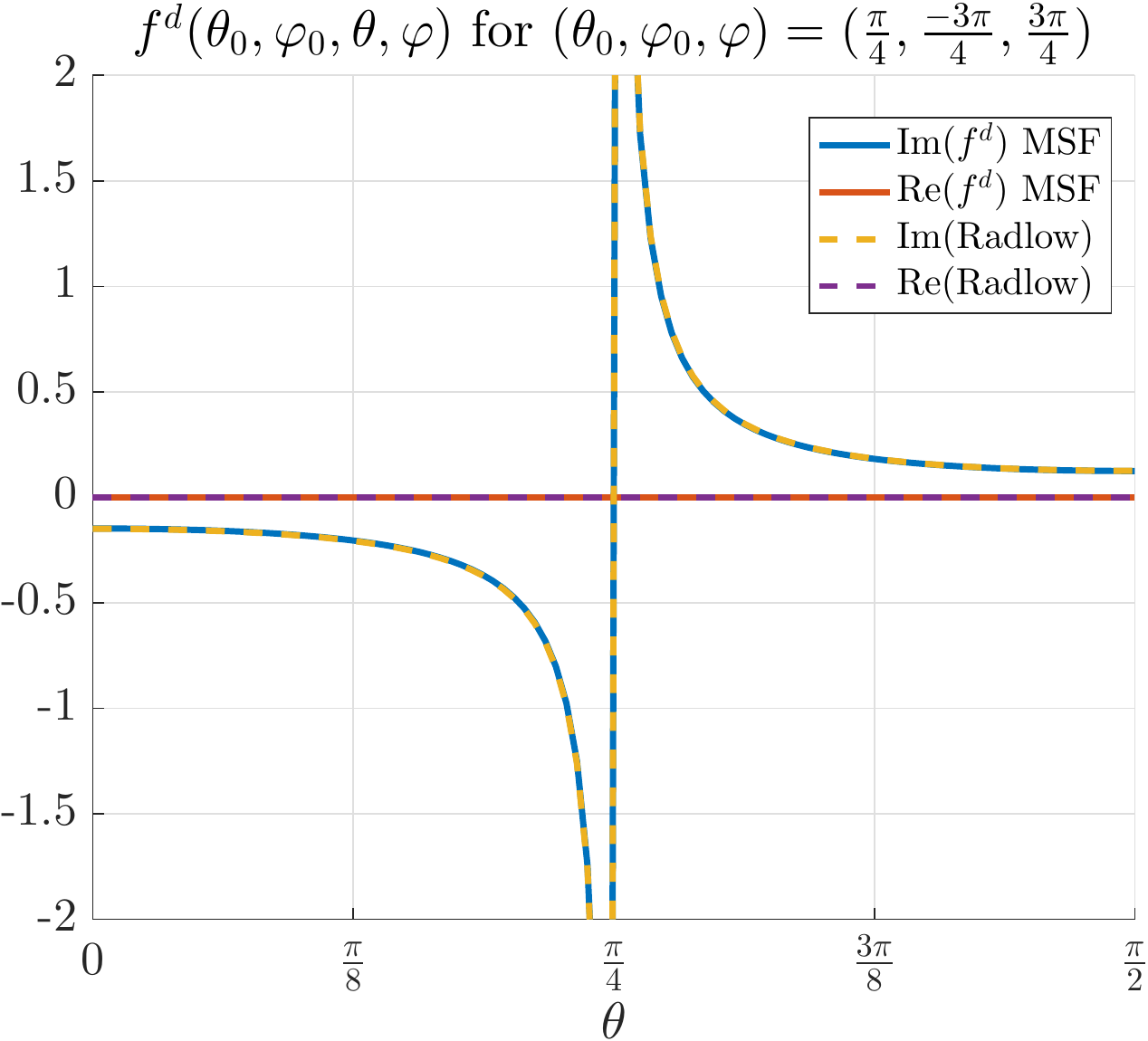}}

(e)\raisebox{-.1\height}{\includegraphics[width=\mywidth\textwidth]{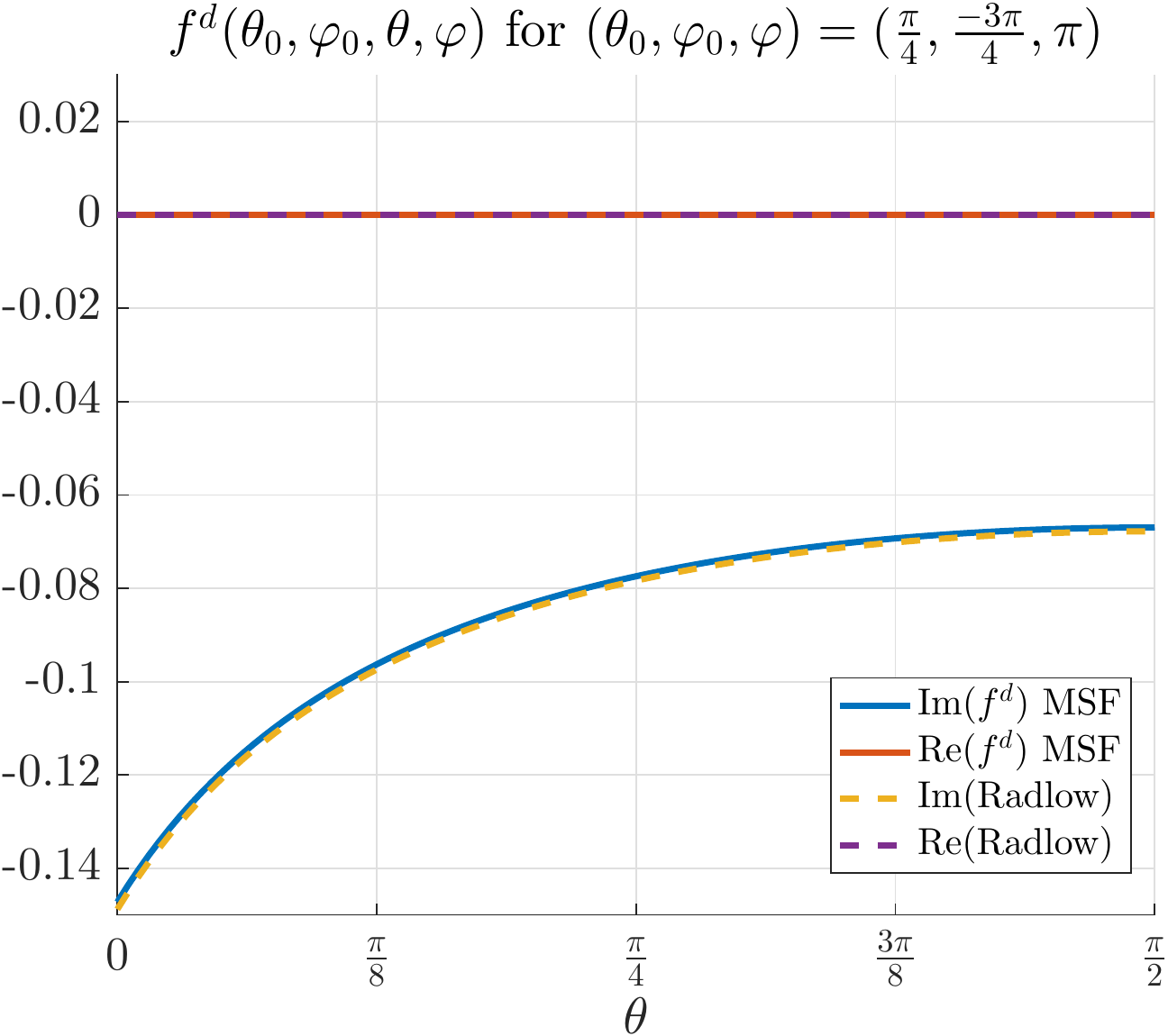}}\quad(f)\raisebox{-.1\height}{\includegraphics[width=\mywidth\textwidth]{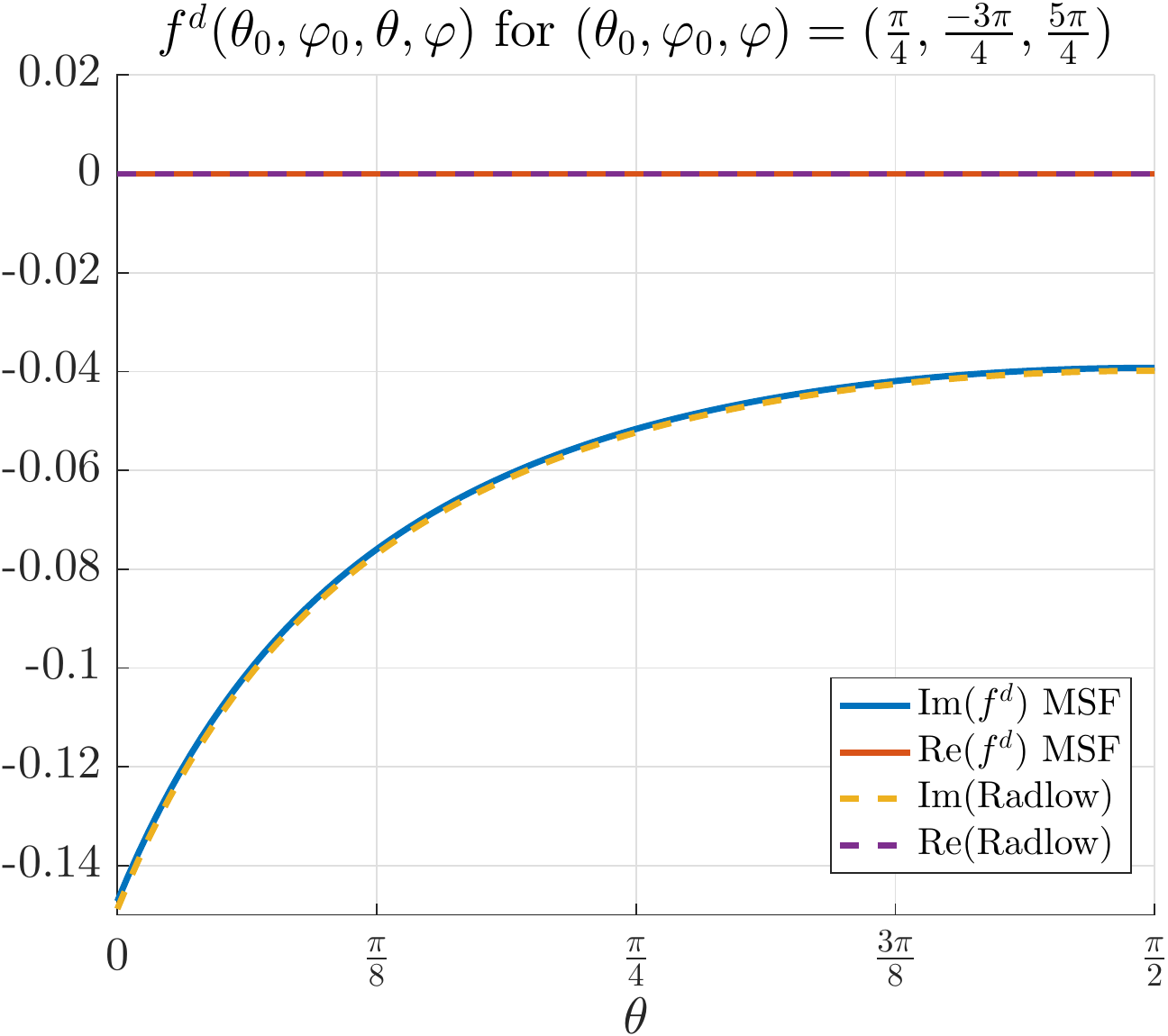}}
  
(g)\raisebox{-.1\height}{\includegraphics[width=\mywidth\textwidth]{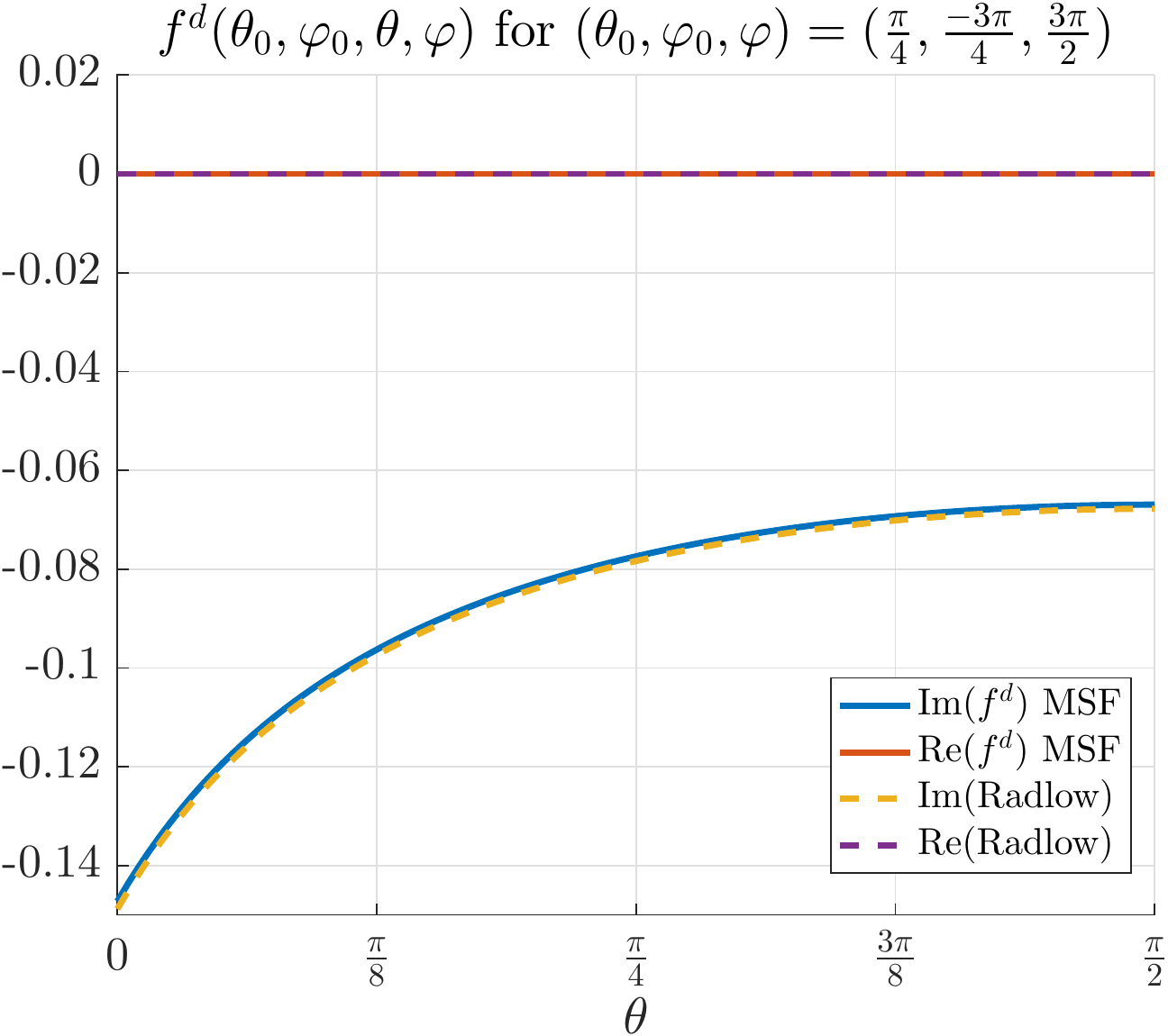}}\quad(h)\raisebox{-.1\height}{\includegraphics[width=\mywidth\textwidth]{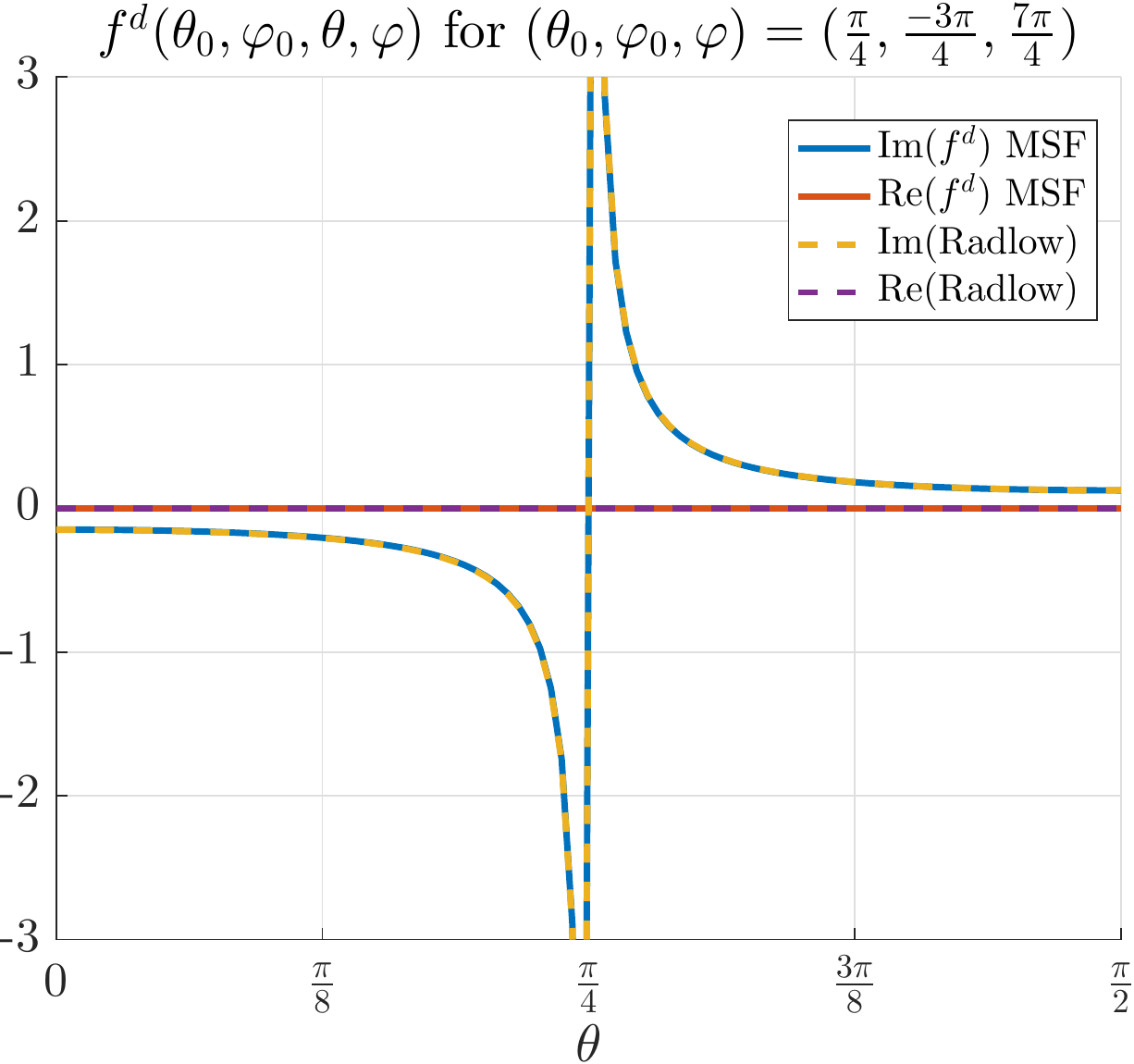}}
  \caption{Diffraction coefficient for incidence $(\theta_0, \varphi_0) =
  \left( \frac{\pi}{4}, \frac{- 3 \pi}{4} \right)$, i.e. we have $a_1 < 0$ and
  $a_2 < 0$, with polar observation angle $\theta \in \left[ 0, \frac{\pi}{2}
  \right]$ and various values of the azimuthal observation angles $\varphi = 0,
  \frac{\pi}{4}, \frac{\pi}{2}, \frac{3 \pi}{4},\pi,
  \frac{5 \pi}{4}, \frac{3 \pi}{2}, \frac{7 \pi}{4}$ (from (a) to (h)).}
\label{fig:run0}
\end{figure}

  

\begin{figure}[htbp]
\centering
 (a)\raisebox{-.1\height}{\includegraphics[width=\mywidth\textwidth]{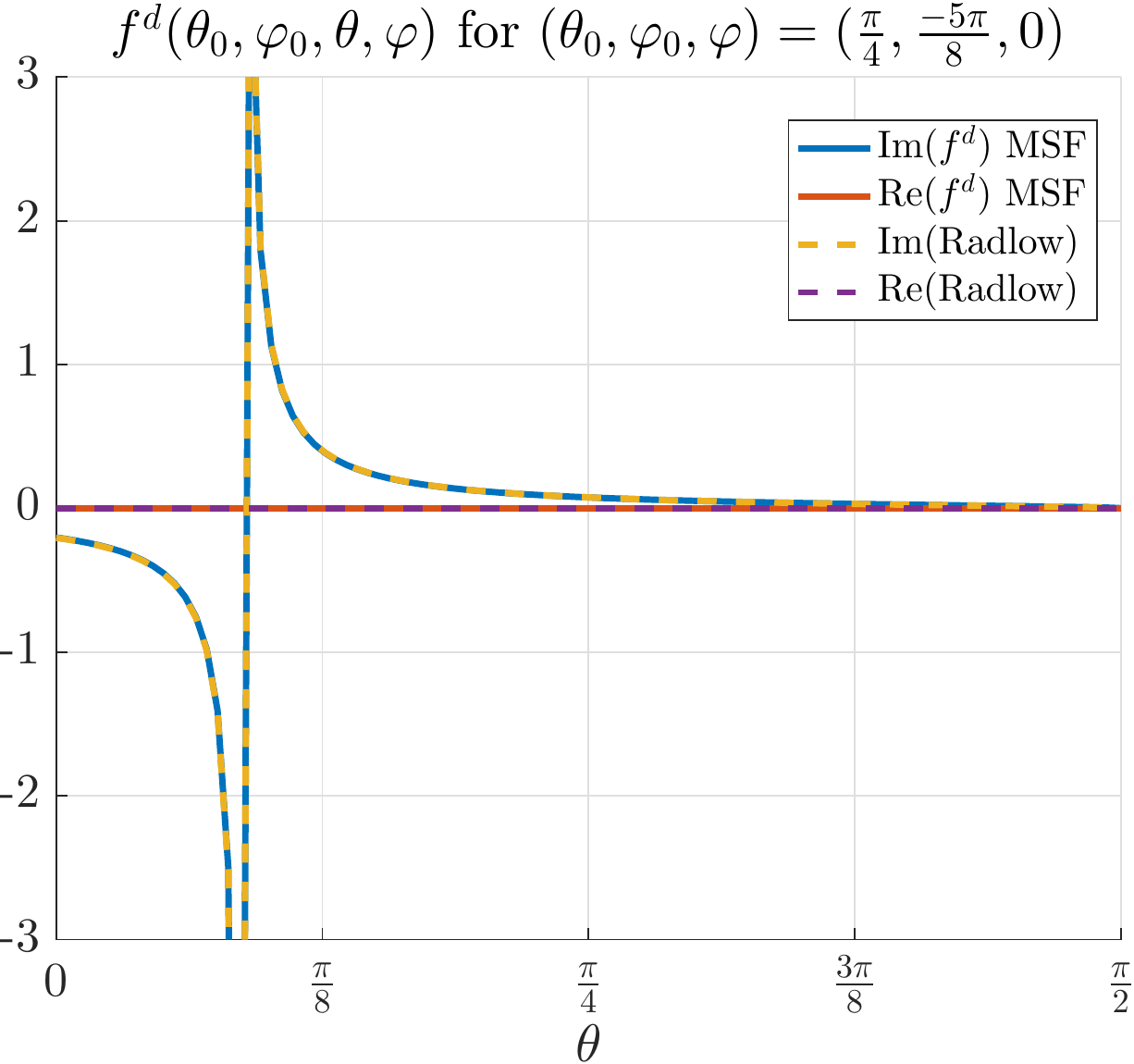}}\quad(b)\raisebox{-.1\height}{\includegraphics[width=\mywidth\textwidth]{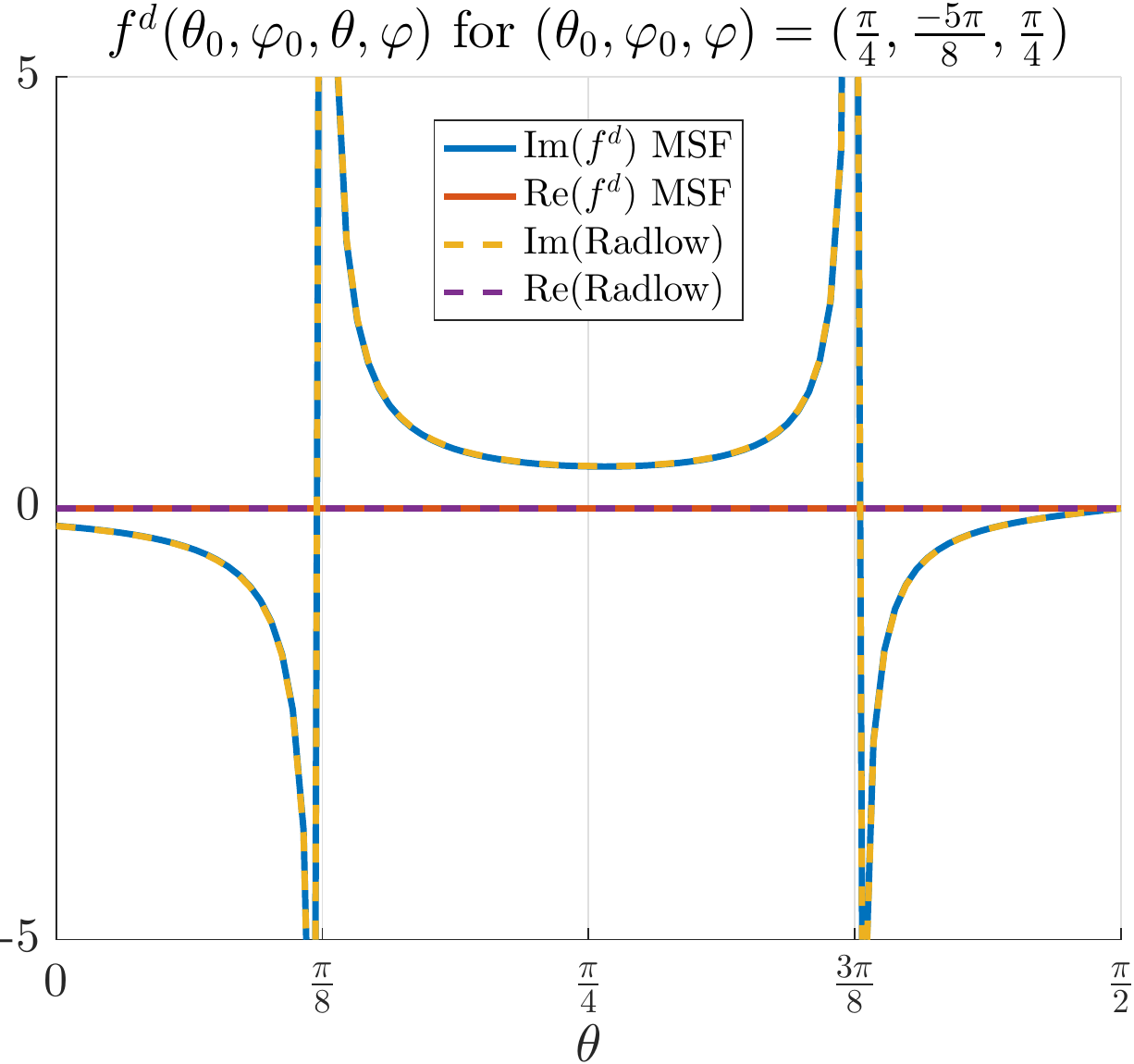}}
 (c)\raisebox{-.1\height}{\includegraphics[width=\mywidth\textwidth]{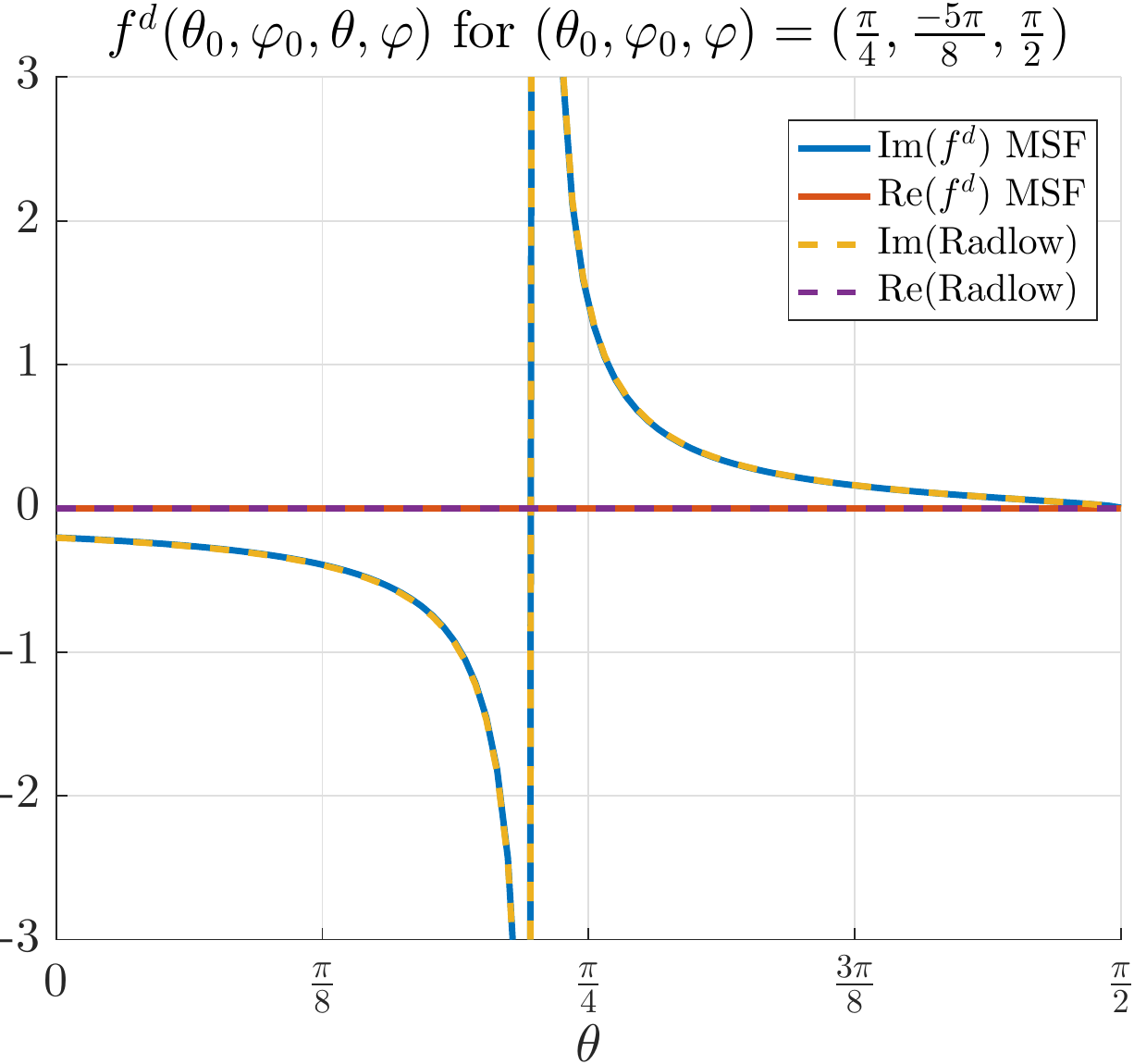}}\quad(d)\raisebox{-.1\height}{\includegraphics[width=\mywidth\textwidth]{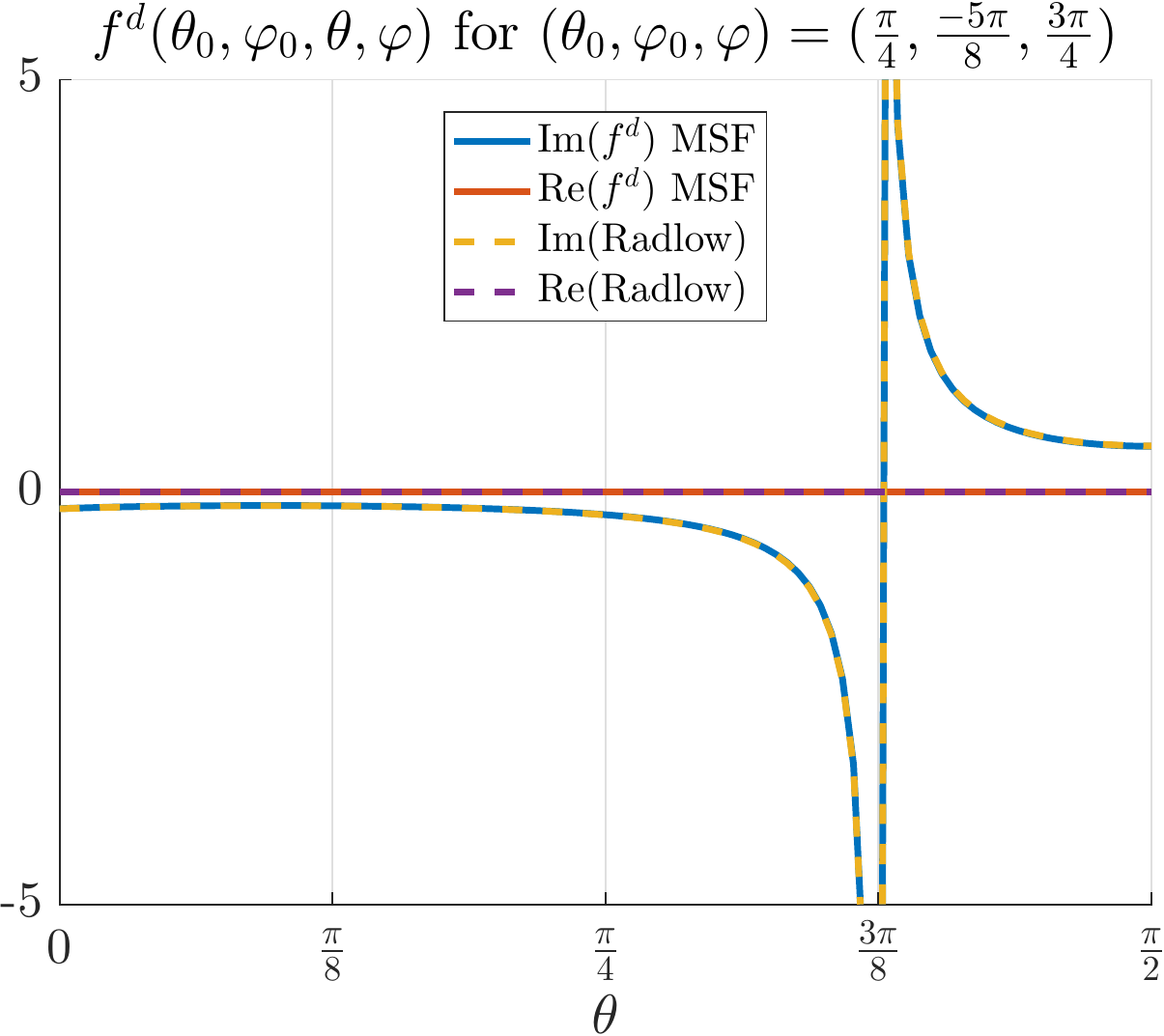}}
 (e)\raisebox{-.1\height}{\includegraphics[width=\mywidth\textwidth]{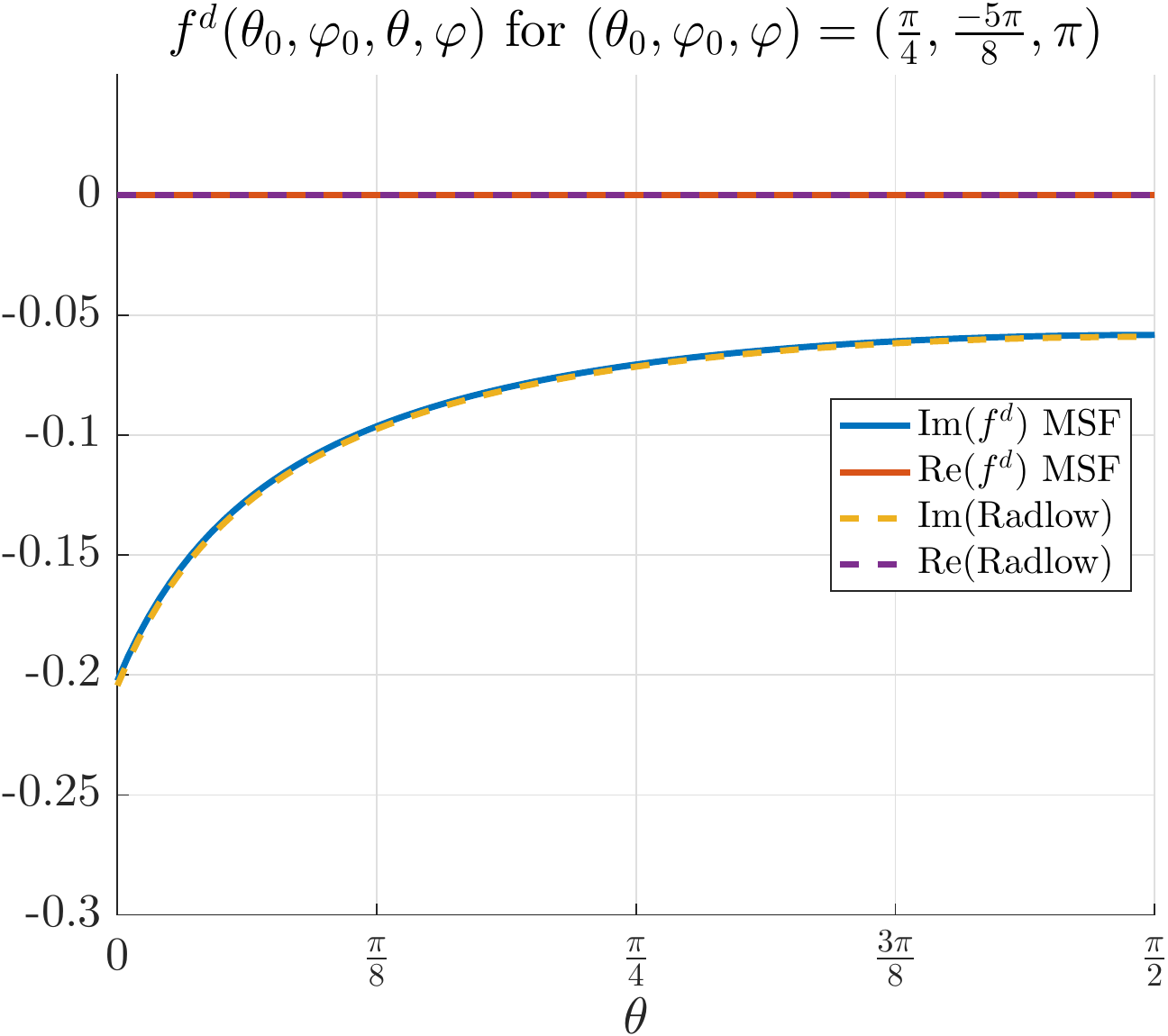}}\quad(f)\raisebox{-.1\height}{\includegraphics[width=\mywidth\textwidth]{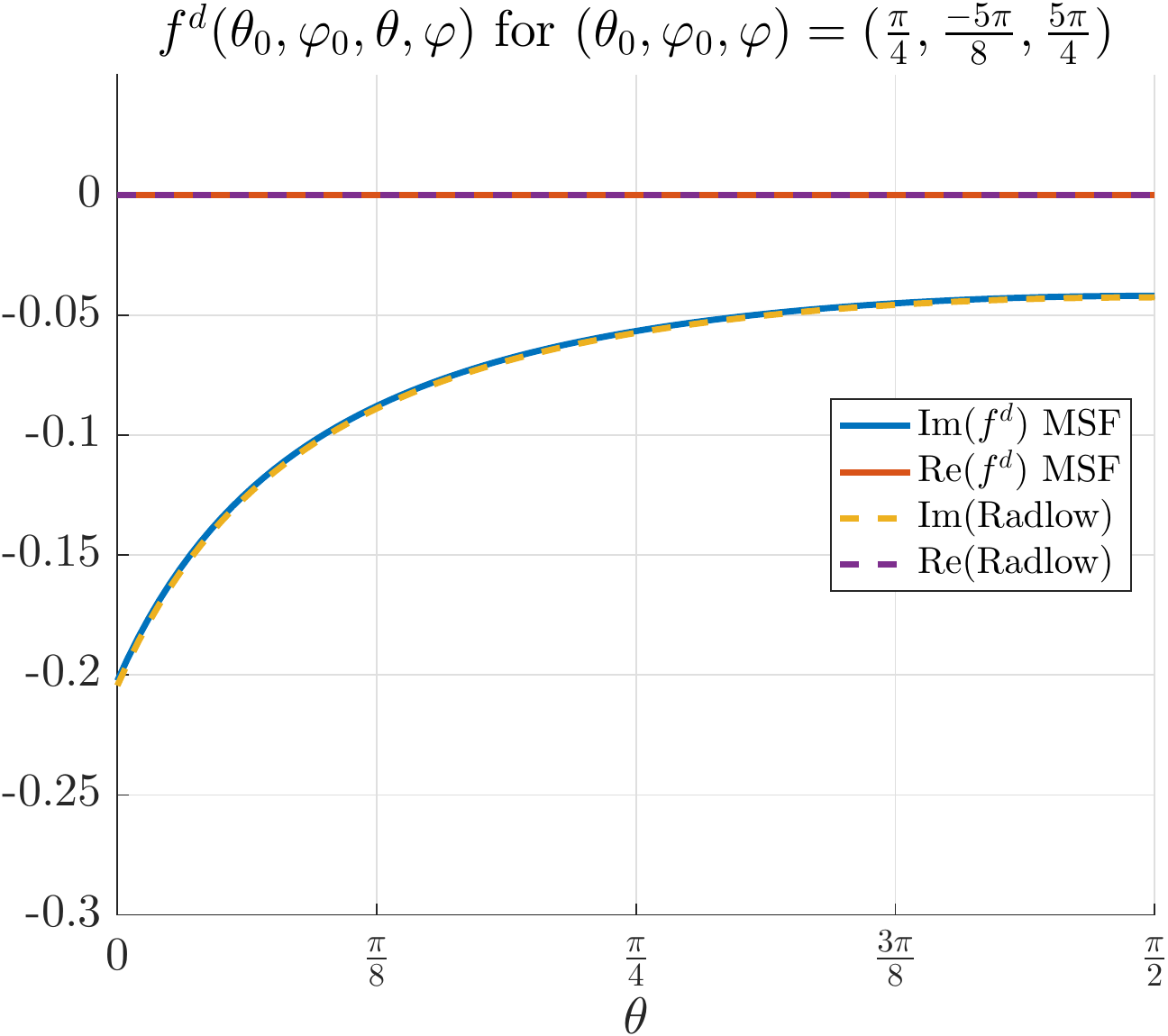}}
 (g)\raisebox{-.1\height}{\includegraphics[width=\mywidth\textwidth]{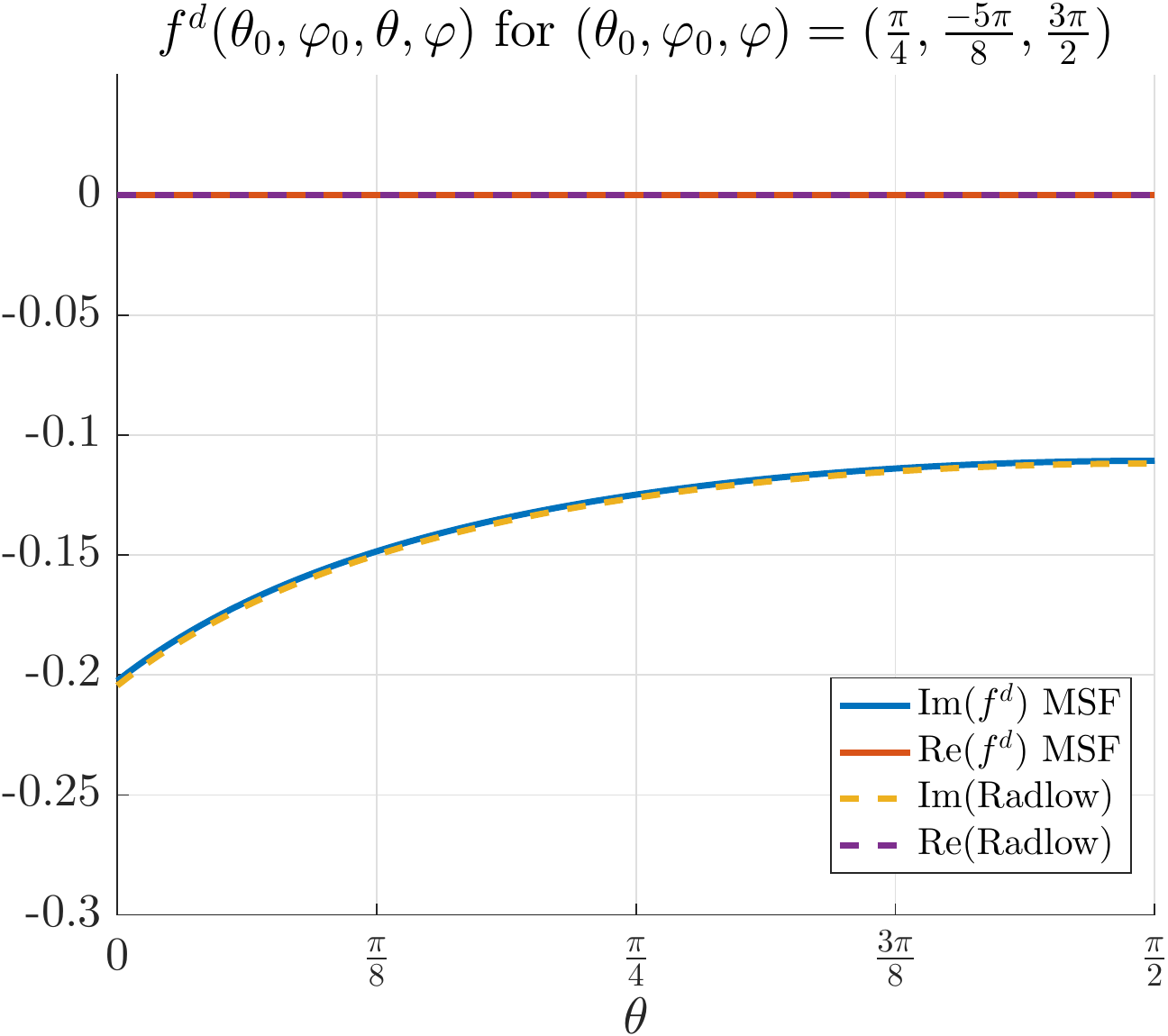}}\quad(h)\raisebox{-.1\height}{\includegraphics[width=\mywidth\textwidth]{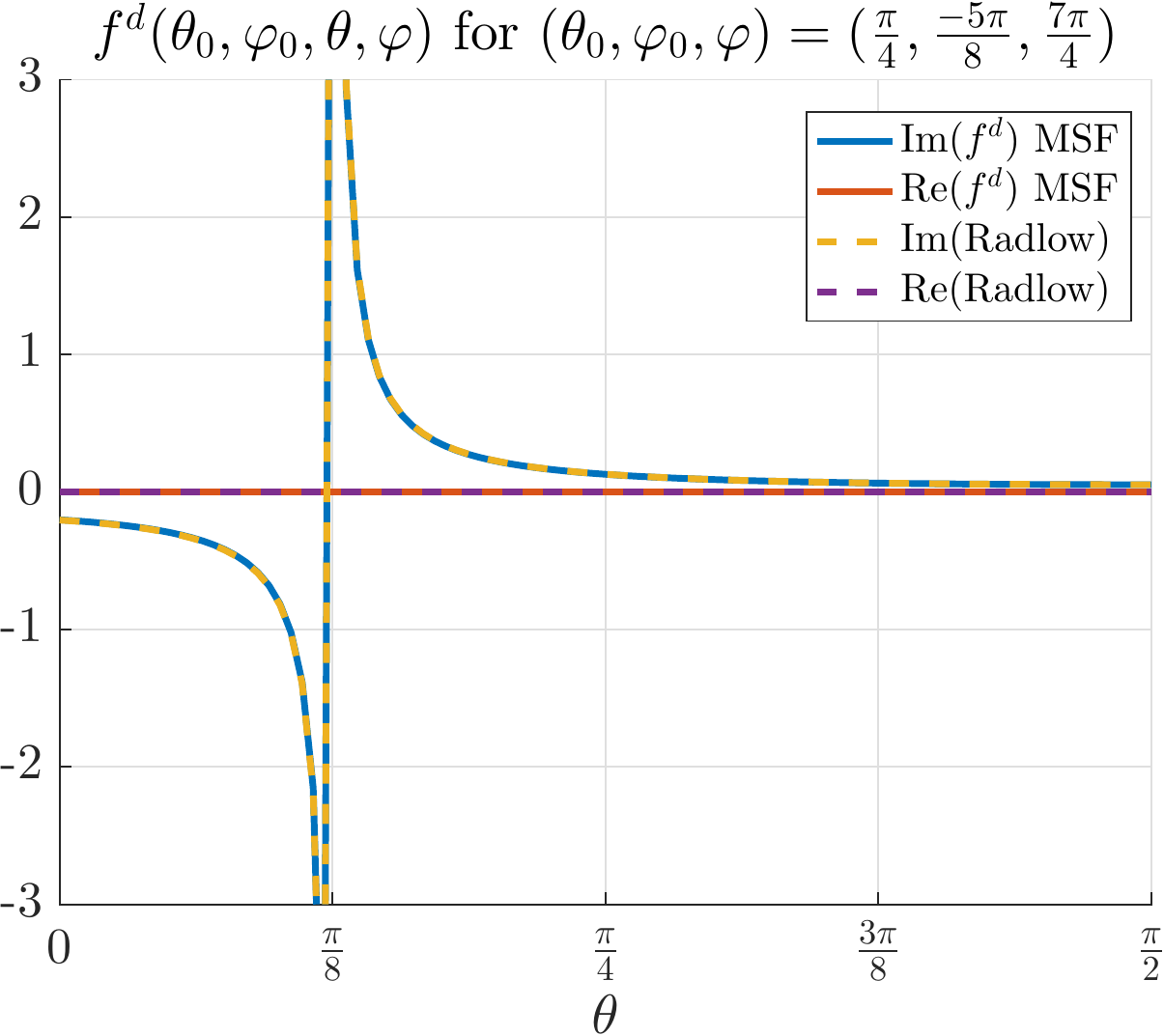}}
  \caption{Diffraction coefficient for incidence $(\theta_0, \varphi_0) =
  \left( \frac{\pi}{4}, \frac{- 5 \pi}{8} \right)$, i.e. we have $a_1 < 0$ and
  $a_2 < 0$, with polar observation angle $\theta \in \left[ 0, \frac{\pi}{2}
  \right]$ and various values of the azimuthal observation angles $\varphi = 0,
  \frac{\pi}{4}, \frac{\pi}{2}, \frac{3 \pi}{4},\pi,
  \frac{5 \pi}{4}, \frac{3 \pi}{2}, \frac{7 \pi}{4}$ (from (a) to (h)).}
\label{fig:run2}
\end{figure}

%

\begin{figure}[htbp]
\centering
  (a)\raisebox{-.1\height}{\includegraphics[width=\mywidth\textwidth]{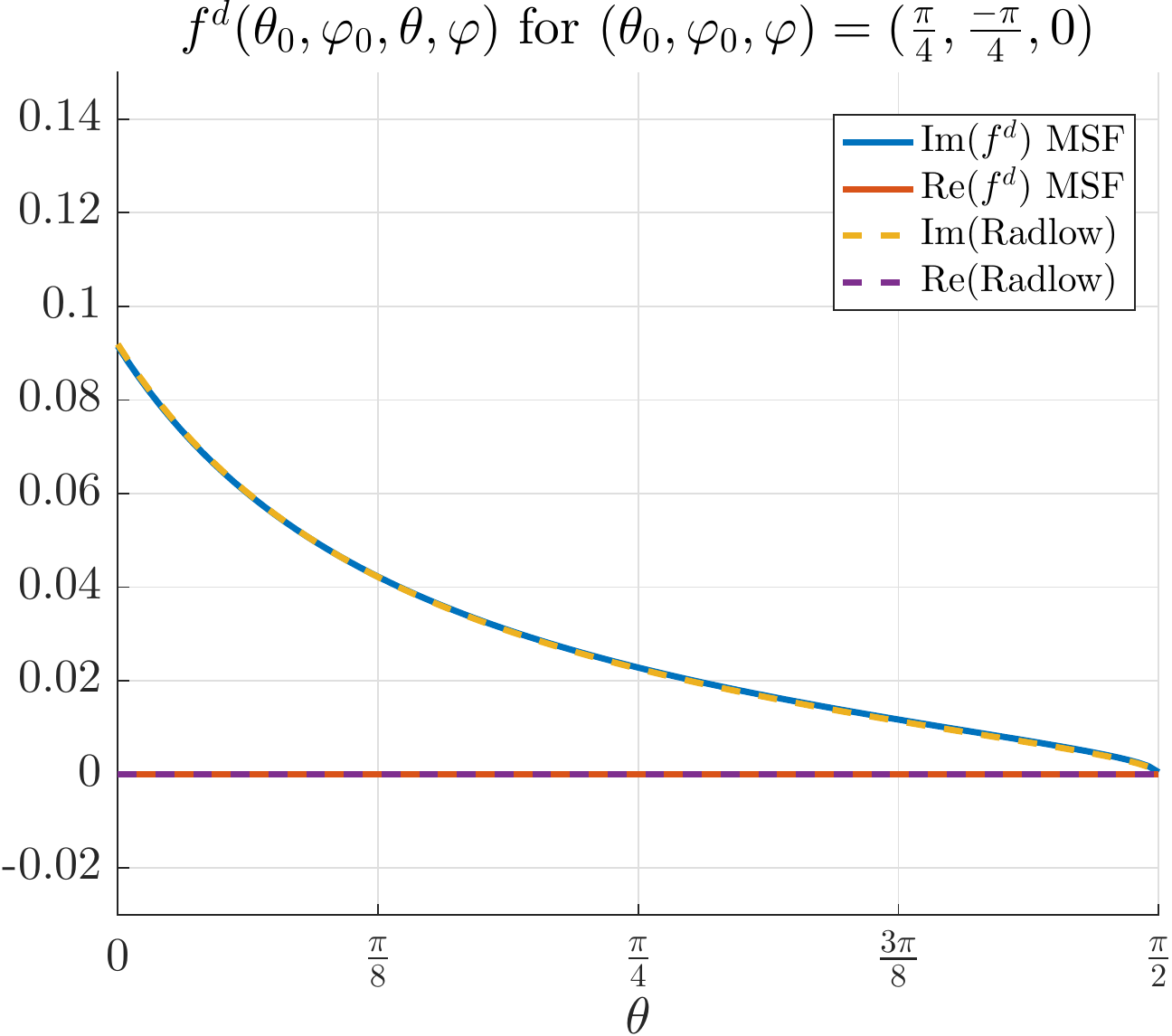}}\quad(b)\raisebox{-.1\height}{\includegraphics[width=\mywidth\textwidth]{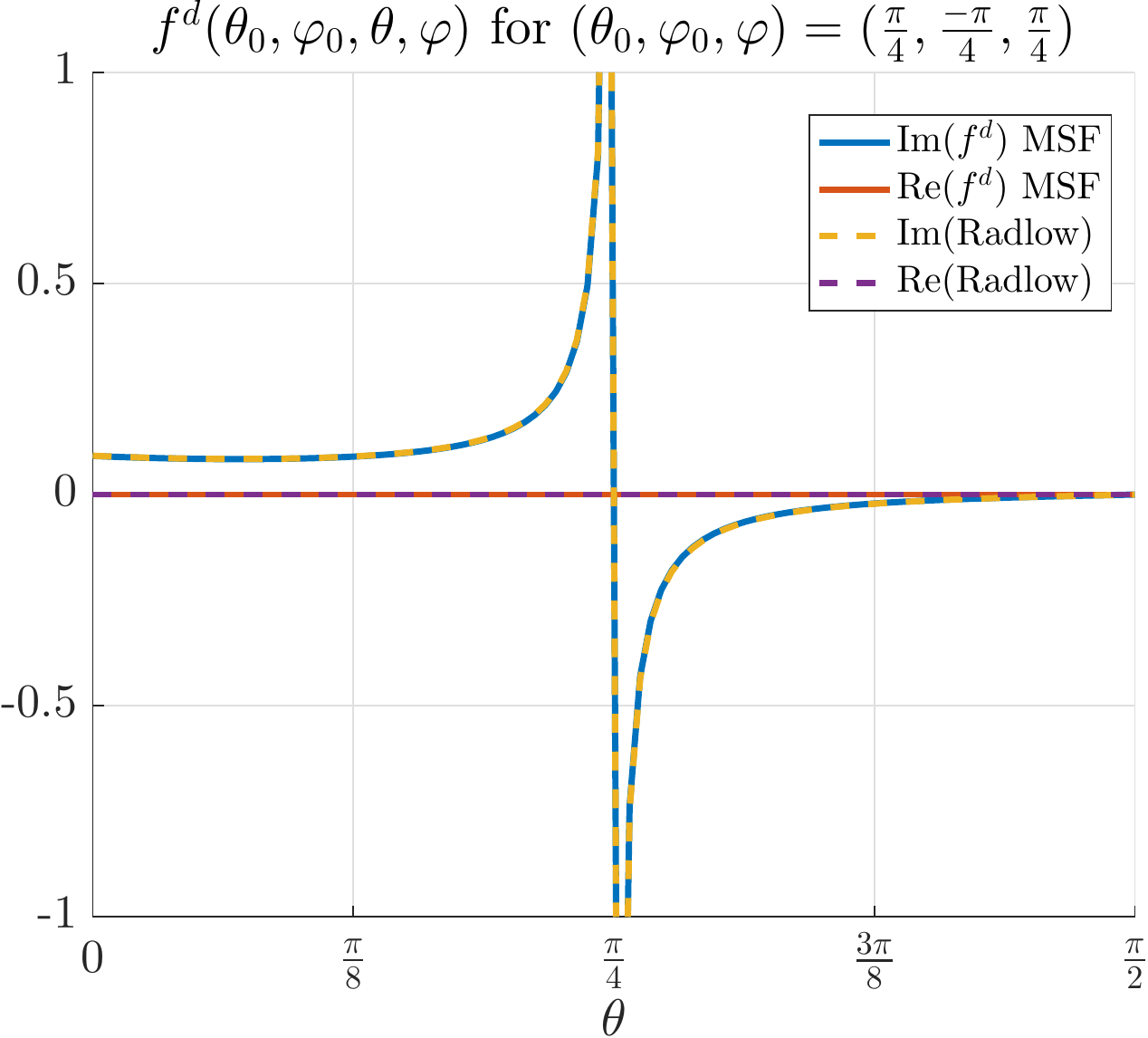}}
  (c)\raisebox{-.1\height}{\includegraphics[width=\mywidth\textwidth]{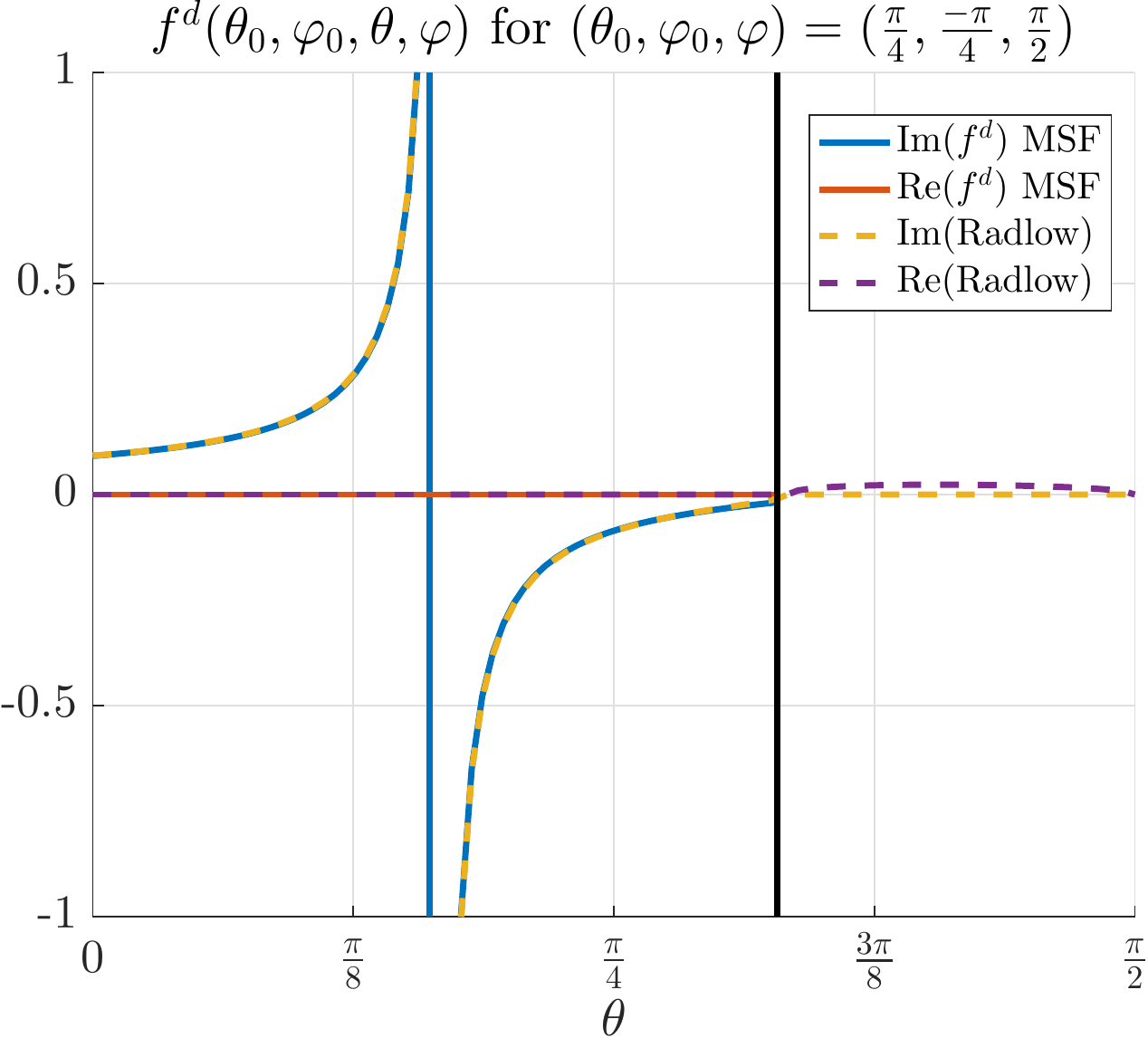}}\quad(d)\raisebox{-.1\height}{\includegraphics[width=\mywidth\textwidth]{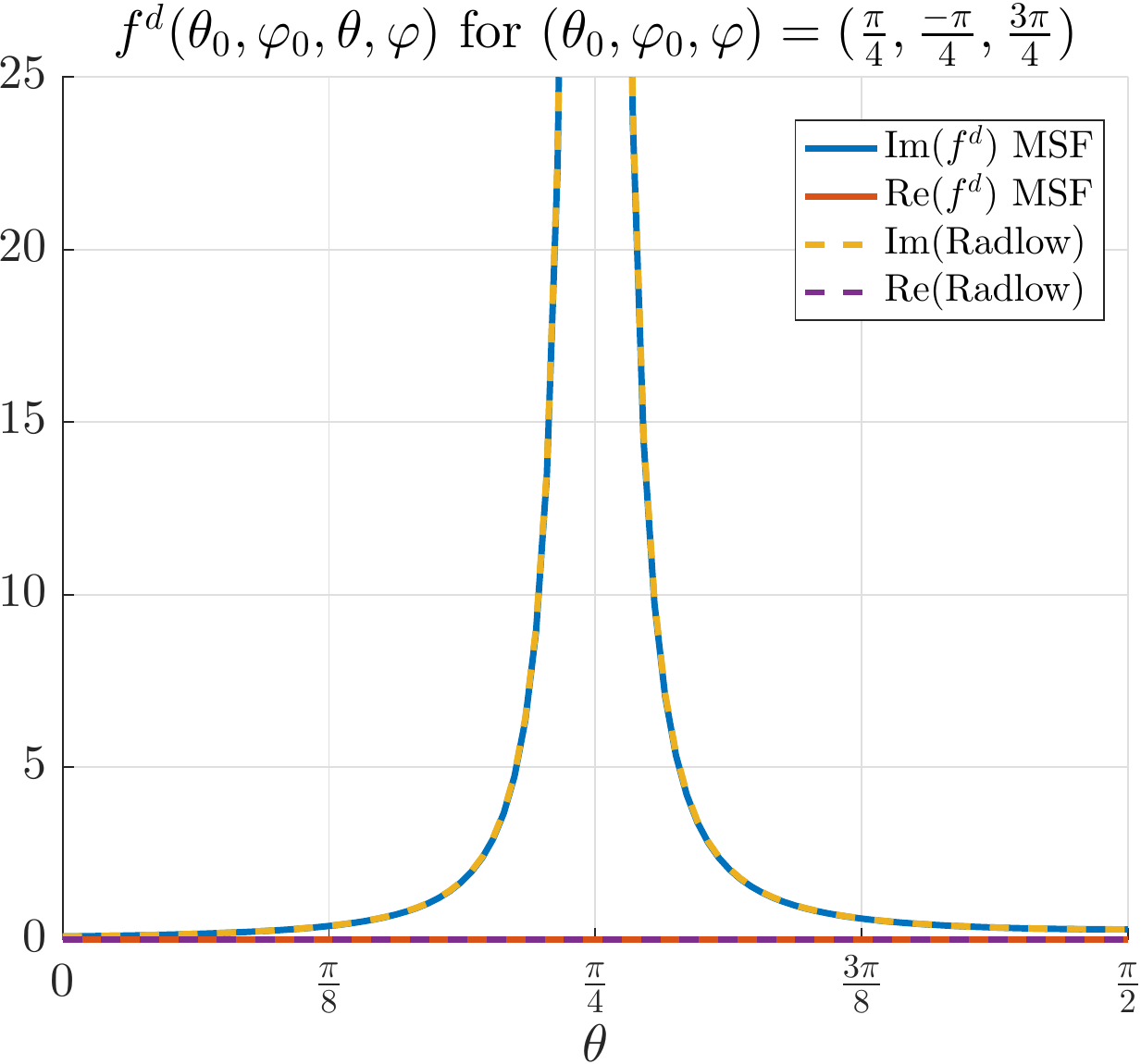}}
  (e)\raisebox{-.1\height}{\includegraphics[width=\mywidth\textwidth]{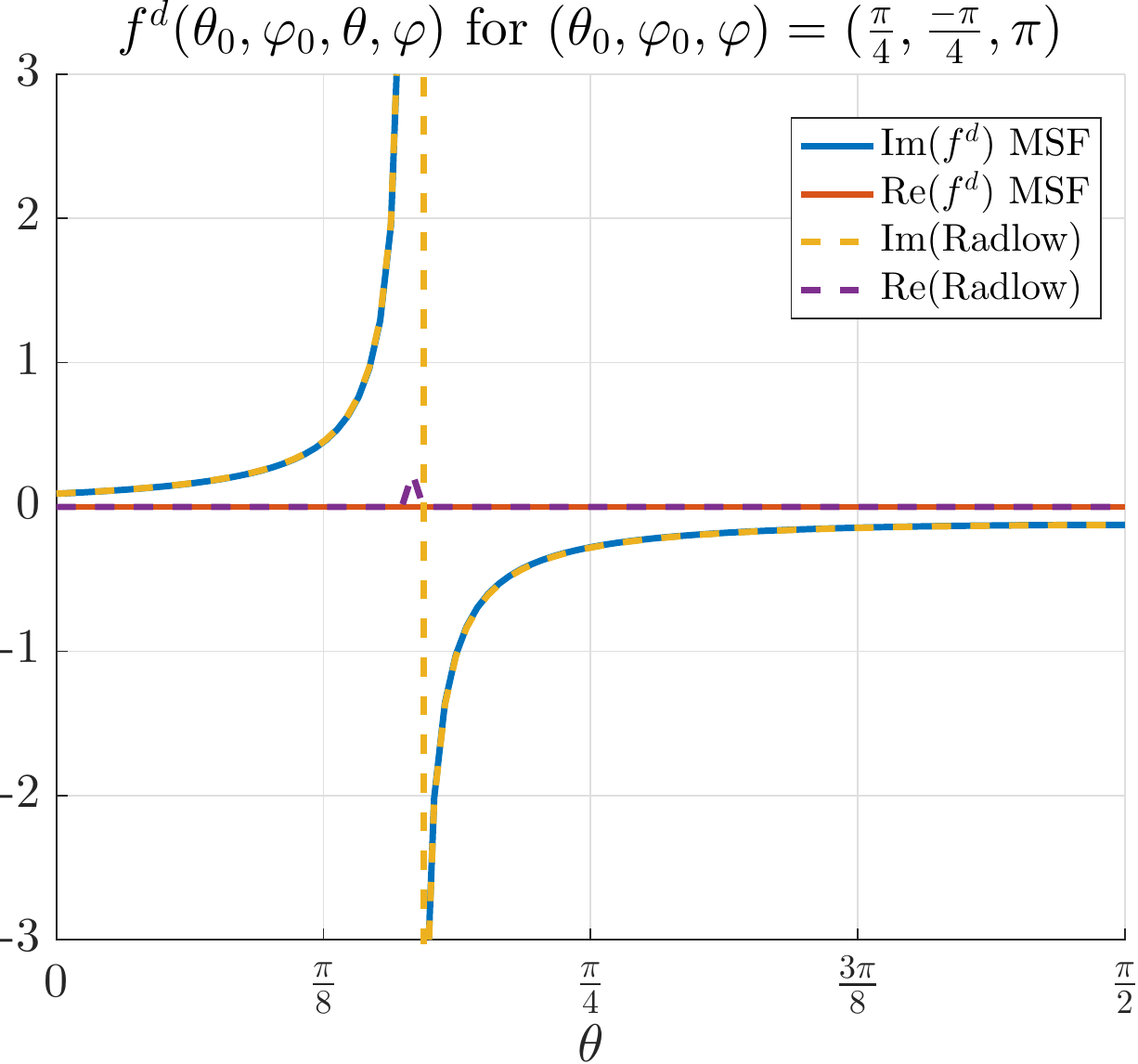}}\quad(f)\raisebox{-.1\height}{\includegraphics[width=\mywidth\textwidth]{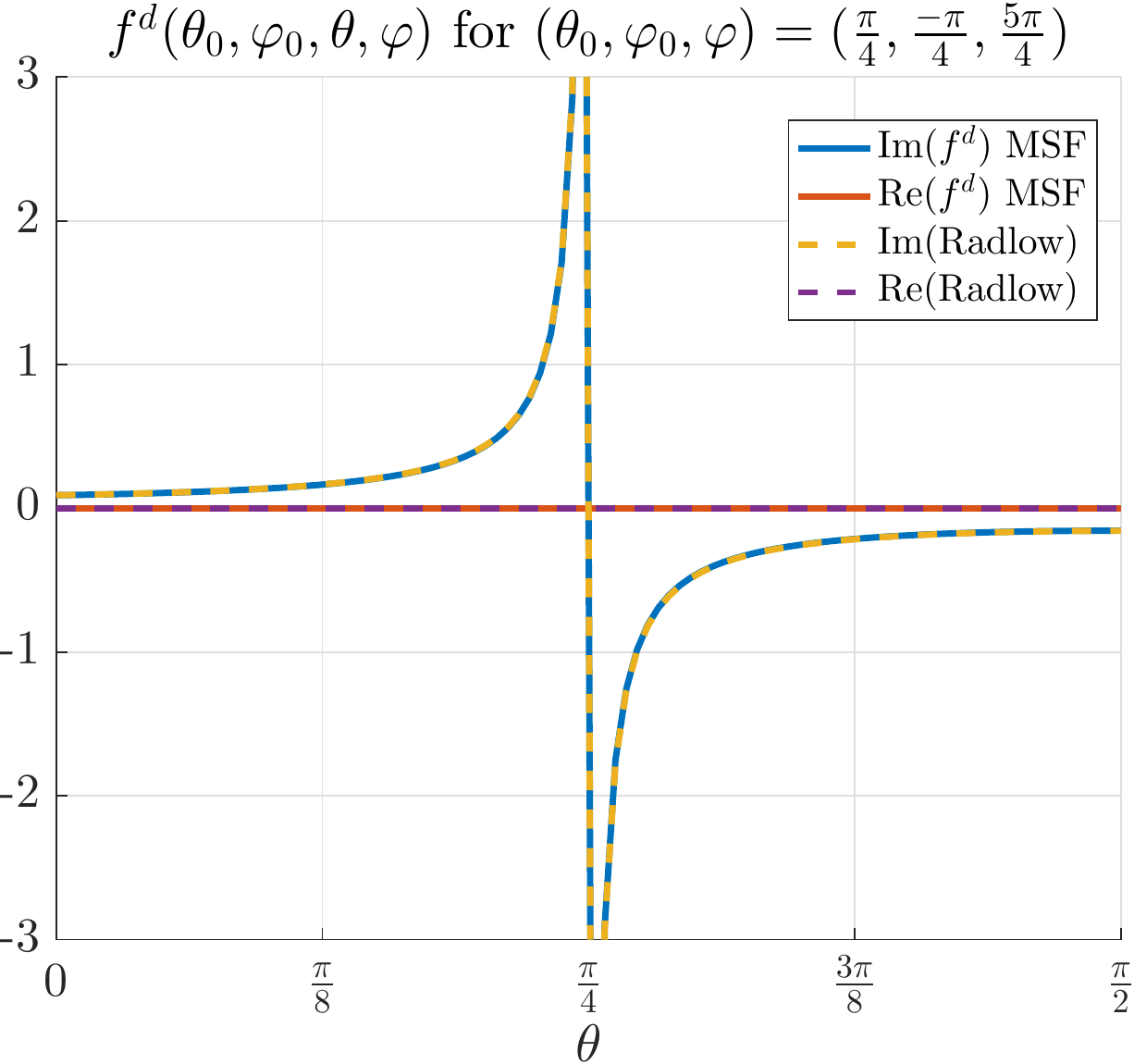}}
  (g)\raisebox{-.1\height}{\includegraphics[width=\mywidth\textwidth]{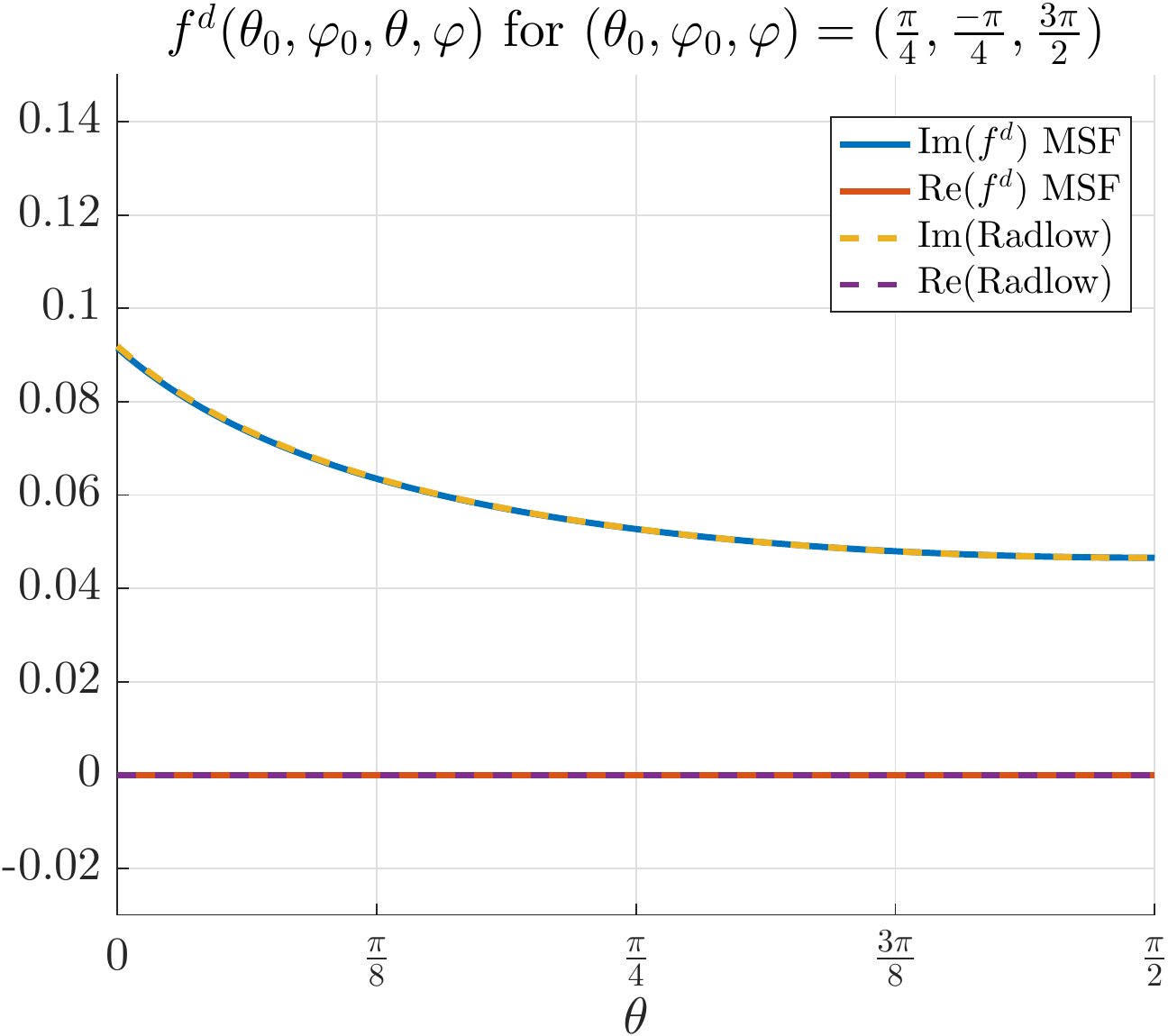}}\quad(h)\raisebox{-.1\height}{\includegraphics[width=\mywidth\textwidth]{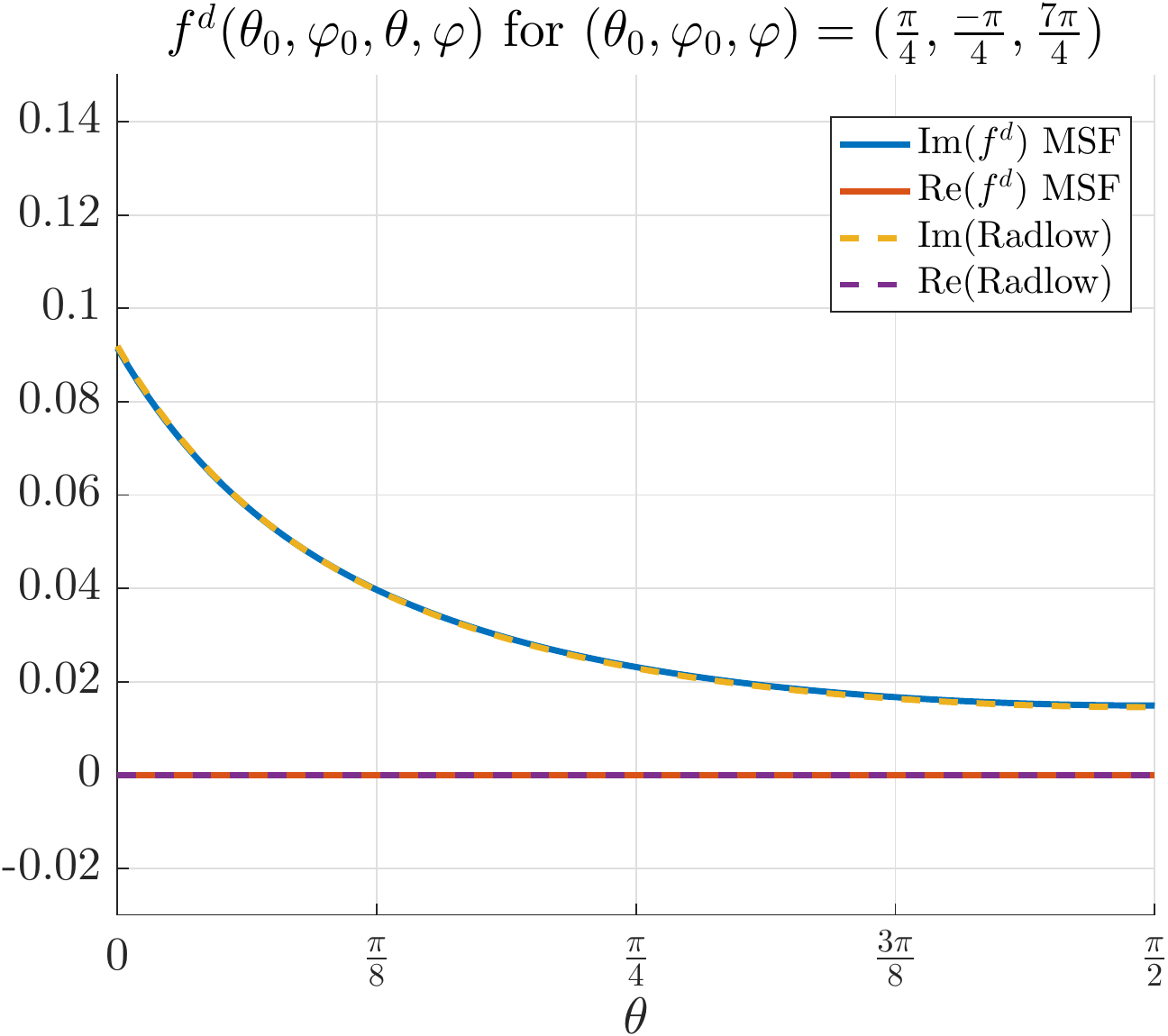}}
  \caption{Diffraction coefficient for incidence $(\theta_0, \varphi_0) =
  \left( \frac{\pi}{4}, \frac{- \pi}{4} \right)$, i.e. we have $a_1 > 0$ and
  $a_2 < 0$, with polar observation angle $\theta \in \left[ 0, \frac{\pi}{2}
  \right]$ and various values of the azimuthal observation angles $\varphi = 0,
  \frac{\pi}{4}, \frac{\pi}{2}, \frac{3 \pi}{4},\pi,
  \frac{5 \pi}{4}, \frac{3 \pi}{2}, \frac{7 \pi}{4}$ (from (a) to (h)). In (c), the black vertical line represents the limit of validity of the
  MSF.}
\label{fig:run3}
\end{figure}

%

\begin{figure}[htbp]
\centering
 (a)\raisebox{-.1\height}{\includegraphics[width=\mywidth\textwidth]{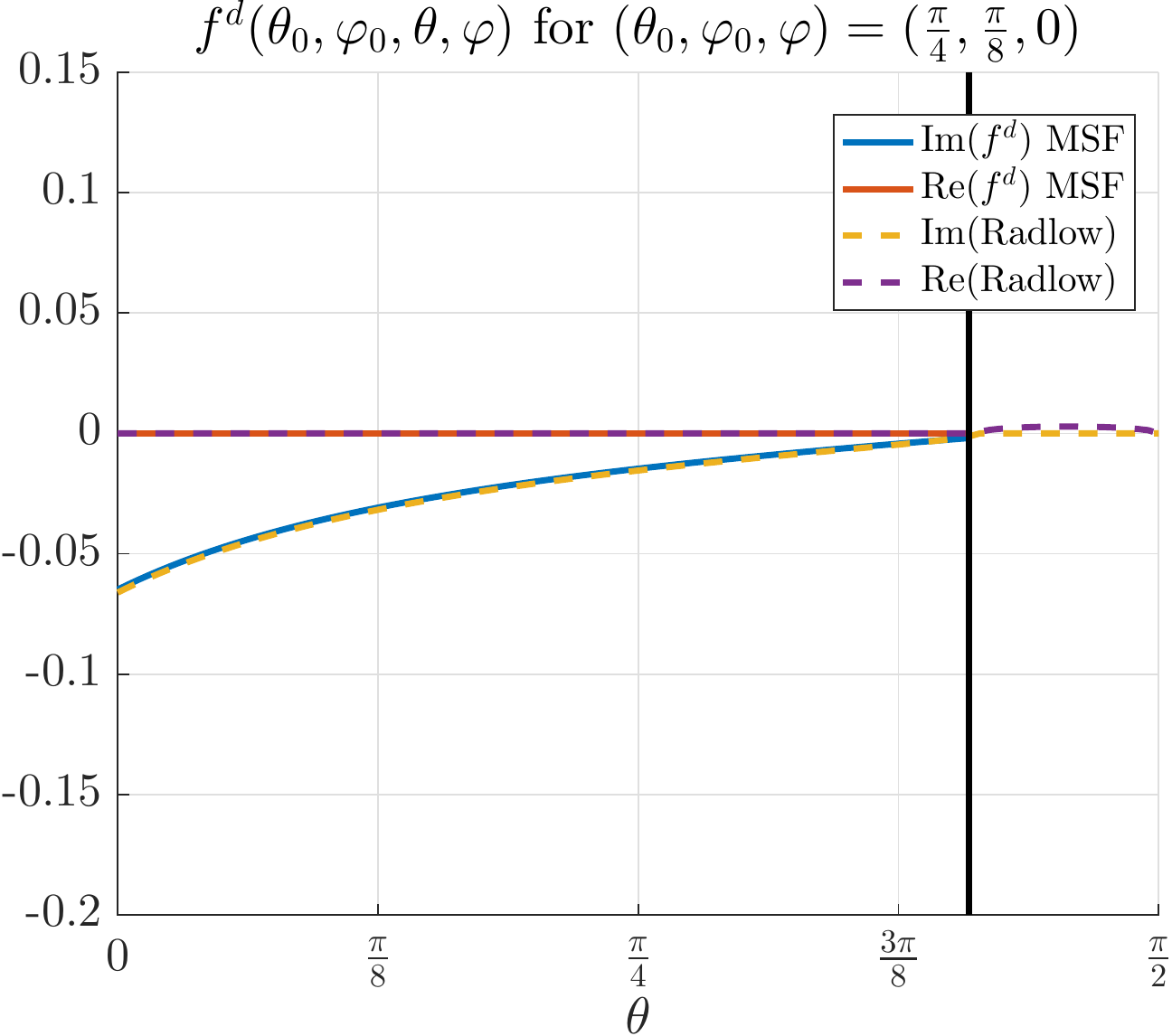}}\quad(b)\raisebox{-.1\height}{\includegraphics[width=\mywidth\textwidth]{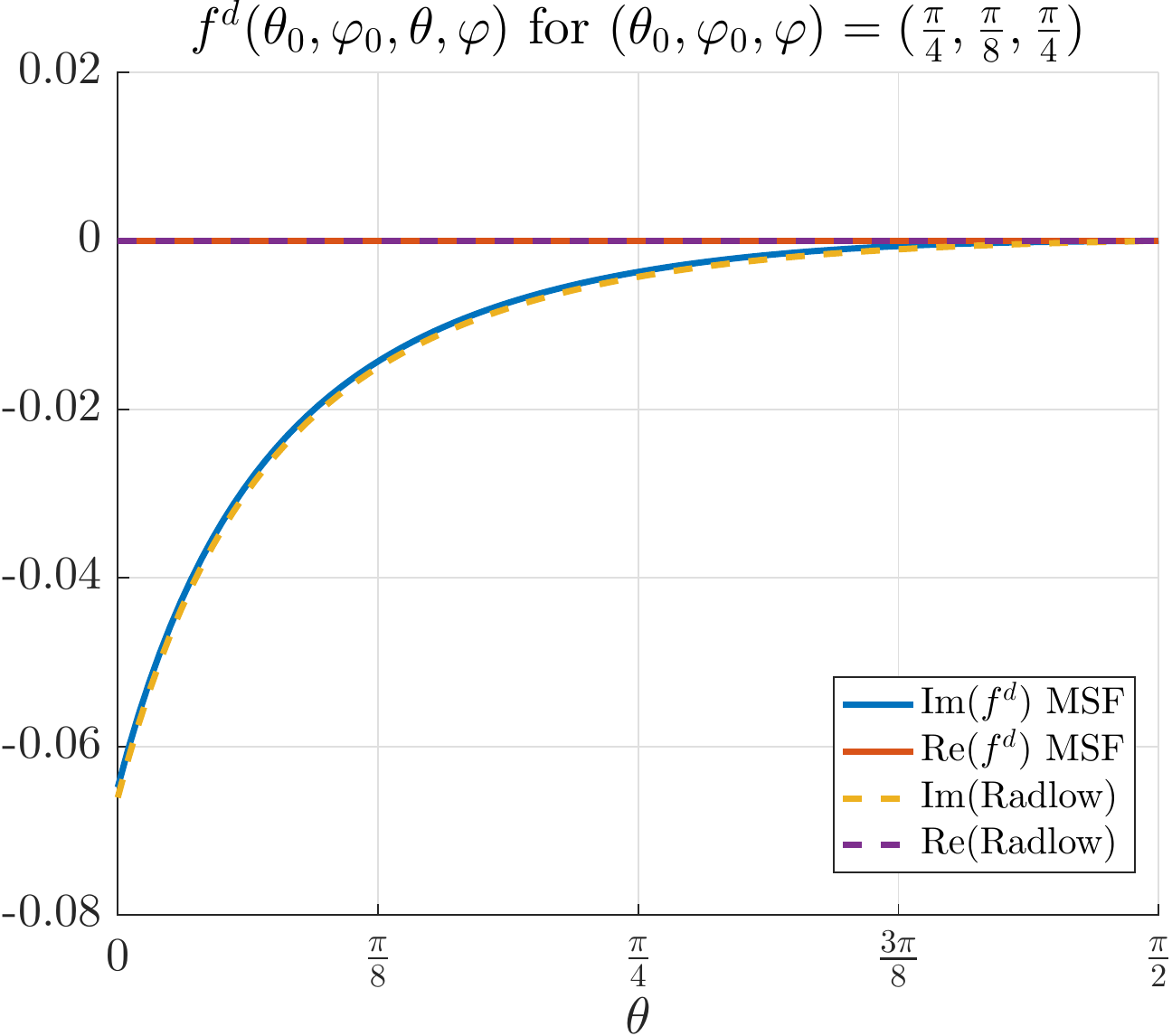}}
 (c)\raisebox{-.1\height}{\includegraphics[width=\mywidth\textwidth]{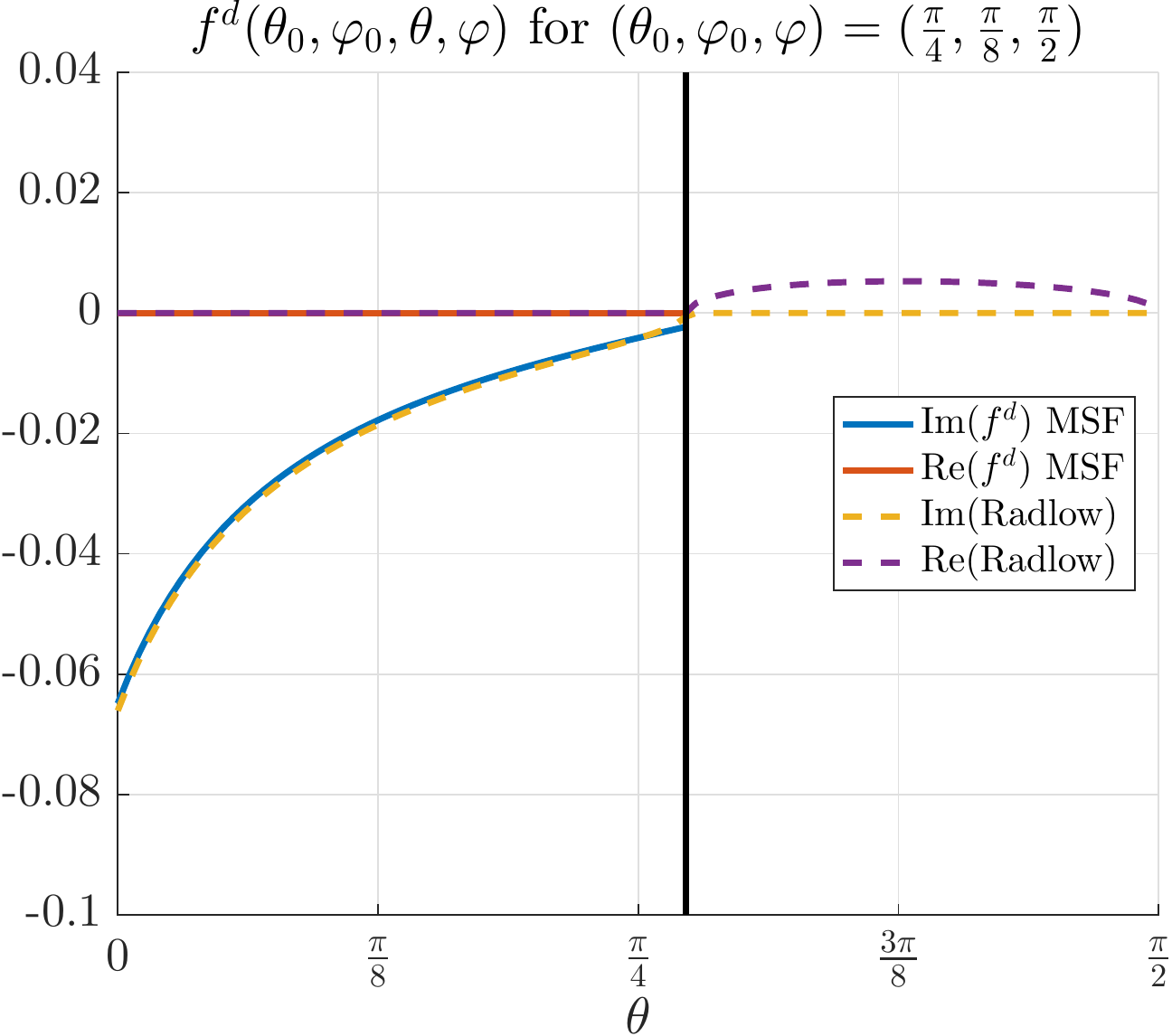}}\quad(d)\raisebox{-.1\height}{\includegraphics[width=\mywidth\textwidth]{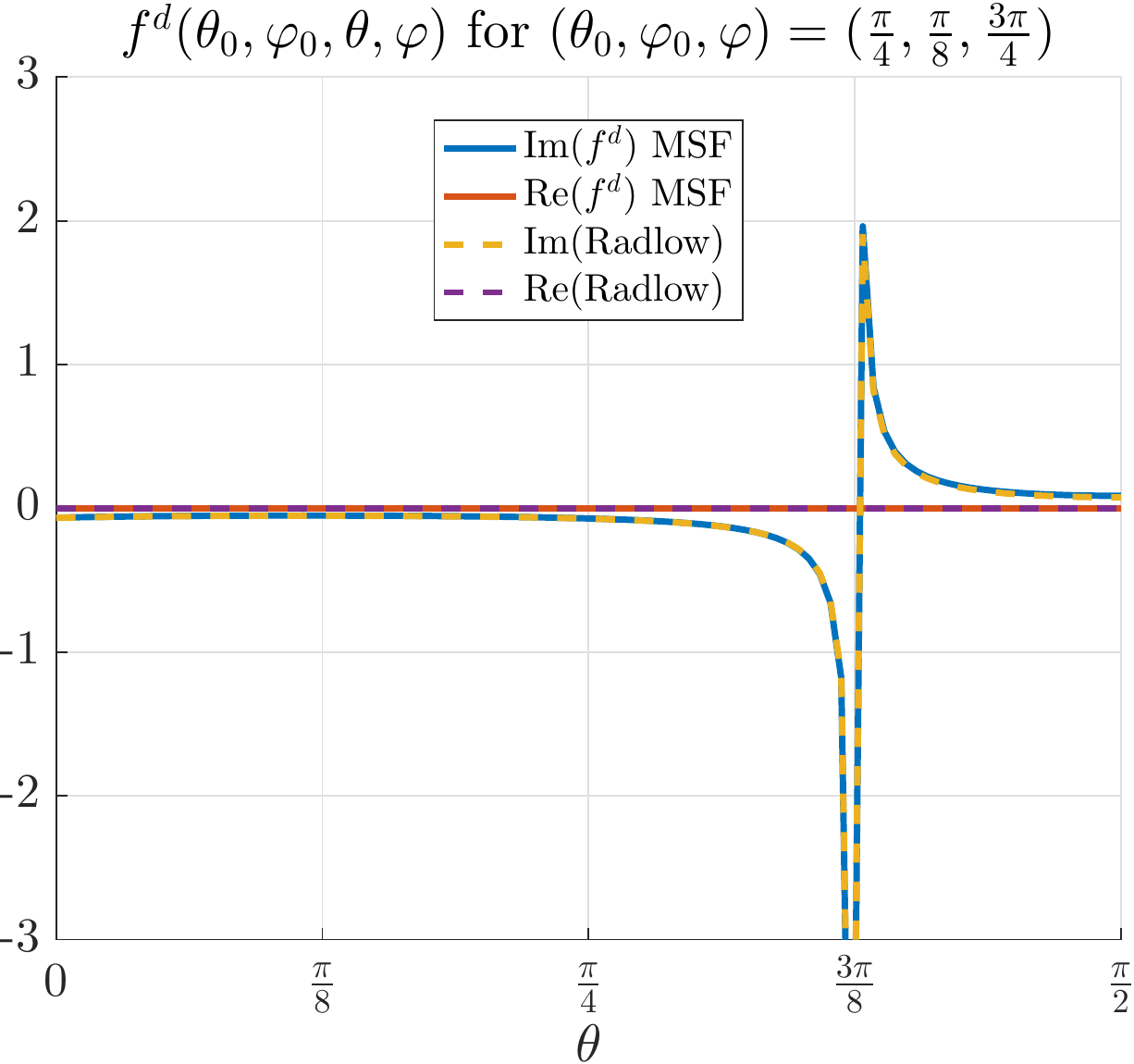}}
 (e)\raisebox{-.1\height}{\includegraphics[width=\mywidth\textwidth]{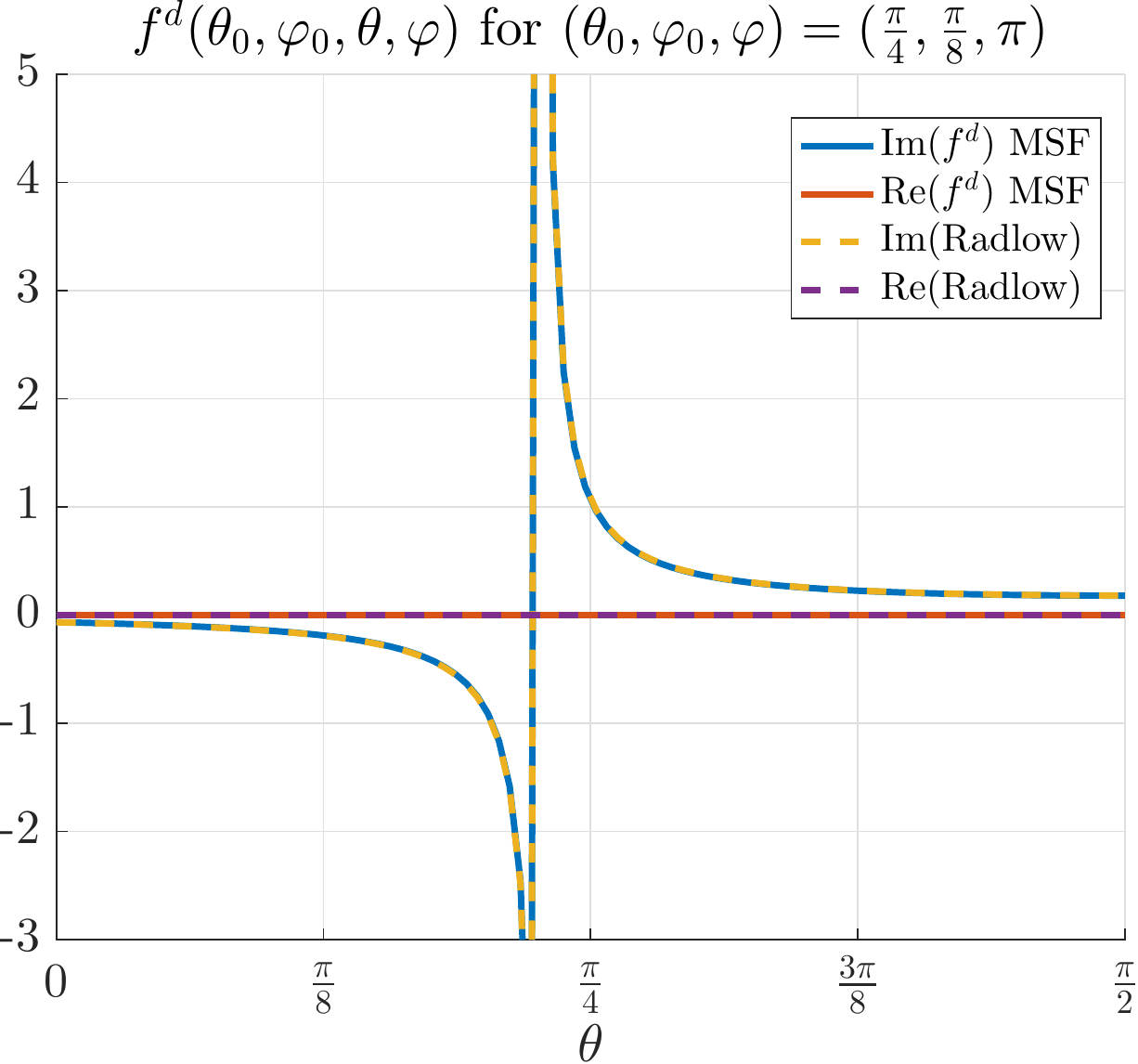}}\quad(f)\raisebox{-.1\height}{\includegraphics[width=\mywidth\textwidth]{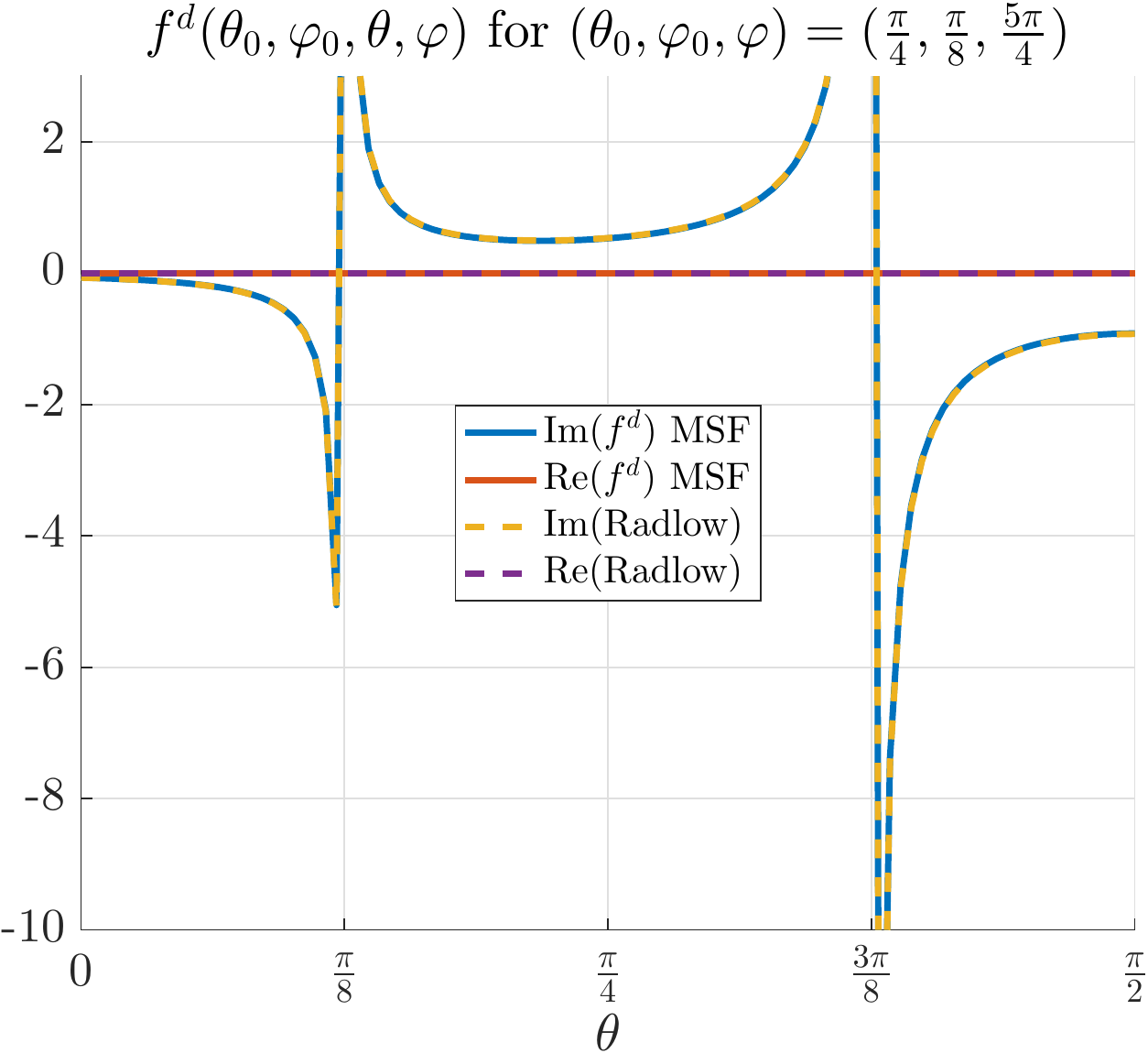}}
 (g)\raisebox{-.1\height}{\includegraphics[width=\mywidth\textwidth]{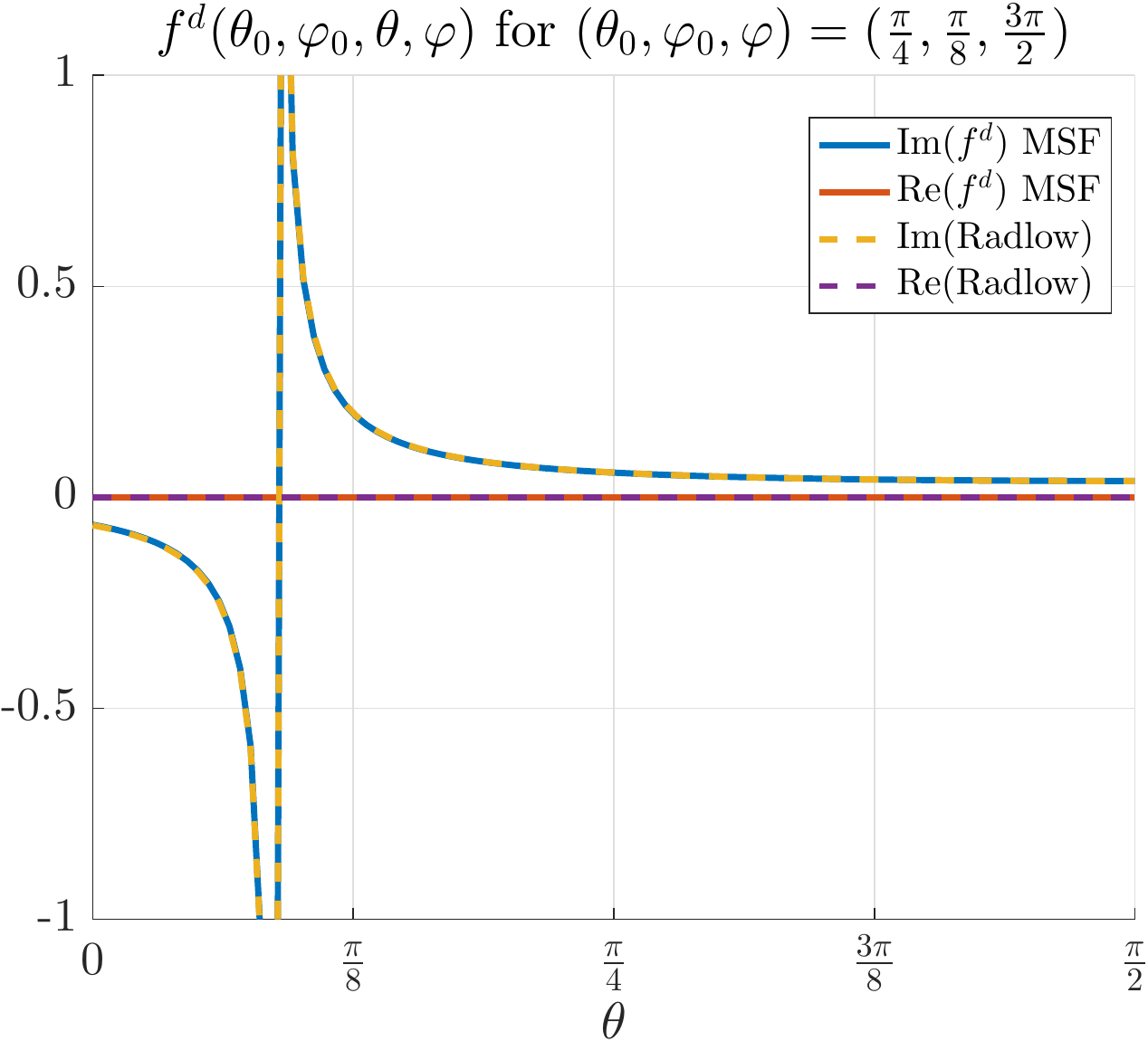}}\quad(h)\raisebox{-.1\height}{\includegraphics[width=\mywidth\textwidth]{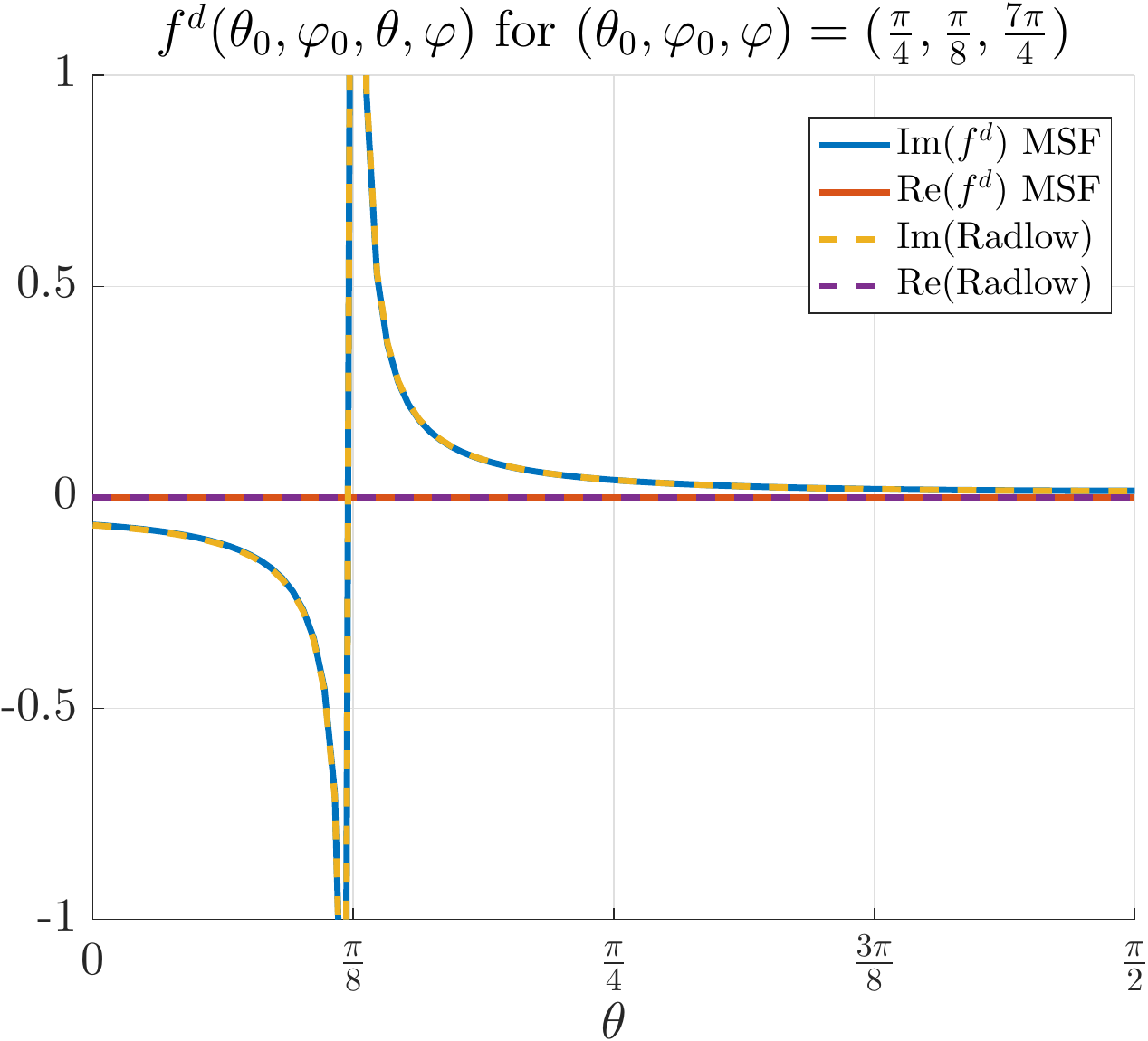}}
  \caption{Diffraction coefficient for incidence $(\theta_0, \varphi_0) =
  \left( \frac{\pi}{4}, \frac{\pi}{8} \right)$, i.e. we have $a_1 > 0$ and
  $a_2 > 0$, with polar observation angle $\theta \in \left[ 0, \frac{\pi}{2}
  \right]$ and various values of the azimuthal observation angles  $\varphi = 0,
  \frac{\pi}{4}, \frac{\pi}{2}, \frac{3 \pi}{4},\pi,
  \frac{5 \pi}{4}, \frac{3 \pi}{2}, \frac{7 \pi}{4}$ (from (a) to (h)). In (a) and (c), the vertical black line represents the limit of validity of the
  MSF.}
\label{fig:run4}
\end{figure}

  


The fact that Radlow's ansatz produces extremely accurate results for the
diffraction coefficient is indeed surprising, but such possibility was not
ruled out in Albani's work {\cite{albani}}. Indeed Albani's approach to
showing that Radlow's ansatz (let us call it $F^{\tmop{Ra}}_{+ +}
(\tmmathbf{\alpha})$) was incorrect was to demonstrate that the resulting physical
field, $u^{\tmop{Ra}} (x_1, x_2, x_3)$ did not satisfy the boundary condition,
i.e. was not equal to $- e^{- i (a_1 x_1 + a_2 x_2)}$ on the quarter-plane
$x_{1, 2} > 0$. An interesting point, however, was that he showed that as both
$x_1$ and $x_2$ tend to infinity simultaneously, we have
\begin{eqnarray}
u^{\tmop{Ra}} (x_1 > 0, x_2 > 0, 0) - (- e^{- i (a_1 x_1 + a_2 x_2)}) & = &
\mathcal{O} ((x_1^2 + x_2^2)^{- 3 / 2}),
\end{eqnarray}
implying that in a way, the boundary conditions are {\tmem{asymptotically
		satisfied}} away from the vertex and the edges. The rapidity of the decay (one
over the cube of the distance to the vertex) being much higher than the decay
of the spherical wave (one over this distance) may be the beginning of an
explanation as to why Radlow's ansatz performs so well in that case.
\RED{It has to be said however that the agreement between the two methods cannot be perfect. Indeed, if it were to be, then $F_{++}$ and $F_{++}^{Ra}$ would have to be exactly the same on a non-isolated region, and hence, due to the theory of analytic functions, they will have to be the same everywhere, which as we showed would violate the compatibility condition. There must hence exist a numerical discrepancy between the two methods. In order to find it, we made sure that the MSF and the Radlow's ansatz where accurately evaluated up to a relative error of the order $\mathcal{O}(10^{-5})$ and looked at the pointwise difference between the two methods for the particular testcase of Figure~\ref{fig:run0}(g). The results are displayed in Figure \ref{fig:pointwise-difference}, and one can see that the relative error is of the order $\mathcal{O}(10^{-3})$, two orders of magnitude higher than the precision with which both methods were computed. We can hence conclude that this is an actual discrepancy between the two methods, and not a numerical artefact.}

\begin{figure}[htbp]
	\centering
	\includegraphics[width=0.5\textwidth]{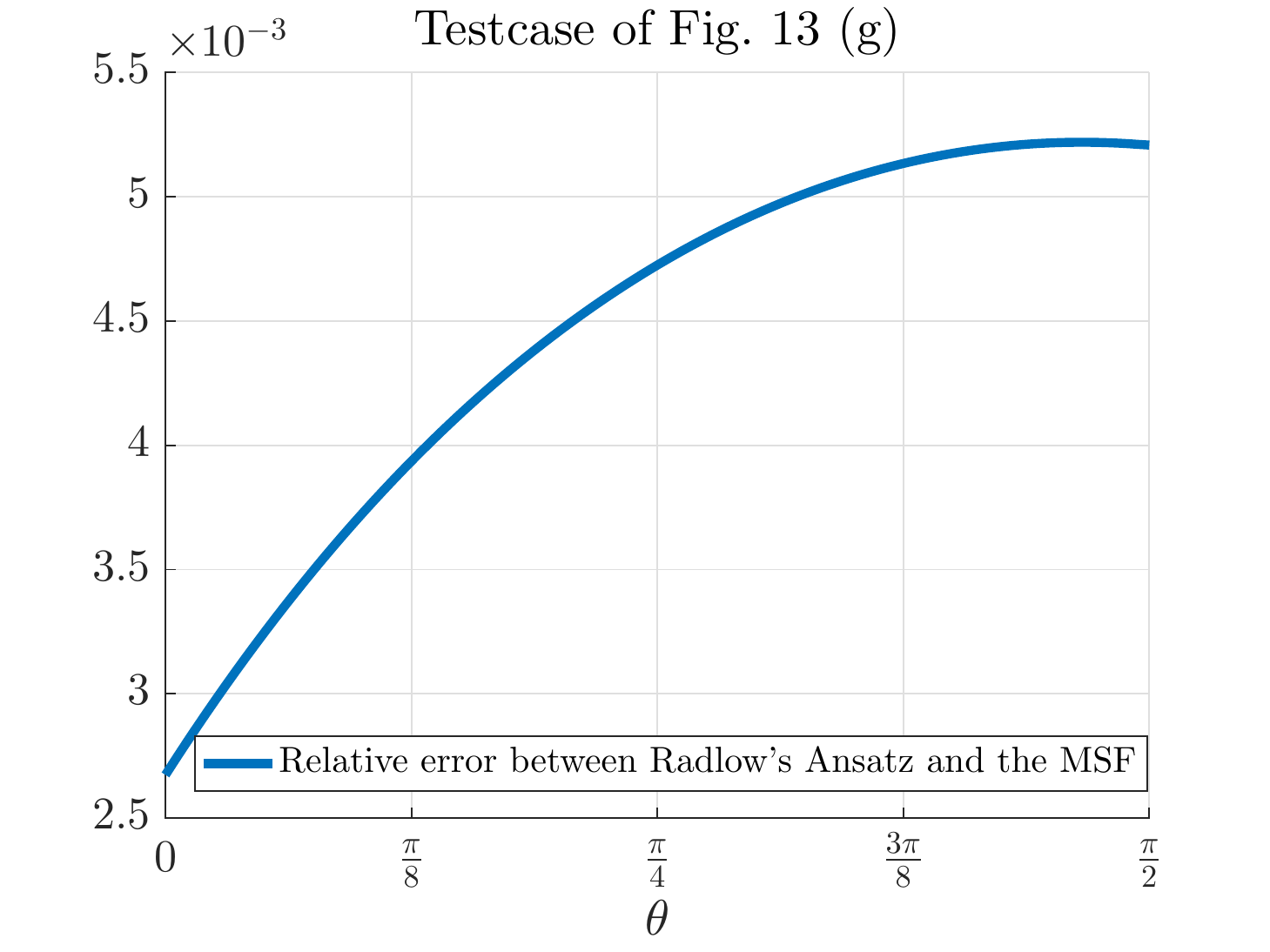}
	\caption{\RED{Pointwise relative error between the diffraction coefficient obtained by the MSF and by Radlow's ansatz for the testcase of Figure~\ref{fig:run0}~(g).}}
	\label{fig:pointwise-difference}
\end{figure}


\section{Conclusion}

In this paper, we revisited Radlow's double Wiener-Hopf approach for the
Dirichlet quarter-plane problem. We have tried to add more clarity and precision to his innovative approach, with an aim to obtain a \textit{constructive} method of solution of this canonical boundary value problem. The inverse Fourier transform
(\ref{eq:explicituFK}), gives the solution in terms of an unknown function
$F_{+ +}$, that depends on two complex variables. We
reduced the problem to two equations, one, (\ref{eq:genericWHeq1}), expresses
$F_{+ +}$ as the sum of two terms, one containing the unknown function $G_{+-}$ and the other being Radlow's ansatz. This, on the one hand, gives a constructive way of obtaining the
ansatz, and on the other hand, offers yet another reason why this ansatz
cannot be the true solution. The second equation, (\ref{eq:genericWHeq2}),
called the {\tmem{compatibility equation}}, involves solely
the unknown function $G_{+ -}$ and could be key to determining this
crucial unknown function.

Finally, following a steepest-descent analysis, we have related $F_{+ +}$ to
the diffraction coefficient $f^d$. Numerical results show that when choosing
$F_{+ +}$ as per Radlow's ansatz, we obtain surprisingly accurate results for the diffraction coefficient. In fact, the
results seem to agree \RED{very well} 
with those obtained by the established Modified Smyshlyaev Formulae, where this method is valid. \RED{Theoretically, it is however impossible for this agreement to be perfect, and we have shown that there exists a small discrepancy between the two methods, with a relative error of order $\mathcal{O}(10^{-3})$.} It should be
noted that the MSF is a very quick way of evaluating the diffraction
coefficient; however, Radlow's ansatz, and the factorisation formulae
provided herein, is even faster (computing the Radlow result for each graph of
Section \ref{sec:results} takes about 1s on a standard laptop). This observation
naturally opens some interesting questions:
\begin{itemize}[leftmargin=.15in]
	\item is the diffraction coefficient arising from  Radlow's ansatz \RED{a very good far-field approximation}, even in the region
	inaccessible by the MSF;
	
	\item why does the near-field have seemingly no influence on the far-field
	behaviour;
	
	\item can we find a constructive method for determining the function $G_{+ -}$, and hence a unique formulation
	reconciling near-field and far-field;
	
	\item {can we take a similar approach in the Neumann case?}
\end{itemize} 
We hope to be able to answer these points in our future work, several of which could have profound consequences on how we approach diffraction problems in general. 

\appendix

\section{Factorisation of $K_{- \circ}$ and $K_{+
		\circ}$}\label{app:factoK-o}

Let us show how the factorisation of $K_{- \circ}$ is obtained. The
factorisation of $K_{+ \circ}$ is obtained in a very similar way. Introduce
the auxiliary function $\mathfrak{K}_{- \circ}$ as
\begin{align*}
\mathfrak{K}_{- \circ} (\tmmathbf{\alpha})  &=  \kappa (k, \alpha_2) K_{-
	\circ}^2 (\tmmathbf{\alpha}) = \tfrac{\kappa (k, \alpha_2)}{\kappa (k,
	\alpha_2) - \alpha_1} = \tfrac{1}{1 - \frac{\alpha_1}{\kappa (k, \alpha_2)}} \, \cdot
\end{align*}
Naturally, for a given $\alpha_2$ in $\mathcal{A}_2$, $\mathfrak{K}_{- \circ}
(\tmmathbf{\alpha})$ remains a {\tmem{minus function}} when seen as a function
of $\alpha_1$. Plots of the auxiliary function $\mathfrak{K}_{- \circ}
(\tmmathbf{\alpha})$ are provided in Figure \ref{fig:weirdkminuso}.

\begin{figure}[htbp]
	\centering
	\includegraphics[width=0.4\textwidth]{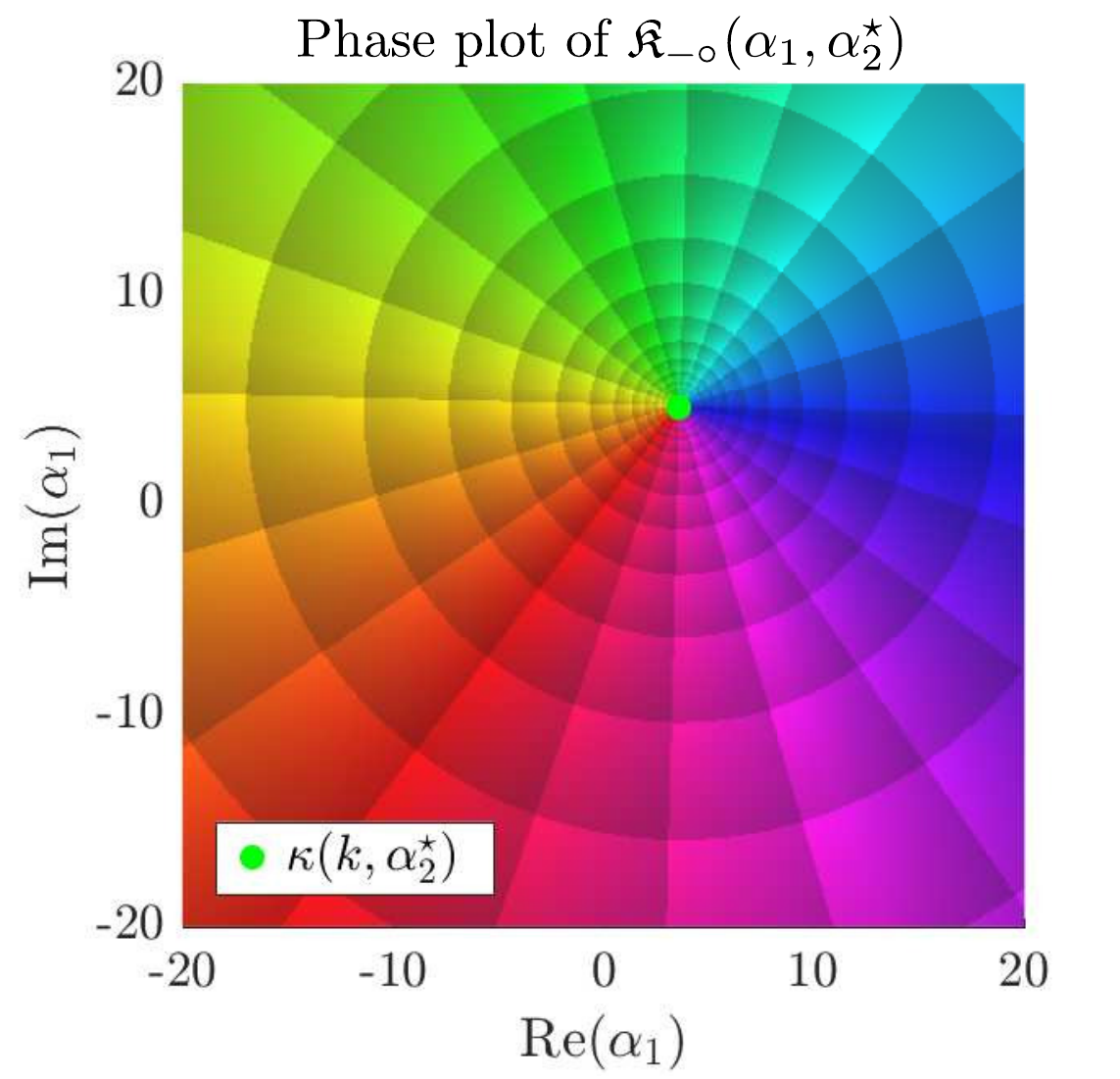}\includegraphics[width=0.4\textwidth]{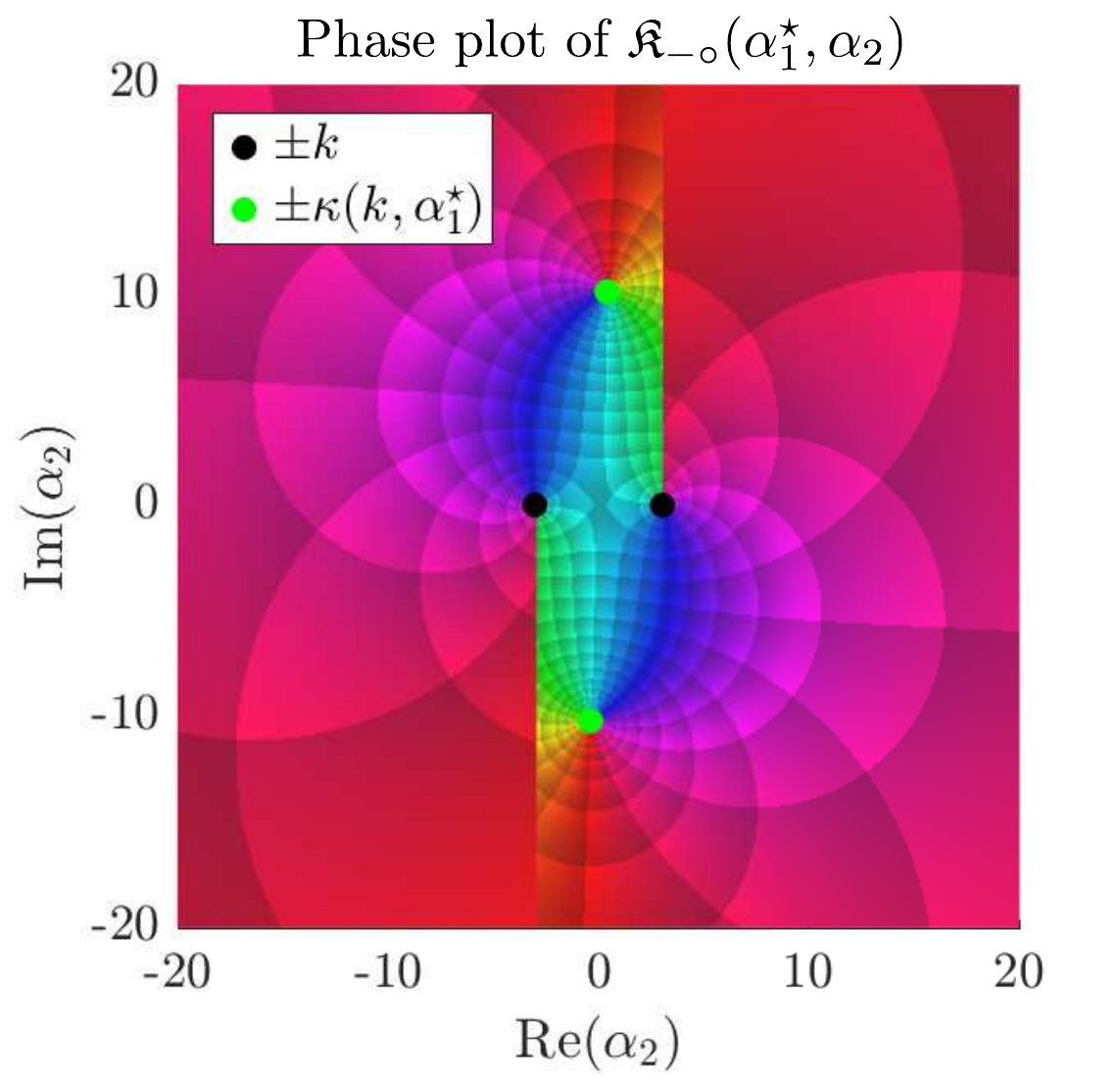}
	\caption{Left: Phase plot of the function $\mathfrak{K}_{- \circ} (\alpha_1,
		\alpha_2^{\star})$ for $\alpha_2^{\star} =\mathcal{A}_2 (5)$ in the
		$\alpha_1$ complex plane. Right: Phase plot of the function $\mathfrak{K}_{-
			\circ} (\alpha_1^{\star}, \alpha_2)$ for $\alpha_1^{\star} =\mathcal{A}_1
		(10)$ in the $\alpha_2$ complex plane. }
	\label{fig:weirdkminuso}
\end{figure}

Note that for $\mathfrak{K}_{- \circ} (\alpha_1, \alpha_2^{\star})$ (Figure
\ref{fig:weirdkminuso}, left) the point $\alpha_1 = \kappa (k,
\alpha_2^{\star})$ is not a branch point anymore, but just a simple pole. For
$\mathfrak{K}_{- \circ} (\alpha_1^{\star}, \alpha_2)$ (Figure
\ref{fig:weirdkminuso}, right), as expected, $\alpha_2 = \pm k$ are branch
points, while $\alpha_2 = \pm \kappa (k, \alpha_1^{\star})$ now correspond to
two simple poles.

Let us now set $\alpha_1 \in \tmop{LHP}_1$. Now for a given $\alpha_2$ in
$\mathcal{A}_2$ (where $\mathfrak{K}_{- \circ} (\tmmathbf{\alpha})$ is
analytic when considered as a function of $\alpha_2$), we can make use of
Corollary \ref{cor:CauchyFac} to write $\mathfrak{K}_{- \circ} (\tmmathbf{\alpha}) = \mathfrak{K}_{- -}
(\tmmathbf{\alpha}) \mathfrak{K}_{- +} (\tmmathbf{\alpha})$,
the equality being valid on $\mathcal{D}_{- \circ}$, where $\mathfrak{K}_{-
	-} (\tmmathbf{\alpha})$ is analytic in $\tmop{LHP}_2$ and $\mathfrak{K}_{- +}
(\tmmathbf{\alpha})$ is analytic in $\tmop{UHP}_2$ when both are considered as
functions of $\alpha_2$. And, these are given by
\RED{\begin{eqnarray*}
		\mathfrak{K}_{- +} (\alpha_1, \alpha_2)  =  e^{ \frac{1}{2 i \pi}
			\int_{\mathcal{A}_{\varepsilon}^b} \frac{\overset{\swarrow}{\log}
				(\mathfrak{K}_{- \circ} (\alpha_1, z))}{z - \alpha_2} \, \mathd z } &\text{ and }& \mathfrak{K}_{- -} (\alpha_1, \alpha_2)  = e^{ \frac{- 1}{2 i
				\pi} \int_{\mathcal{A}_{\varepsilon}^a}
			\frac{\overset{\swarrow}{\log} (\mathfrak{K}_{- \circ} (\alpha_1, z))}{z -
				\alpha_2} \, \mathd z },
\end{eqnarray*}}
where $\overset{\swarrow}{\log}$ was defined in Section
\ref{sec:usefulfunctions}. The choice of this particular logarithm is in fact
extremely important in order to avoid crossings between branch cuts and the
contour of integration. Using the exact expression of $\mathfrak{K}_{- \circ}
(\tmmathbf{\alpha})$, this can be simplified to
\RED{\begin{eqnarray*}
		\mathfrak{K}_{- +} (\alpha_1, \alpha_2) =  e^{ \frac{- 1}{2 i
				\pi} \int_{\mathcal{A}_{\varepsilon}^b}
			\frac{\overset{\swarrow}{\log} \left( 1 - \frac{\alpha_1}{\kappa (k, z)}
				\right)}{z - \alpha_2} \, \mathd z } &\text{ and }&
		\mathfrak{K}_{- -} (\alpha_1, \alpha_2)  =  e^{ \frac{1}{2 i \pi}
			\int_{\mathcal{A}_{\varepsilon}^a} \frac{\overset{\swarrow}{\log}
				\left( 1 - \frac{\alpha_1}{\kappa (k, z)} \right)}{z - \alpha_2} \, \mathd z} .
\end{eqnarray*}}
Going back to $K_{- \circ} (\tmmathbf{\alpha})$, we have
\begin{eqnarray*}
	& K^2_{- \circ} (\tmmathbf{\alpha}) = \frac{\mathfrak{K}_{- \circ}
		(\tmmathbf{\alpha})}{\kappa (k, \alpha_2)} = \frac{\mathfrak{K}_{- +}
		(\alpha_1, \alpha_2)}{\sqrt[\downarrow]{k + \alpha_2}} 
	\frac{\mathfrak{K}_{- -} (\alpha_1, \alpha_2)}{\sqrt[\downarrow]{k -
			\alpha_2}} \, \cdot & 
\end{eqnarray*}
Note that $\sqrt[\downarrow]{k + \alpha_2}$ is a {\tmem{plus function}} in the
$\alpha_2$-plane (branch point at $\alpha_2 = - k$) and $\sqrt[\downarrow]{k -
	\alpha_2}$ is a {\tmem{minus function}} in the $\alpha_2$-plane (branch point
at $\alpha_2 = + k$). Hence the function $\mathfrak{K}_{- +} /
\sqrt[\downarrow]{k + \alpha_2}$ is a {\tmem{plus function}} and the function
$\mathfrak{K}_{- -} / \sqrt[\downarrow]{k - \alpha_2}$ is a {\tmem{minus
		function}}. We can then write $ K_{- \circ} (\tmmathbf{\alpha}) = K_{- -} (\tmmathbf{\alpha}) K_{- +}
(\tmmathbf{\alpha})$, where $K_{- -} (\tmmathbf{\alpha})$ is analytic in $\tmop{LHP}_1 \times
\tmop{LHP}_2$ and $K_{- +} (\tmmathbf{\alpha})$ is analytic in $\tmop{LHP}_1
\times \tmop{UHP}_2$ when both are considered as functions of $\alpha_2$ and
given by
\begin{eqnarray}
& K_{- +} (\tmmathbf{\alpha}) = \left(
\frac{\mathfrak{K}_{- +} (\alpha_1, \alpha_2)}{\sqrt[\downarrow]{k +
		\alpha_2}} \right)^{1 / 2} =
\frac{1}{\sqrt[\downarrow]{\sqrt[\downarrow]{k + \alpha_2}}} \exp \left\{
\frac{- 1}{4 i \pi} \int_{\mathcal{A}_{\varepsilon}^b}
\frac{\overset{\swarrow}{\log} \left( 1 - \frac{\alpha_1}{\kappa (k, z)}
	\right)}{z - \alpha_2} \, \mathd z \right\}, &  \label{app:K-+}
\end{eqnarray}
and
\begin{eqnarray}
& K_{- -} (\tmmathbf{\alpha}) =\left(
\frac{\mathfrak{K}_{- -} (\alpha_1, \alpha_2)}{\sqrt[\downarrow]{k -
		\alpha_2}} \right)^{1 / 2} =
\frac{1}{\sqrt[\downarrow]{\sqrt[\downarrow]{k - \alpha_2}}} \exp \left\{
\frac{1}{4 i \pi} \int_{\mathcal{A}_{\varepsilon}^a}
\frac{\overset{\swarrow}{\log} \left( 1 - \frac{\alpha_1}{\kappa (k, z)}
	\right)}{z - \alpha_2} \, \mathd z \right\}, &  \label{eq:appK--}
\end{eqnarray}
recovering (\ref{eq:integralK-+}) and (\ref{eq:integralK--}). This choice of
realising the square root of the numerator by solely halving the inside of the
exponential ensures that no spurious branch cuts occur. This would have been
the case if instead we chose to take $\sqrt{\phantom{k}}$ or even $\sqrt[\downarrow]{\phantom{k} }$
of the numerator. The second square root of the denominator does not affect
its branch cut structure. These functions are very fast to evaluate since the
integrand now decays like $x^{- 2}$ along $\mathcal{A}_{
	\varepsilon}^a (x)$ as $x \rightarrow \pm \infty$. 


\section{On the application of Liouville's theorem}

\subsection{A useful result}

The following lemma is \COM{establishing} a link between the decay of a
function $\Phi (\alpha)$ and the decay of its respective plus and minus
sum-split parts $\Phi_+ (\alpha)$ and $\Phi_- (\alpha)$.

\begin{lemma}
	\label{lem:Cauchydecay}Let $\Phi (\alpha)$ be a function analytic on some
	strip. And consider its sum-split $\Phi (\alpha) = \Phi_+ (\alpha) + \Phi_-
	(\alpha)$, where $\Phi_+$ and $\Phi_-$ are analytic in the UHP and LHP
	respectively.
	\begin{enumerate}
		\item[a)] If $\Phi (\alpha) =
		\mathcal{O} (1 / | \alpha |^{\COM{\lambda}})$ as $|\alpha| \rightarrow \infty$ within the strip, with $\COM{\lambda} > 1$, then
		$\Phi_{\pm} (\alpha)$ are decaying at least like $1 / | \alpha |$ as $|
		\alpha | \rightarrow \infty$ within their respective half-plane\COM{s}.
		
		\item[b)] If $\Phi (\alpha) =
		\mathcal{O} (1 / | \alpha |)$ as $|\alpha| \rightarrow \infty$ within the strip, then $\Phi_{\pm} (\alpha)$ are decaying at
		least like $\ln | \alpha | / | \alpha |$ as $| \alpha | \rightarrow \infty$ within their respective half-plane\COM{s}.
		
		\item[c)] If $\Phi (\alpha) =
		\mathcal{O} (1 / | \alpha |^{\COM{\lambda}})$ as $|\alpha| \rightarrow \infty$ within the strip, with $0 < \COM{\lambda} < 1$, then
		$\Phi_{\pm} (\alpha)$ are decaying at least like $1/ |
		\alpha |^{\COM{\lambda}}$as $| \alpha | \rightarrow \infty$ within their respective half-plane\COM{s}.
	\end{enumerate}
\end{lemma}

These results are classic. The leading order results (as presented here) can be found for example in \cite{Woolcock1967}, while full asymptotic expansions are given in \cite{McClure1978} and \cite{Wong1980}.

\subsection{For the $\alpha_1$ plane factorisation}\label{app:Liouville1}

Let us show that the top \RED{(resp. bottom)} line of (\ref{eq:EalphaLiouville1}) tends to zero as $|
\alpha_1 | \rightarrow \infty$ within $\tmop{UHP}_1$ \RED{(resp. $\tmop{LHP}_1$)}\COM{, while $\alpha_2 \in \mathcal{A}_2$ is fixed}. First of all, \COM{due to (\ref{eq:exactfunctionG++}),} it is
clear that
\begin{eqnarray*}
	G_{+ +} (\alpha_1, \alpha_2) / K_{- \circ} (a_1, \alpha_2) &
	\overset{\alpha_2 \tmop{fixed}}{\underset{| \alpha_1 | \rightarrow
			\infty}{=}} & \mathcal{O} (1 / | \alpha_1 |) \, .
\end{eqnarray*}
\COM{The condition on the $(x_1\!=\!0,x_2\!>\!0)$ edge implies that for a fixed $x_2>0$, for $x_3=0^+$, we have $\tfrac{\partial u}{\partial x_3}=\mathcal{O}(x_1^{-1/2})$ as $x_1\to 0^+$, while $u=\mathcal{O}((-x_1)^{1/2})$ as $x_1\to 0^-$. Because $F_{++}\propto\mathfrak{F}[\tfrac{\partial u}{\partial x_3}]$ and $G_{-+}=\mathfrak{F}[u_2]$ (see (\ref{eq:G+-definition})), the Abelian theorems \cite{Noble} imply that}
\begin{eqnarray}
F_{+ +} \overset{\alpha_2 \tmop{fixed}}{\underset{| \alpha_1 | \overset{\RED{\tmop{UHP}}}{\rightarrow}
		\infty}{=}} \mathcal{O} (1 / | \alpha_1 |^{1 / 2}) & \RED{\text{ and }} &  \RED{G_{- +} \overset{\alpha_2 \tmop{fixed}}{\underset{| \alpha_1 | \underset{\RED{\tmop{LHP}}}{\rightarrow}
			\infty}{=}} \mathcal{O} (1 / | \alpha_1 |^{3 / 2})}.
\label{appeq:raphB1}
\end{eqnarray}
\COM{For a fixed $x_2<0$ and $x_3=0^+$, the field is well-behaved as $x_1\to 0$, hence $u=\mathcal{O}(1)$. Since $G_{--}=\mathfrak{F}[u_3]$ and $G_{+-}=\mathfrak{F}[u_4]$, the Abelian theorems imply that}
\begin{eqnarray}
G_{+ -} \overset{\alpha_2 \tmop{fixed}}{\underset{| \alpha_1 | \overset{\RED{\tmop{UHP}}}{\rightarrow}
		\infty}{=}} \mathcal{O} (1 / | \alpha_1 |) &\RED{\text{ and }} &  \RED{G_{- -}  \overset{\alpha_2 \tmop{fixed}}{\underset{| \alpha_1 |  \underset{\RED{\tmop{LHP}}}{\rightarrow}
			\infty}{=}} \mathcal{O} (1 / | \alpha_1 |)}.
\label{appeq:raphB2}
\end{eqnarray}
Moreover we have
\begin{eqnarray}
K_{+ \circ} (\alpha_1, \alpha_2) \overset{\alpha_2 \tmop{fixed}}{\underset{|
		\alpha_1 | \rightarrow \infty}{=}} \mathcal{O} (1 / | \alpha_1 |^{1 / 2}) &
\text{ and } & K_{- \circ} (\alpha_1, \alpha_2) \overset{\alpha_2
	\tmop{fixed}}{\underset{| \alpha_1 | \rightarrow \infty}{=}} \mathcal{O} (1
/ | \alpha_1 |^{1 / 2}) \, .
\label{appeq:raphB3}
\end{eqnarray}
Hence, \RED{using (\ref{appeq:raphB1}), (\ref{appeq:raphB2}) and (\ref{appeq:raphB3}),} we know that
\begin{eqnarray*}
	F_{+ +} K_{+ \circ} \overset{\alpha_2 \tmop{fixed}}{\underset{| \alpha_1 |
			\overset{\RED{\tmop{UHP}}}{\rightarrow} \infty}{=}} \mathcal{O} \left(\tfrac{1}{| \alpha_1 |}\right),\,& 
	\frac{G_{+ -}}{K_{- \circ}} \overset{\alpha_2 \tmop{fixed}}{\underset{| \alpha_1 | \overset{\RED{\tmop{UHP}}}{\rightarrow}
			\infty}{=}}\mathcal{O} \left(\tfrac{1}{ | \alpha_1 |^{1 / 2}}\right),\, & \frac{G_{- \circ}}{K_{- \circ}} \overset{\alpha_2 \tmop{fixed}}{\underset{| \alpha_1 | \underset{\RED{\tmop{LHP}}}{\rightarrow}
			\infty}{=}}\mathcal{O} \left(\tfrac{1}{| \alpha_1 |^{1 / 2}}\right).
	\label{eq:decayfixedalpha2g+-}
\end{eqnarray*}
Finally, using the Lemma \ref{lem:Cauchydecay}c) in the $\alpha_1$ plane, we conclude
that we have (at least)
\begin{eqnarray*}
	\left[ \frac{G_{+ -}}{K_{- \circ}} \right]_{+ \circ} \overset{\alpha_2
		\tmop{fixed}}{\underset{| \alpha_1 | \overset{\RED{\tmop{UHP}}}{\rightarrow} \infty}{=}}  \mathcal{O}
	\left( \frac{1}{| \alpha_1 |^{1 / 2}} \right) & \text{ and } & \left[ \frac{G_{+ -}}{K_{- \circ}} \right]_{- \circ} \overset{\alpha_2
		\tmop{fixed}}{\underset{| \alpha_1 | \underset{\RED{\tmop{LHP}}}{\rightarrow}\infty}{=}}  \mathcal{O}
	\left( \frac{1}{| \alpha_1 |^{1 / 2}} \right) .
\end{eqnarray*}
This shows that \RED{the terms of the} top \RED{(resp. bottom)} line of (\ref{eq:EalphaLiouville1}) go to
zero as $| \alpha_1 | \rightarrow \infty$ within $\tmop{UHP}_1$ \RED{(resp. $\tmop{LHP}_1$)}. Hence, Liouville's theorem can safely be applied.

\subsection{For the $\alpha_2$ plane factorisation}\label{app:Liouville2}

\RED{Here we wish to show that $E_{+2}$ tends to zero when $\alpha_1$ is fixed in $\text{UHP}_1$ and $|\alpha_2|\to \infty$.}

\subsubsection{\RED{The terms without brackets}}\label{app:Liouville2-1}

Let us show that \RED{the terms without brackets in} the top line\footnote{\RED{We can show in a very similarly fashion (omitted for brevity) that the terms without brackets in the bottom line of (\ref{eq:Eplus2alphaLiouville}) do also tend to zero for fixed $\alpha_1 \in \text{UHP}_1$ and $|\alpha_2| \to \infty$ within $\text{LHP}_2$.}} of (\ref{eq:Eplus2alphaLiouville}) tend to zero as
$\alpha_1$ is fixed \RED{in $\text{UHP}_1$} and $| \alpha_2 | \rightarrow \infty$ within
\RED{$\text{UHP}_2$}. First of all, using the \COM{expression (\ref{eq:exactfunctionG++}) for} $G_{+ +}$, it is
straightforward to see that
\begin{eqnarray}
\frac{G_{+ +} (\tmmathbf{\alpha})}{K_{- -} (a_1, a_2) K_{+ -} (\alpha_1,
	a_2)} & \overset{\alpha_1 \tmop{fixed}}{\underset{| \alpha_2 | \rightarrow
		\infty}{=}} & \mathcal{O} (1 / | \alpha_2 |) . 
\end{eqnarray}
Moreover, using the definition of $K_{+ +}$ and $K_{- +}$, we can see that
\begin{eqnarray*}
	K_{+ +} (\alpha_1, \alpha_2) \overset{\alpha_1 \tmop{fixed}}{\underset{|
			\alpha_2 | \rightarrow \infty}{=}} \mathcal{O} (1 / | \alpha_2 |^{1 / 4}) &
	\text{  and } & K_{- +} (\RED{a_1}, \alpha_2) \RED{\underset{| \alpha_2 | \rightarrow \infty}{=}} \mathcal{O} (1
	/ | \alpha_2 |^{1 / 4}) .
\end{eqnarray*}
\COM{The $(x_1>0,x_2=0)$ edge condition implies that for $x_3=0^+$ and a fixed $x_1>0$, 
	$\tfrac{\partial u}{\partial x_3}=\mathcal{O}(x_2^{-1/2})$ as $x_2\to 0^+$, which,  by the Abelian theorems, requires that $F_{+ +} (\alpha_1, \alpha_2)=\mathcal{O} (1 / | \alpha_2 |^{1 / 2})$ for $\alpha_1 \in \text{UHP}_1$ as $|\alpha_2|\to\infty$ within $\text{UHP}_2$. This leads to}
\begin{eqnarray}
F_{+ +} (\tmmathbf{\alpha}) K_{+ +} (\tmmathbf{\alpha}) K_{- +} (a_1,
\alpha_2) & \overset{\alpha_1 \tmop{fixed}}{\underset{| \alpha_2 |
		\rightarrow \infty}{=}} & \mathcal{O} (1 / | \alpha_2 |) . \label{eq:B5secondrevraph}
\end{eqnarray}

\COM{\subsubsection{The terms with brackets}\label{app:Liouville2-2}
	Because of Lemma \ref{lem:Cauchydecay}, in order to prove that the last term on the top (resp. bottom) line of (\ref{eq:Eplus2alphaLiouville}) tends to zero as $|\alpha_2|\to \infty$ within $\text{UHP}_2$ (resp. $\text{LHP}_2$), it is sufficient to show that $\tfrac{K_{-+}(a_1,\alpha_2)}{K_{+-}}[\tfrac{G_{+-}}{K_{-\circ}}]_{+ \circ}$ tends to zero as a power of $|\alpha_2|$ as $|\alpha_2|\rightarrow\infty$ while on $\mathcal{A}_2$. In order to show this\footnote{\COM{Thank you to the anonymous reviewer for this
			suggestion!}}, rewrite (\ref{eq:secondsplitinter1}), which is valid for $\alpha_2 \in \mathcal{A}_2$,
	as
	\begin{align*}
	\frac{K_{- +} (a_1, \alpha_2)}{K_{+ -}} \left[ \frac{G_{+
			-}}{K_{- \circ}} \right]_{+ \circ}  =  \underbrace{F_{+ +} K_{+ +} K_{-
			+} (a_1, \alpha_2)}_{\overset{\alpha_1 \tmop{fixed}}{\underset{| \alpha_2 |
				\rightarrow \infty}{=}} \mathcal{O} (1 / | \alpha_2 |)} -
	\underbrace{\frac{G_{+ +}}{K_{- -} (a_1, \alpha_2) K_{+
				-}}}_{\overset{\alpha_1 \tmop{fixed}}{\underset{| \alpha_2 | \rightarrow
				\infty}{=}} \mathcal{O} (1 / | \alpha_2 |^{1 / 2})} \cdot
	\end{align*}
	The first estimate on the RHS is already given in the paper in (\ref{eq:B5secondrevraph}), whilst
	the second comes naturally from the asymptotic behaviours
	\begin{align*}
	G_{+ +} \overset{\alpha_1 \tmop{fixed}}{\underset{| \alpha_2 | \rightarrow
			\infty}{=}} \mathcal{O} (1 / | \alpha_2 |), \quad K_{- -} \overset{\alpha_1
		\tmop{fixed}}{\underset{| \alpha_2 | \rightarrow \infty}{=}} \mathcal{O} (1
	/ | \alpha_2 |^{1 / 4}), \quad K_{+ -} \overset{\alpha_1
		\tmop{fixed}}{\underset{| \alpha_2 | \rightarrow \infty}{=}} \mathcal{O} (1
	/ | \alpha_2 |^{1 / 4}) .
	\end{align*}
	The first of these results is obvious given the exact expression (\ref{eq:exactfunctionG++}) for $G_{+
		+}$, whilst the second and third come directly from the integral
	representations (\ref{eq:integralK--}) and (\ref{eq:integralK+-}). Hence, we have proved that the sufficient
	condition is satisfied and that, consequently, Liouville's theorem can safely be applied to obtain $E_{+2}\equiv0$.}

{
\bibliographystyle{plain}
\bibliography{library}
}

\end{document}